\documentclass[11pt]{article}

\usepackage{amssymb}
\usepackage{latexsym}
\usepackage[mathscr]{euscript}
\usepackage{graphicx}

\usepackage{amscd}
\usepackage{amsmath}
\usepackage{amssymb}
\usepackage{fullpage}

\usepackage{amssymb}
\usepackage{latexsym}
\usepackage[mathscr]{euscript}
\usepackage{pb-diagram}
\usepackage{latexsym}
\usepackage{amsfonts}
\usepackage{epsfig}

\newcommand\tover[1]{\column{\otimes}{#1}}

\newcommand\pr{^{\prime}}

\newcommand\commentout[1]{\marginpar{\tiny $\backslash$commentout}}

\newcommand\qed{\hfill$\square$}
\def\column#1#2{\mathrel{\mathop{#1}\limits_{#2}}}

\def\compcirc {\mbox{\hspace{.05cm}}\raisebox{.04cm}{\tiny  {$\circ$ }}}

\newtheorem{Lemma}{Lemma}[section]
\newtheorem{Theorem}{Theorem}[section]
\newtheorem{Proposition}{Proposition}[section]
\newtheorem{Definition}{Definition}[section]
\newtheorem{Corollary}{Corollary}[section]

\newenvironment{Proof}{\par\noindent\textbf{Proof:}}
{\qed}

\usepackage{graphics}
\title{Project Description}

\title{A geometric reformulation of the representational assembly conjecture for the $K$-theory of fields with pro-$l$ absolute Galois group}
\author{Gunnar Carlsson \footnote{Research supported in part by NSF DMS-0406992}  \\ Department of Mathematics, Stanford University  \\Stanford, California 94305}
\begin{document}
\maketitle
\section{Introduction} 
In recent years there have been striking developments in the study of the algebraic $K$-theory of fields (see \cite{voevodsky}, \cite{friedlander}, \cite{survey}).  These developments have resulted in the identification of the $E_2$-term of a spectral sequence converging to the algebraic $K$-theory of fields (rather a mod-$p$ version of it) in terms of Galois cohomology, and have resolved long standing conjectures of Milnor and Bloch-Kato in the area.    In particular the $E_2$-term depends only on the absolute Galois  group of the field.  These results raise a pair of questions.   
\begin{itemize}
\item{For a field $F$, is there a construction depending only on the absolute Galois  group $G_F$  of $F$, which reconstructs the full $K$-theory spectrum of $F$, rather than an $E_2$-term of a spectral sequence?  Can such a construction be made functorial for inclusions of subfields and the corresponding inclusions of absolute Galois  groups? Is the construction readily describable in terms of well understood constructions on the absolute Galois group?  }
\item{The conjectures of Milnor and Bloch-Kato are cohomological statements about the absolute Galois  group, specifically that the cohomology is generated in degree one and  that the relations are quadratic.  Are there properties more closely connected to the structure of the absolute Galois  group as a group which are relevant?  }
\end{itemize} 

A conjectural answer to the first question in the case of geometric fields (i.e. fields containing an algebraically  closed subfield $k$)  was proposed in \cite{repassembly}, and proved in the case of abelian absolute Galois  groups.  This conjecture is formulated as follows.  We assume we are given a field $F$ containing an algebraically closed subfield $k \subseteq F$, let $\overline{F} $ denote its algebraic closure, and let $G_F$ denote its absolute Galois group $Gal(\overline{F} /F)$, in general a profinite group.  One can construct a version of the equivariant $K$-theory spectrum,  $K^G$,  which applies to profinite groups, and observe that $K^G(\overline{F}) \cong K(F)$.  On may also consider the natural map $K^G(k) \rightarrow K^G(\overline{F})$ of spectra.  This map is far from being an equivalence since, for example, $\pi _0 K^G(k)$ is the free abelian group on the continuous $k$-linear representations of $G_F$, where $k$ is given the discrete topology, and on the other hand 
$$ \pi _0 K^G(\overline{F}) \cong \pi _0 K(F) \cong \Bbb{Z}
$$
However, the spectra $K^G(k)$ and $K^G(\overline{F})$ both carry multiplicative structures, which makes $K^G(\overline{F})$ into a module over the commutative ring spectrum $K^G(k)$.  See \cite{ekmm} and \cite{symmetric} for a discussion of these multiplicative theories.   In \cite{completion}, a derived version of the completion construction for commutative rings was developed, which associates to a homomorphism of commutative ring spectra $f:A \rightarrow B$ and an $A$-module $M$ a completion of $M$ along the homomorphism $f$, denoted by $M^{\wedge}_f$, which agrees with the standard completion construction for finitely generated modules over Noetherian commutative rings.   We will denote by $\varepsilon$ the homomorphism $K^G(k) \rightarrow \Bbb{H}_l$ which is induced by the functor which assigns to every representation of $G_F$ its dimension mod $l$, where $\Bbb{H}_l$ is the mod-$l$ Eilenberg-MacLane spectrum, regarded as the result of an infinite loop space machine applied to the symmetric monoidal category with object set equal to $\Bbb{Z}/l\Bbb{Z}$, and with only identity morphisms.  At this point, we may construct derived completions along the homomorphism $\varepsilon$, to obtain a map 
$$\alpha _F : K^G(k)^{\wedge}_{\varepsilon} \rightarrow K^G(\overline{F})^{\wedge}_{\varepsilon}$$
called the {\em representational assembly}, and the conjecture is that this map is an equivalence of spectra. 

The difficulty in approaching this conjecture is that we do not have simple direct methods for computing with the completions.  We will instead approach the problem indirectly as follows.  We will construct an affine scheme  ${\cal E}{G_F}$, with associated ring equipped with the discrete topology, and equipped with a continuous  action by the topological group $G_F$.  Equivariant $K$-theory can be suitably defined for such actions.  The scheme  ${\cal E}{G_F}$ has the following properties.  

\begin{enumerate}
\item{The maps $$K^{G_F}(Spec(k))^{\wedge}_{\varepsilon} \rightarrow K^{G_F}({\cal E}{G_F})^{\wedge} _{\varepsilon}$$ and 
$$K^{G_F}(Spec(\overline{F}))^{\wedge}_{\varepsilon} \rightarrow K^{G_F}(Spec(\overline{F}) \times {\cal E}{G_F})^{\wedge}_{\varepsilon}$$ are equivalences.}
\item{There is a well-defined notion orbit scheme for ${\cal E}F$, and there are equivalences $$K({\cal E}{G_F}/G_F) \rightarrow K^{G_F}({\cal E}{G_F})$$ and 
$$K(Spec(\overline{F}) \column{\times}{G_F} {\cal E}{G_F}) \rightarrow K^{G_F}(Spec(\overline{F}) \times {\cal E}{G_F})$$  This follows from a freeness property of ${\cal E}{G_F}$. }
\item{ There is a natural equivalence $K^{G_F}({{\cal E}G_F)} ^{\wedge}_{\epsilon} \cong K^{G_F}({\cal E}G_F) ^{\wedge}_{l}$, where $( - )^{\wedge}_l$ denotes Bousfied-Kan $l$-adic completion.    }
\end{enumerate} 
These properties give rise to the following diagram.

$$
\begin{diagram} \node{K^{G_F}(Spec(k))^{\wedge}_{\varepsilon}} \arrow{s,t}{\alpha _F} \arrow{e} \node{K^{G_F}({\cal E}{G_F})^{\wedge}_{\varepsilon}} \arrow{s} 
\node{K({\cal E}{G_F} /G_F) ^{\wedge}_{l}} \arrow{w} \arrow{s,b}{\eta} \\
\node{K^{G_F}(Spec(\overline{F})) ^{\wedge}_{\varepsilon}} \arrow{e}
 \node{K^{G_F}(Spec(\overline{F})
  \times {\cal E}{G_F})^{\wedge} _{\varepsilon}} \node{K(Spec(\overline{F})\column{\times}{G_F} {\cal E}{G_F})^{\wedge}_l} \arrow{w}
\end{diagram} 
$$
The goal is to prove that $\alpha _F$ is an equivalence.  Properties (1) and (3) above show that all the horizontal maps are equivalences, and it therefore follows that it will suffice to prove that the right hand vertical map $\eta$ is an equivalence.  Proving that the map $\eta$ is an equivalence can be thought of as a parametrized rigidity theorem.   We will use this reduction in the paper \cite{royg} to prove the conjecture.  Some comments on the construction are in order. 
\begin{itemize}
\item{${\cal E}G$ is not functorial for homomorphisms of groups, and it depends on a number of choices.   Moreover, ${\cal E}G_F$ is not even $\Bbb{A}^1$-contractible, as its algebraic $K$-theory is not isomorphic to the algebraic $K$-theory of the base field $k$.    What can be proved is that it is potentially $l$-adically contractible, in that $K({\cal E}G_F)$ has $l$-adic completion equivalent to the $l$-adic completion of $K(k)$.   We believe that there should ultimately be canonical constructions for the universal space of a profinite group, and that our construction is a particular model which is $l$-adically equivalent to it.  }
\item{The equivariant $K$-theory of ${\cal E}G_F$ is interesting in that the natural action of the representation ring of $G_F$, $R[G_F]$, on $K_*({\cal E}G_F)$ has the property that the augmentation ring of $R[G_F]$ acts trivially on it. This fact is crucial in the proof of  Property (3) above.  }
\item{ The construction of ${\cal E}{G_F}$ relies on a structural property shared by all absolute  Galois groups of geometric fields, that of {\em total torsion freeness}, which asserts that the abelianization of any closed subgroup of $G_F$ is torsion free.  This is a surprisingly restrictive property for groups, and permits the construction of free actions on varieties of the form $\Bbb{G}_m^n$, which turn out to be extremely useful.  
}
\item{${\cal E}G_F $ is of the form $Spec(A)$, where $A = \bigcup_n k[t_{s}^{\pm \frac{1}{l^n}}; s \in S]$, where $S$ is some indexing set.  }
\item{There has been a great deal of work on reconstructing the arithmetic of a fields from its Galois group, see \cite{pop} and \cite{bogomolov}, for example.  In the abstract, this asserts that one is able to reconstruct the $K$-theory spectrum of the field from its absolute Galois group.  The point of the present work is to identify the $K$-theory spectrum in terms of well understood constructions on the absolute Galois group, with the expectation that such understanding will make applying $K$-theory to arithmetic problems will be made easier.  }

\item{There is another formulation of the main conjecture, which would replace the ring spectra which occur by their counterparts within equivariant spectra,  and would therefore naturally replace all homotopy groups which occur by their corresponding Mackey functors, and further replace the representation rings by their Green functor counterparts (see \cite{Bouc}).  The two conjectures are not a priori equivalent, since the relationship between the completion in the category of rings and the completion in the Green functor context is not at all clear.  Indeed, our results are proved for fields whose absolute Galois group is a pro-$l$ group, although we expect to deal with the general case in future work.  Moreover, the evaluation of the Mackey functor completion would be expected in general to be very difficult, since very simple questions about, for example, derived functors of tensor products are quite difficult in the Green functor situation \cite{ventura}.  C. Barwick has developed methods for proving the Mackey/Green functor version of the conjecture. Presumably one will find that our methods will permit the computation of the derived completions in the Mackey/Green functor version, in terms of derived completions of ordinary rings.  }
\end{itemize}

The paper is organized as follows. Section 2 records preliminary material, including a discussion of actions of profinite groups on affine schemes which will be important for the construction of ${\cal E}{G_F}$.  Section 3 defines totally torsion free  profinite groups, and proves that absolute  Galois groups all enjoy this property.  Section 4 develops a method for constructing profinite group actions on certain (non-Noetherian) affine schemes.  The schemes in question are the inverse limits of systems of the form 
$$  \cdots \stackrel{\times l}{\longrightarrow} \Bbb{G}_m^n  \stackrel{\times l}{\longrightarrow} \Bbb{G}_m^n  \stackrel{\times l}{\longrightarrow} \Bbb{G}_m^n  \stackrel{\times l}{\longrightarrow} \Bbb{G}_m^n  
$$
and ultimately ${\cal E}{G_F}$ will be constructed as an infinite product of such schemes. We are able to construct actions on such varieties from affine representations over $\Bbb{Z}_l$  of $G_F$.   Section 5 actually constructs ${\cal E}{G_F}$, using the totally torsion free property to construct sufficiently many such affine representations.  Section 6 is a technical section, which studies the representation rings of profinite $l$-groups which are an extension of a torsion free, topologically finitely generated, abelian profinite $l$-group by a finite $l$-subgroup of the symmetric group.   Sections 7 and 8 now use this preparatory work to prove the main result.

The author wishes to thank B. Conrad, P. Diaconis, D. Dugger, L. Hesselholt, T. Lawson, G. Lyo,  I. Madsen, D. Ramras, R. Vakil, and K. Wickelgren for numerous helpful conversations on this work.  

\section{Preliminaries }
\subsection{ General}

We will assume that the reader is familiar with the theory of ring spectra and modules over them, as well as the notion of a commutative ring spectrum. We will use the version of this theory which uses symmetric spectra, which is presented in \cite{schwede}.  In particular, we will assume the reader is familiar with the standard choices of model structures on this category, and the associated notions of cofibrant objects, as well as $Hom_A(M,N)$, $M \column{\wedge}{A} N$ and universal coefficient and K\"{u}nneth spectral sequences converging to them.    We will also assume familiarity with homotopy limits and colimits in this category.  Throughout the paper, a spectrum will mean a symmetric spectrum, a ring spectrum will mean a symmetric ring spectrum, and a commutative ring spectrum will mean a commutative symmetric ring spectrum.  


\subsection{$K$-theory} \label{ktheory}

There are numerous versions of the construction of the $K$-theory spectra associated with  various kinds of categorical input data.  We will require two of them.  The first is the construction to A. Elmendorf and M. Mandell \cite{elmendorf}.  It takes as input small permutative categories, and produces from that construction a spectrum.  It advantage is that the category theoretic information required to impose a commutative ring spectrum structure on the resulting spectrum, and to produce module spectra over such spectra,  is codified in \cite{elmendorf}.  This kind of multiplicative structure is at the core of the results of this paper.  The second construction (due to Waldhausen \cite{waldhausen}) takes as input data a category with cofibrations and weak equivalences, and produces a $K$-theory spectrum from this data. It has the advantage that it has the excision properties proved to exist by Quillen (\cite{quillen}), such as localization and devissage.   Each Waldhausen category can be regarded as a small permutative category, to which the Elmendorff-Mandell construction applies,   and there is a sequence of natural equivalences of spectra relating the two constructions  which are described in \cite{repassembly}, section 4.  We will denote the Elmendorff-Mandell construction  by $\frak{K}$ and the Waldhausen construction by $K$.  We also define a weakening of the notion of a ring spectrum as follows. 

\begin{Definition} By a {\em weak ring spectrum}, we mean a spectrum $R$ equipped with a map $R \wedge R \rightarrow R$ which is homotopy associative in the evident sense.  By a {\em module} over a weak ring spectrum $R$, we will mean a spectrum $M$ equipped with a map $R \wedge M \rightarrow M$, where the diagrams involving the structure maps commute up to homotopy.  The Waldhausen construction readily produces weak ring structures from tensor products on modules, as well as modules from appropriate notions of tensor products. 
\end{Definition}

{\bf Remark:} The reason we need both constructions is the following.  In order to define our derived completion construction, we will need to have the precise commutative ring structure on our spectra afforded by the Elmendorf-Mandell model.  On the other hand, for computation behavior of homotopy groups, it is particularly important to realize constructions such as homotopy fibers as the $K$-theory of categories themselves, and this can be carried out effectively using the Waldhausen construction. We note that the structure of a weak ring spectrum is sufficient to compute behavior on homotopy groups.   The natural equivalences allow us to go back and forth between models to obtain the required information.  

\begin{Proposition}\label{requirements} Let $\underline{C}$ an abelian category, equipped with a strictly associative symmetric monoidal product, which turns $\underline{C}$ into a bipermutative category in the sense of \cite{elmendorf}.  Then the following hold. 
\begin{enumerate}
\item{$\frak{K}(\underline{C})$ is naturally a commutative ring spectrum in the sense of \cite{elmendorf}.  }
\item{Let $\underline{D} \subseteq \underline{C}$ be a Serre subcategory, with the property that it is closed under tensor product with any object in $\underline{C}$.  Then the quotient category $\underline{C}/\underline{D}$ is also bipermutative, and $\frak{K}(\underline{C}/\underline{D})$ is also a commutative ring spectrum.  }
\item{The natural equivalence between $\frak{K}$ and $K$ respects the weak ring spectrum structures, where the weak ring spectrum structure on $\frak{K}$ is simply the forgetful one obtained from the full ring spectrum structure. }
\item{Given $\underline{C}$ and $\underline{D}$ as above, the homotopy fiber of the evident map $\frak{K}(C) \rightarrow \frak{K}(\underline{C}/\underline{D})$ is a $\frak{K}(\underline{C})$-module.  As a weak $K(\underline{C})$-module, it is equivalent to $K(\underline{D})$.}

\end{enumerate} 
\end{Proposition}

\subsection{Equivariant algebraic $K$-theory} \label{equivariant}

We will fix an algebraically  closed ground field $k$ throughout this paper.    In this section we will define the equivariant algebraic-$K$-theory of finite group actions on certain commutative rings. We will be addressing the fact that in order to build the schemes we require, we will have to leave the Noetherian setting.  The suitable setting will be that of coherent commutative rings, for which the category of finitely presented modules is an abelian category. We will let $R[G]$ denote the representation ring of the group $G$ in an algebraically closed field $k$, which is a fixed ground field throughout the paper.  We also let $I_G$ denote the augmentation ideal, i.e. the kernel of the natural augmentation $R[G] \rightarrow \Bbb{Z}$.    We recall from \cite{glaz} the following definition. 

\begin{Definition} A commutative  ring $A$ is said to be {\em coherent} if for every homomorphism $f : F_0 \longrightarrow F_1$ of finitely generated free $A$-modules, the kernel of $f$ is a finitely generated $A$-module.  
\end{Definition}

A standard fact about coherent rings is the following. 
\begin{Proposition} The category of finitely presented modules over a coherent ring $A$ is an abelian category.  
\end{Proposition}
The following result is an extension of    Proposition 4.5 of \cite{repassembly}.  
\begin{Proposition}  \label{coherent}Let $D$ be a directed set, and let $F$ be a functor from the associated category to the category of rings.  Suppose further that $F(d \leq d \pr)$ is faithfully flat for all $d \leq d \pr$, and that all the rings $F(d)$ are coherent.  Then the ring $\column{colim}{D} F$ is coherent.   
\end{Proposition} 
\begin{Proof} The proof is identical to the proof of Proposition 4.5 of \cite{repassembly}, and we omit it. 
\end{Proof}

\begin{Definition}  We say a ring is {\em Noetherian approximable} if it can be described as a colimit over a directed set of Noetherian rings, with the transition maps in the colimit being faithfully flat.  
\end{Definition} 
Given a coherent commutative ring $A$, it is not known whether or not a polynomial ring $A[x_1, \ldots, x_n]$ is also coherent.  However, this does hold for Noetherian approximable rings $A$. 
\begin{Proposition}\label{noetherianapproximable}  Let $A$ be a Noetherian approximable commutative ring.  Then the ring $B = A[x_s;s \in S]$, where $S$ is some indexing set, is Noetherian approximable and therefore coherent. 
\end{Proposition}
\begin{Proof}  $B$ can be described as the colimit of rings of the form $R[x_1, \ldots , x_t]$, where $R$ is Noetherian, and where the transition maps in the colimit are faithfully flat.  $R[x_1, \ldots , x_n ]$ is Noetherian, and the maps in the system are clearly faithfully flat, hence the result. 
\end{Proof}

We let $A$ be a coherent  commutative ring, equipped with an action by a finite group $G$.  

\begin{Definition}  By a {\em $G$-twisted $A$-module}, we will mean a finitely presented $A$-module $M$ equipped with a $k$-linear $G$-action so that $g \cdot (\alpha \cdot m ) = \alpha ^g  \cdot (g \cdot m)$ for all $g \in G, \alpha \in A$, and $m \in M$.   The $G$-twisted $A$-modules form a category whose morphisms are the equivariant $A$-linear maps, which we denote $\mbox{\em Mod}^G(A)$.  
\end{Definition} 

\begin{Proposition}  Given $A$ and $G$ as above, the category $\mbox{\em Mod}^G(A)$ is an abelian category.
\end{Proposition} 
\begin{Proof}  Clear.  
\end{Proof} 

Although we will primarily be considering affine schemes, we will use ring and scheme terminology interchangeably.




\begin{Definition}  We will define the $G$-equivariant $K$-theory on the category of coherent commutative rings  to be the algebraic $K$-theory of the full subcategory of $\mbox{Mod}^G(A)$ on the finitely presented $G$-twisted modules, and denote it by $K^G(A)$.  $K^G$ is functorial for flat homomorphisms of coherent commutative rings with $G$ action.  In the case of Noetherian rings, this is what would usually be referred to as $G$-theory or $K^{\prime}$-theory.   
\end{Definition}

  In the case of a point, i.e. $A = k $, we have that $K^G(k)$ is simply the category of finite dimensional $k$-linear representations of $G$.  $\pi _0 (K^G(k))$ is then the $k$-linear representation ring $R[G]$.  In particular, if  either (1) $Char(k) = 0$ or (2) $Char(k)$ is relatively prime to $\# (G)$, then $\pi _0K^G(k)$ is isomorphic to the usual complex representation ring or character ring  of $G$.  $K^G(k)$ is a commutative  ring spectrum.     It is also equipped with an augmentation $\epsilon : K^G(k) \rightarrow K(k)$, which forgets the $G$ action, and for which applying $\pi _0$ induces the usual augmentation on representation rings, which we also refer to as $\epsilon$. 
 
 
 We need a brief discussion of orbit schemes.  From \cite{exposefive}, it is known that for finite group actions on affine schemes, there is always a categorical orbit scheme which is the spectrum of the fixed point subring of the action.  
 
  
 \begin{Proposition}  \label{free} Suppose that the action of $G$ on an affine scheme $X = Spec(A) $ is free in the sense that it is free on the set  of all points,  where a point denotes a morphism $Spec(\Omega) \rightarrow X$  for some algebraically closed  field $\Omega$ containing $k$.  Then there is a natural equivalence $K^G(A) \cong K(A^G)$.   It is induced by the composite
 $$ K(A^G) \rightarrow K^G(A^G) \rightarrow K^G(A)
 $$
 Equivalently, for free $G$-actions we have an equivalence $K(S/G) \rightarrow K^G(X)$.  

 \end{Proposition}

 \begin{Proof}  The equivalence of categories from which this equivalence results is proved in \cite{grothendieck}.   \end{Proof} 
 
  We call group actions satisfying the above conditions {\em etale}. 
 The following is well known. 
 
 \begin{Proposition}  The category of affine schemes admits arbitrary inverse  limits.  The ring of all global sections of an inverse limit is the corresponding colimit of the rings of global sections of the component schemes.  
 \end{Proposition}
 
 \begin{Proposition} \label{inverseetale}  Let $G$ be a finite group, and let $\underline{X}$ denote a diagram of affine schemes with $G$-action.  If all the $G$ actions composing the diagram are etale, then so is the $G$-action on $lim \underline{X}$.   
 \end{Proposition} 
 \begin{Proof}  For any algebraically closed field $\Omega$, the set of maps $Spec(\Omega) \rightarrow lim \underline{X}$ is the inverse limit of the sets of maps from $Spec( \Omega)$ to the component schemes.  An inverse limit of free $G$-sets over a directed set is always free.  
 \end{Proof} 
 
 
 
 Given an action of a finite group $G$ on a scheme $X$, and a $G$-invariant subscheme $Y \subseteq X$, we define the category $\mbox{Nil}^G(X,Y)$ to be the category of $G$-twisted $X$-modules whose support lies entirely in $Y$.  
 
 \subsection{Profinite group actions on commutative rings}

We begin with the definition of  profinite groups. The basic material on profinite groups is contained in \cite{ribes}.
\begin{Definition}
A topological group $G$ is a {\em profinite group} if it is Hausdorff, compact, and totally disconnected.  It follows from this definition that $G$ is the inverse limit of the inverse system of finite groups $\{ G/N \}_N$, where the indexing ranges over the closed normal subgroups of $G$ of finite index.  
\end{Definition}

\begin{Proposition} \label{opensubgroup}
Let $G$ be a profinite group.  Then a subgroup of $G$ is open only if it is closed and of finite index. 
\end{Proposition}
We next consider continuous actions of a profinite group $G$ on discrete sets $X$. 

\begin{Proposition} \label{criterion}  Let $X$ be a discrete set, and suppose that we are given an action of a profinite group $G$ on $X$.  Then the action is continuous if and only if the stabilizers $G_x$ are closed subgroups of finite index of $G$ for all $x \in X$.  
\end{Proposition} 
\begin{Proof}  From the discreteness of $X$, it follows that $\alpha ^{-1}(x)$ is an open set in $G \times X$ for all $x \in X$, where $\alpha: G \times X \rightarrow X$ is the action map.  A subset of $G \times X$ is  open if and only if it is of the form 
$$ \coprod _{x \in X} U_x \times \{ x \}
$$
where each $U_x$ is an open set in $G$.  On the other hand, 
$$\alpha ^{-1}(x)   =  \coprod _{\gamma \in \Gamma}  G_x \gamma^{-1} \times \gamma \cdot x
$$
It is now clear that this set is open if and only if $G_x$ is an open subgroup of $G$, i.e. that it is closed and of finite index by Proposition \ref{opensubgroup}. 
\end{Proof} 

We also have the following result on transitive $G$-spaces. 

\begin{Proposition} \label{classify}
Let $X$ be a topological space with a  continuous transitive action of a profinite group $G$, and let $G_x$ denote the stabilizer of a point $x \in X$.   Then if $X$ is Hausdorff, we have that $X$ is $G$-homeomorphic to the quotient space $G/G_x$, and $G_x$ is a closed subgroup of $G$. \end{Proposition} 
\begin{Proof} It is immediate that $X$ is compact, since it is the surjective image of the compact space $G$. The natural map $G/G_x \rightarrow X$ is a bijective continuous map from one compact Hausdorff space to another, and is therefore a homeomorphism.   It is direct that $G/G_x$ is Hausdorff if and only if $G_x$ is closed.  
\end{Proof}

\begin{Definition}  Let $G$ be a profinite group.  By a {\em continuous affine $G$-scheme $A$}, we will mean a continuous action of $G$ on a commutative ring $A$ equipped with the discrete topology.  
\end{Definition} 

{\bf Example:} Let $G$ and $G_{sep}$  be the  Galois groups of $\overline{k}$, the algebraic closure of a field  $k$, and $k_{sep}$,  its absolute closure.  Then the actions of $G$ and $G_{sep}$ on $Spec(\overline{k})$ and $Spec(k_{sep})$ are continuous affine $G$ (respectively $G_{sep}$)-schemes.  

{\bf Example:} Let $G$ denote the $l$-adic integers, and let $A = \bigcup _n k[t^{\pm\frac{1}{l^n}}]$.  The group of automorphisms of $A$ leaving $k[t^{\pm 1}]$ elementwise fixed is isomorphic to $\Bbb{Z}_l$, and the action of $\Bbb{Z}_l$ on $A$ is  a  continuous affine $\Bbb{Z}_l$-scheme.  

Let $G$ be a profinite group, and let $A$ be a continuous  affine $G$-scheme.      Let $B$ denote the fixed point subring $A^G$.  For every closed normal subgroup $N$  of finite index in $G$, we may consider the scheme $Spec(A^N)$. These schemes fit together into an inverse system of affine schemes 
$\{ Spec(A^N) \}_N$, and we may regard the inverse system as an inverse system of sets, by simply considering the points of the schemes.  We also have the group actions 
$$G/N \times Spec(A^N ) \longrightarrow Spec(A^N )
$$
which fit together into an action of a pro-group on a pro-set.   For every $\frak{p} \in Spec(B)$, we may consider the fiber over $\frak{p}$ in 
$\column{lim}{\longleftarrow} Spec(A^N)$, denoted by $F_{\frak{p}}$, and find that we obtain a continuous action of $G$ on $F_{\frak{p}}$.  

\begin{Proposition} Let $A$ be a continuous affine $G$-scheme, with $B = A^G$.  Then the map 
$$ \pi: Spec(A) \rightarrow Spec(B)
$$
is surjective, and all the fibers $\pi ^{-1}(\frak{p})$ are transitive $G$-sets.  Moreover, $\pi ^{-1}(\frak{p})$ is homeomorphic to a $G$-space of the form $G/G^{\prime}$, where $G^{\prime}$ is a closed subgroup of $G$.  
\end{Proposition} 
\begin{Proof}  For every closed subgroup of finite index $K \subseteq G$, we have the subring $A^K$ of elements of $A$ fixed by $K$.  The ring $A$ is equal to the union $\bigcup A^K$ as $K$ varies over smaller and smaller closed subgroups of finite index in $G$, because the stabilizer of every element of $A$ is a closed subgroup of finite index in $G$ by Proposition \ref{criterion}.  We will now show that the natural map 
$$  \alpha: Spec(A) \longrightarrow \mbox{ }\column{lim}{\longleftarrow} Spec(A^K)
$$
is a bijection.  To show that $\alpha$ is injective, consider two prime ideals $\frak{p}$ and $\frak{q}$ in $A$.  If $\frak{p} \neq \frak{q}$, then there is an element $a \in A$ with $a$ an element of one of $\frak{p}$ and $\frak{q}$ and not in the other, say $a \in \frak{p}$ and $a \notin \frak{q}$.  But we have seen that $a \in A^K$ for some closed subgroup $K$ of finite index in $G$.  It follows that the images of $\frak{p}$ and $\frak{q}$ in $ Spec(A^K)$ are distinct, which gives the injectivity of $\alpha$.   For surjectivity, we first note that elements $\column{lim}{\longleftarrow} Spec(A^K)$ correspond exactly to families $\{ \frak{p} _K \}_K$, where the system varies over the closed subgroups of finite index, where each $\frak{p}_K$ is a prime ideal in $A^K$, and where we have the consistency condition $\pi ^K_{K^{\prime}} (\frak{p}_K) = \frak{p}_{K^{\prime}}$ whenever $K \subseteq K^{\prime}$.  Here $\pi ^K_{K^{\prime}}$ denotes the restriction map $Spec(A^K) \rightarrow Spec(A^{K^{\prime}})$.  Given a family ${\cal F} = \{ \frak{p} _k \} _K$ as above,  we define a subset $\frak{p}_{{\cal F}} \subseteq A$ to be the set of all elements $a \in A$ such that there exists a closed subgroup $K \subseteq G$ of finite index so that $a \in \frak{p}_K \subseteq A^K$.  We will verify that $\frak{p}_{{\cal F}}$ is a prime ideal in $A$.  We first make the elementary observation that if $a \in \frak{p}_K$ for some $K$, then $a \in \frak{p}_{K^{\prime}}$ for all closed subgroups of finite index $K^{\prime} \subseteq K$.   We verify that $\frak{p}_{{\cal F}}$ is closed under addition.  Suppose $a,a^{\prime} \in \frak{p}_{{\cal F}}$, and fix closed subgroups of finite index $K$ and $K^{\prime}$ so that $a \in \frak{p}_K$ and $a^{\prime} \in \frak{p}_{K^{\prime}}$.  Now, the subgroup $K \cap K^{\prime}$ is also a closed subgroup of finite index, and both $a$ and $a^{\prime}$ are elements of $\frak{p}_{K \cal K^{\prime}}$, consequently $a + a^{\prime} \in \frak{p}_{K \cal K^{\prime}}$, which means that $a + a^{\prime} \in \frak{p}_{{\cal F}}$.  To see that $\frak{p}_{{\cal F}}$ is an ideal, we consider $a \in A$ and $i \in \frak{p}_{{\cal F}}$.  Select closed subgroups of finite index $K$ and $K^{\prime}$ so that $a \in A^K$ and $i \in \frak{p}_{K^{\prime}}$.   We then have $a \in A^{K \cap K^{\prime}}$ and $i \in \frak{p}_{K \cap K^{\prime}}$, and therefore that $a\cdot i \in \frak{p}_{K \cap K^{\prime}}$.   Finally, to see that $\frak{p}_{{\cal F}}$ is prime, consider $a,a^{\prime} \in A$ for which $a \cdot a^{\prime} \in \frak{p}_{{\cal F}}$.  Select a closed subgroup of finite index $K$ so that $a,a^{\prime} \in A^{K}$.  Since $a\cdot a^{\prime} \in \frak{p}_{{\cal F}}$, there is another closed subgroup of finite index $K^{\prime}$ so that $a \cdot a^{\prime} \in \frak{p}_{K^{\prime}}$.   We then have that $a,a^{\prime}$, and $a \cdot a^{\prime} $ all lie in $A^{K \cap K^{\prime}}$ and that $a \cdot a^{\prime} \in \frak{p}_{K \cap K^{\prime}}$.  Because of the primeness of $\frak{p}_{K \cap K^{\prime}}$, it follows that either $a$ or $a^{\prime}$ is an element of $\frak{p}_{K \cap K^{\prime}}$, and consequently that either $a$ or $a^{\prime}$ lies in $\frak{p}_{{\cal F}}$.   Under the identification 
$$ Spec(A) \longrightarrow \column{lim}{\longleftarrow} Spec(A^N)
$$
The fiber $\pi ^{-1}(\frak{p})$ is identified with $F_{\frak{p}}$ above.  

We next note that  the maps $\pi ^K_{K^{\prime}} :Spec(A^K) \rightarrow Spec(B)$ are surjective for any closed normal subgroup of finite index $K \subseteq G$.  This result is immediate from Proposition 1.1 of \cite{exposefive}, applying it to the action of the finite group $G/K$ on $Spec(A^K)$. The surjectivity statement above now follows immediately. 

For the transitivity of the $G$ action on the fibers of the map $\pi$, we proceed as follows.  We first  note that $Spec(A) $ can also be described as  
$$  \column{lim}{\longleftarrow} Spec(A^N)
$$
where $N$ ranges over the closed normal subgroups of finite index in $G$, due to the finality of the closed normal subgroups of finite index in the set of all closed subgroups of finite index.  For each such $N$, and pair of elements $x,y \in Spec(A)$, we consider the subset $\Phi _N(x,y) \subseteq G$ consisting of all $g \in G$ such that $gx_N = y_N$, where $x_N$ and $y_N$ denote the projections of $x$ and $y$ in $Spec(A^N)$.  We are letting $G$ act on $Spec(A^N)$ in the obvious way through the projection $G \rightarrow G/N$.  The set $\Phi_N(x,y)$ is non-empty due to the transitivity statement for finite group actions on rings made in Proposition 1.1 of \cite{exposefive}.  Moreover, if we have $N^{\prime} \subseteq N$, then $\Phi _{N^{\prime}}(x,y) \subseteq \Phi _N(x,y)$.  It is clear that the sets $\Phi_N(x,y)$ are closed subsets of $G$, and  for a fixed choice of $x$ and $y$, have the finite intersection property. Therefore,  due to the compactness of $G$, 
$$ \bigcap _N \Phi _N(x,y) \neq \emptyset
$$
Any element in  $\bigcap _N \Phi _N(x,y)$ will  now carry $x$ to $y$, giving the transitivity statement. We now know that $\pi^{-1}(\frak{p})$ is a Hausdorff transitive continuous $G$-set, and conclude from Proposition \ref{classify} that $\pi^{-1}(\frak{p})$ is $G$-homeomorphic to a $G$-space of the form $G/G^{\prime}$, where $G^{\prime} $ is a closed subgroup of $G$.  
\end{Proof} 

\begin{Corollary}  For continuous affine $G$-schemes,  a categorical orbit space exists, and $Spec(A^G)$ is a choice of such an orbit space. 
\end{Corollary}
\begin{Proof} This is proved as in the finite case, performed in  \cite{exposefive}. 
\end{Proof} 



Let $A$,$B$, and $G$ be as above.  Let $\frak{p}$ be a prime ideal of $A$, and let $G_{\frak{p}}$ denote its stabilizer, which will be called its {\em decomposition group}. $G_{\frak{p}}$ is a closed subgroup of $G$.  Let $k({\frak{p}})$ denote the field of fractions of $A/\frak{p}$.  Similarly, let $\frak{q} = \frak{p} \cap B$, and write $k(\frak{q})$ for the field of fractions of $B/\frak{q}$.  There is an evident inclusion $k(\frak{q}) \subseteq k(\frak{p})$.    The group $G_{\frak{p}}$ acts continuously on $k(\frak{p})$, when $k(\frak{p})$ is equipped with the discrete topology.  There is therefore a homomorphism from $G_{\frak{p}}$ to the Galois  group  $Gal(k(\frak{p}), k(\frak{q}))$, which is in general a profinite group.  This homomorphism is surjective.  This is stated explicitly in \cite{exposefive} in the case where $G$ is finite, and the result in this context follows immediately from this case by passing to inverse limits over Hausdorff finite quotients.  The kernel of this homomorphism is called the {\em inertia group} of $\frak{p}$, and is clearly a closed subgroup of $G_{\frak{p}}$.  We denote it by $I_{\frak{p}}$.  

We will also study the notion of etale actions in the context of profinite actions. 

\begin{Definition}
Let $G$ be a profinite group, and let $A$ be a continuous affine $G$-scheme.   We say the action is {\em etale} if each of the inclusions $A^G \hookrightarrow A^K$ is an etale map for every closed subgroup of finite index $K$.  This is clearly equivalent to the condition that each inclusion $A^G \hookrightarrow A^N$ is a Galois cover with group $G/N$, where $N$ varies over all the closed normal subgroups of finite index.  
\end{Definition} 

We have the following analogue of  Corollary 2.4 of \cite{exposefive} for verifying that an action is etale.

\begin{Proposition}
 Let $A$ be a commutative ring equipped with the discrete topology, and let a profinite group $G$ act continuously on $A$.  Then the action is etale if and only if the inertia group $I_{\frak{p}}$ is trivial for all prime ideals $\frak{p}$ of $A$.  
\end{Proposition}
\begin{Proof} This is a straightforward consequence of the case of a finite group $G$, for which a proof is given in \S 2 of  \cite{exposefive}. 
\end{Proof}

There is a version of this criterion which will be useful to us.  Let $A$ and $G$ be as above, and let $\frak{q}$ denote any prime ideal of $A^G$.  Let $\Omega$ denote any algebraically closed field which admits an embedding $i: k(\frak{q}) \rightarrow \Omega$.  We let $\frak{H}(A,\Omega, i)$ denote the set of all ring homomorphisms $A \rightarrow \Omega$ extending the composite 
$$A^G \rightarrow A^G/\frak{q} \hookrightarrow k(\frak{q}) \stackrel{i}{\longrightarrow}\Omega$$
For any closed normal subgroup $N$ of $G$ of finite index, we similarly have the set $\frak{H}(A^N, \Omega, i)$ and we have restriction maps 
$$ \frak{H}(A^{N^{\prime}}, \Omega, i) \longrightarrow \frak{H}(A^N, \Omega, i)
$$
whenever $N^{\prime} \subseteq N$.  It is readily verified that 
$$ \frak{H}(A, \Omega, i) \cong \column{lim}{\leftarrow} \frak{H}(A^N, \Omega, i)
$$
and therefore further that $\frak{H}(A, \Omega, i)$ is equipped with a continuous $G$ action. 
\begin{Proposition} The $G$-space $\frak{H}(A, \Omega, i)$ is transitive, and is $G$-homeomorphic to the quotient $G/I_{\frak{p}}$, where $\frak{p}$ is the kernel of some ring homomorphism $A \rightarrow \Omega$ extending $i$.  
\end{Proposition}
\begin{Proof} For the finite case, this is proved in \S 2  of \cite{exposefive}.  The extension to the profinite case is direct. 
\end{Proof}
\begin{Corollary} A continuous affine $G$-scheme is etale if and only if the $G$ actions on all sets $\frak{H}(A,\Omega, i)$ are free, for all prime ideals $\frak{p}$ of $A^G$ and homomorphisms $i: k(\frak{p}) \hookrightarrow \Omega$ for all algebraically closed fields $\Omega$.  
\end{Corollary}

In order to define the equivariant algebraic $K$-theory of continuous affine $G$-schemes, we extend the quasi-Noetherian condition to such actions. 

\begin{Definition}  We say that a continuous affine $G$-scheme $(A,G)$ is coherent if $A$ is a coherent ring.  We say it is {\em totally coherent} if  $A^N$ is a coherent  ring for all closed subgroups $N \subseteq G$.  
\end{Definition} 
We now define the category of twisted $G$ modules for coherent continuous affine $G$-schemes. 

\begin{Definition}  Let $(A,G)$ be a ring $A$ equipped with a continuous action by a profinite group $G$, where we assume that $A$ is equipped with the discrete topology.  By a $G$-twisted module over $A$, we mean an $A$ -module $M$, endowed with the discrete topology, together with a continuous action by $G$, which is compatible with the $G$ action on $A$ via the formula 
$$  g \cdot (\alpha m) = \alpha ^g \cdot gm
$$
The  $G$-twisted $A$-modules form an abelian category.  A $G$-twisted $A$-module is said to be {\em finitely presented} if the underlying module is a finitely presented $A$-module.    In the case of  coherent continuous affine $G$-schemes,  the finitely presented  $G$-twisted $A$-modules also form an abelian category.  We will write $\mbox{\em Mod}^G(A)$ for the category of finitely presented $G$-twisted $A$-modules.  
\end{Definition}

We now define the equivariant algebraic $K$-theory of coherent continuous affine $G$-schemes.   
\begin{Definition}   Given a coherent continuous affine  $G$-scheme $(A,G)$, by the equivariant $K$-theory of $A$, denoted $K^G(A)$, we will mean the $K$-theory of the Waldhausen category associated to the abelian category $\mbox{\em Mod}^G(A)$.  
\end{Definition}

With this definition, we would now like to provide an analogue to Proposition \ref{free} for actions of a profinite group. 

Begin with some results on coherent rings.  

\begin{Proposition}  \label{descendproperties} Let $f: A \rightarrow B$ be a faithfully flat morphism of rings.  Then an $A$-module $M$ is finitely generated (respectively finitely presented) if and only if $B \column{\otimes}{A} M$ is finitely generated (finitely presented). 
\end{Proposition}
\begin{Proof}  This is Corollary 1.11 of \cite{grothendieck}
\end{Proof}
\begin{Corollary} \label{descendcoherence} Let $f: A \rightarrow B$ be a faithfully flat morphism of rings.  If $B$ is coherent, then so is $A$.  
\end{Corollary} 
\begin{Proof}  We will need to prove that the kernel of a homomorphism $f: F \rightarrow F^{\prime}$ between finitely generated $A$-modules is finitely generated.  Because of the flatness, we have an exact sequence 
$$ 0 \rightarrow B \column{\otimes}{A} Ker(f)  \rightarrow B  \column{\otimes}{A}F \rightarrow   B  \column{\otimes}{A}F^{\prime} 
$$
The coherence of $B$ now gives that $B \column{\otimes}{A} Ker(f)$ is finitely generated.  Proposition \ref{descendproperties} now gives the result.  
\end{Proof} 

\begin{Corollary}  Let $A$ be a coherent commutative ring equipped with the discrete topology, and let a profinite group $G$ act continuously on $A$.  If the $G$-action is etale, then the fixed point subring $A^G$ is coherent.  
\end{Corollary} 
\begin{Proof}  The action being etale implies that the inclusion  $A^G \hookrightarrow A$ is faithfully flat.  
\end{Proof} 

\begin{Corollary}  For any continuous etale action of a profinite group $G$ on a coherent ring $A$, $A^G$ is coherent  there is a natural map 
$$ K(A^G) \rightarrow K^G(A^G) \rightarrow K^G(A)
$$  In the terminology of schemes $X = Spec(A)$, there is a  map $\eta_X: K(X/G) \rightarrow K^G(X)$  when $X$ is a coherent affine scheme, natural for  flat maps of schemes.  
\end{Corollary} 

We wish to prove that $\eta _X$ is an equivalence.  In the case of a finite group $G$, this result follows from Theorem II.5.1 of \cite{ojanguren}, by proving that the  functors $ M \rightarrow M^G$ from $\mbox{Mod}^G(A) $  to $\mbox{Mod}^{\{ e \}}(A^G)$  and $M \rightarrow B \tover{A}M$ from $\mbox{Mod}^{\{ e \}}(A^G)$ to $\mbox{Mod}^G(A) $  are inverse equivalences of categories.   The fact that both functors preserve the subcategories of finitely presented modules follows from the faithful flatness of the inclusion $A^G  \hookrightarrow A$ together with Proposition \ref{descendproperties} above.  

\begin{Proposition} \label{orbitequiv} The natural transformation $\eta_X$ is an equivalence for any continuous etale action of a profinite group $G$ on an affine scheme $X = Spec(A)$, where $A$ is coherent. 
\end{Proposition} 
\begin{Proof}    We may  construct the categories  $\mbox{Mod}^{G/N}(A^N)$ for any closed normal subgroup $N \subseteq G$ of finite index.  We note that all the rings $A^N$ are coherent by Corollary \ref{descendcoherence}, and that all $N \subseteq N \pr$ the inclusions $A^N \hookrightarrow A^{N \pr }$ are faithfully flat, because they are etale ring extensions.  It follows that we have a directed system $N \rightarrow \mbox{Mod}^{G/N}(A^N)$, where the functors 
$$   \mbox{Mod}^{G/N}(A^N) \rightarrow  \mbox{Mod}^{G/N\pr}(A^{N\pr})
$$
are defined to be  $A^{N \pr} \tover{A^N} -$.  There is therefore a functor 
$$j:\column{colim}{\column{\longrightarrow}{N}} \mbox{Mod}^{G/N}(A^N) \rightarrow \mbox{Mod}^G(A)$$
which we claim is an equivalence of categories.  This will suffice to prove the result, because each of the functors $ \mbox{Mod}^{\{e\}}(A^G) \rightarrow  \mbox{Mod}^{G/N}(A^{N})$ is an equivalence of categories by the validity of the result for finite group actions mentioned above.  

We proceed to prove that $j$ is an equivalence of categories.  We first observe that $j$ is surjective on isomorphism classes of objects.  
In order to do this, we will first give a description of $G$-twisted modules over $A$ whose underlying $A$-module is free  in terms of matrices.  Suppose $F$ is a free module with continuous $G$-action twisted compatibly with the action of $G$ on $A$.  We fix a basis $B = \{ b_1, \ldots , b_t \}$  for $F$.  Because of the continuity of the action, it is clear that there is a closed normal subgroup $N$ of finite index in $G$ so that $N$ fixes the set  $B$ elementwise.  It follows that the action of an element $g \in G$ on an element of $B$ depends only on its coset in $G/N$.  For each coset $\gamma$  in $G/N$, we then obtain a $t \times t$ matrix $M(\gamma)$ with entries in $A$, defined by $M(\gamma)(b_j) = \sum _i M(\gamma) _{ij} b_i$.    Let $N$ be such that  $M(\gamma)_{ij} \in A^N$ for all $\gamma $ and pairs $(i,j)$.  The coefficients of all the matrices $M(\gamma)$ are therefore acted on trivially by $N$, and therefore the action of $G$ on the free $A^N$-submodule generated by $B$ factors through $G/N$.  The matrices $M(\gamma)$ clearly satisfy the relations 
 \begin{equation} \label{relations}
 M(\gamma _1 \gamma _2) = M(\gamma _1) \cdot M(\gamma _2)^{\gamma _1}
 \end{equation}
  where the superscript indicates entrywise action of $G$ on the matrix $M(\gamma _2)$.  Sets of matrices $\{M_{\gamma} \}_{\gamma} \in G$ satisfying these relations in turn determine a $G/N$-twisted  module over $A^N$, and after applying $A \column{\otimes}{A^N}-$, over $A$. It is furthermore clear that the result is isomorphic to the original $G$-twisted $A$-module, so we have proved that any continuous twisted $G$-module over $A$, with free underlying module,  is obtained by extension of scalars from a free  twisted $G/N$-module over $A$, for some closed normal subgroup of finite index  $N$. 

Next, we will prove that given any homomorphism of $G$-twisted $A $-modules 
$$  f: A \column{\otimes}{A^N} F_0 \rightarrow A \column{\otimes}{A^N} F_1
$$
where $F_0$ and $F_1$ are $G/N$-twisted $A^N$-modules whose underlying modules are finitely generated and free, there is an $N \pr \supseteq N$ and homomorphism of $G/N \pr $-twisted $A^{N \pr}$-modules 
$$ \overline{f}: A^{N \pr}\column{\otimes}{A^N} F_0 \rightarrow A^{N \pr } \column{\otimes}{A^N } F_1
$$
so that $f = A \column{\otimes}{A^{N \pr} } \overline{f}$.  To see this, we observe that given bases for $F_0$ and $F_1$, $f$ determines a matrix $M_f$ with entries in $A$.  The $G$-twisted modules $F_0$ and $F_1$ determine families of matrices $ \{ M_0(\gamma)  \}_{\gamma \in G}$ and $ \{ M_1(\gamma)  \}_{\gamma \in G}$ with entries in $A^N$, and the fact that $f$ is a homomorphism of $G$-twisted $A$-modules means that the equations 
\begin{equation} \label{reltwo} M_1({\gamma}) M_f^{\gamma} = M_f M_0(\gamma)\end{equation} hold  for all $\gamma \in G$, when $M_0(\gamma)$ and $M_1(\gamma)$ are regarded as matrices with entries in $A$ via the homomorphism $A^N \rightarrow A$.   As above, there is a $N \pr \supseteq  N$ so that all the entries of $M_f$ lie in the image of the homomorphism $A^{N \pr}  \rightarrow A$, and we let $\overline{M}_f$ denote the matrix $M_f$ regarded as a matrix with entries in $A^{N \pr} $.  We know that the relations (\ref{reltwo}) hold in $A$, but this means that the relations
$$ M_1({\gamma}) \overline{M}_f^{\gamma} = \overline{M}_f M_0(\gamma)$$
hold in $A^{N \pr} $ by the faithful flatness of the homomorphism $A^{N \pr}  \rightarrow A$.  The required  homomorphism $\overline{f}$ is now the homomorphism represented by the matrix $\overline{M}_f$.  

Next consider any finitely presented $G$-twisted $A$-module $M$, with a presentation 
\begin{equation} \label{presentation}  F_1 \stackrel{f}{\rightarrow} F_0 \rightarrow M \rightarrow 0
\end{equation} 
where $F_0$ and $F_1$ are $G$-twisted $A_F$-modules with finitely generated free underlying modules.  Our first argument above shows that there is closed normal subgroup of finite index $N$ so that $F_0$ and $F_1$ are extended from $A^N$, and the second shows that there is a $N \pr \supseteq N$ so that $f $ is extended from $A^{N \pr} $.  Let $\overline{F}_0, \overline{F}_1$ be $G/N \pr$-twisted $A^{N \pr}$-modules (with free underlying modules) and $\overline{f} $ be a homomorphism $\overline{f}: \overline{F}_1 \rightarrow \overline{F}_0$ of $G/N\pr $-twisted $A^{N \pr} $-modules so that the presentation \ref{presentation} is obtained by applying $A \column{\otimes}{A^{N \pr}} -$ to the exact sequence  
$$  \overline{F}_1 \stackrel{\overline{f}}{\rightarrow} \overline{F}_0 \rightarrow C \rightarrow 0
$$ of $G/N \pr $-twisted $A^{N \pr}$-modules, where $C$ denotes the cokernel of $\overline{f}$.  It follows by flatness that $M \cong A \column{\otimes}{A^{N \pr} } C$.  It follows that the functor $j $ above is surjective on isomorphism classes of objects.  

To prove injectivity on isomorphism classes of objects, we suppose that we are given a closed normal subgroup of finite index $N \subseteq G$ and finitely presented $G/N$-twisted $A^N$-modules $M $ and $M \pr $, so that $A \column{\otimes}{A^N} M$ and $A\column{\otimes}{A^N} M\pr$ are isomorphic.  We suppose we are given an isomorphism $f: A \column{\otimes}{A^N} M \rightarrow A\column{\otimes}{A^N} M\pr$ and its inverse $f^{-1}$.  Both $f$ and $f^{-1}$ are  of the form $A \column{\otimes}{A^{N \pr} } \overline{f}$ and $A^{N \pr} \column{\otimes}{A^N} \overline{f^{-1}}$ for some $N\pr$ and  homomorphisms $$\overline{f}: A^{N \pr} \column{\otimes}{A^N} M  \rightarrow A^{N \pr}  \column{\otimes}{A^N } M\pr $$
and 
$$\overline{f^{-1}}: A^{N \pr}  \column{\otimes}{A^N} M \pr   \rightarrow A^{N \pr}\column{\otimes}{A^N} M$$
But $\overline{f}$ and $\overline{f^{-1}}$ are also inverses, because they become so over $A$, and because of the faithful flatness of $A^{N \pr} \rightarrow A$. It follows that $\overline{f}$ is an isomorphism, which gives the injectivity result for isomorphism classes of objects. 

Finally we will need to check that for any pair of objects $x,y$ in $\underline{C} = \column{colim}{\column{\longrightarrow}{N}} \mbox{Mod}^{G/N}(A^N)$, $j$ induces a bijection $$Hom_{\underline{C}}(x,y) \stackrel{\sim}{\rightarrow} Hom_{\mbox{Mod}^G(A)}(j(x), j(y))$$  The injectivity of this map is immediate from the faithful flatness of all the maps $ A^N \rightarrow A$.  To see the surjectivity, we assume that the objects $x$ and $y$ are represented by $G/N$-twisted $A^N$-modules $M$ and $M\pr$, respectively, equipped with presentations
$$  F_1 \stackrel{\alpha}{\rightarrow} F_0 \rightarrow M  \rightarrow 0
$$
and 
$$ F_1\pr \stackrel{\beta}{\rightarrow} F_0 \pr  \rightarrow M \pr \rightarrow 0
$$
A homomorphism $f$  from $M$ to $M \pr$ can be represented by a pair of homomorphisms $f_0 : F_0 \rightarrow F_0 \pr$ and $f_0 : F_1 \rightarrow F_1 \pr$ so that $ f_0 \alpha = \beta f_1$.  These homomorphisms can be represented by matrices, and by arguments similar to the ones above we find that these matrices can be represented over  $A^N$ for some closed normal subgroup of finite index $N \subseteq G$.  This gives the result.  
\end{Proof}

 \subsection{Homotopy properties}

Let $G$ be any finite group, and let $G$ act by automorphisms on a vector space $V$ over a field $k$ via a representation $\rho$f.  We suppose that $\# (G)$ is relatively prime to $char(k)$. We may form the skew polynomial ring $ k[V]_{\rho} [G]$ (see \cite{goodearl}).  The following is now a formal observation.

\begin{Proposition}  \label{veryformal} The category $\mbox{Mod}^G(V)$ ($V$ is regarded as an affine scheme) is equivalent to the category of finitely generated modules over  the regular Noetherian ring $ k[V]_{\rho} [G]$. 
\end{Proposition} 
  Let $\mbox{Rep}[G]$ denote the category of finite dimensional $k$-linear representations of $G$, which can be identified with the category of finitely generated modules over the group ring $k[G]$.  We have an inclusion $i: k[G] \hookrightarrow k[V]_{\rho} [G]$ and a projection $\pi : k[V]_{\rho} [G] \rightarrow k[G]$ of rings, sending the vector space $V$ identically to $0$.  The following result is now an immediate consequence of Theorem 7 of \cite{quillen}.

\begin{Theorem} \label{eqhomotopy} The ring homomorphisms $i$ and $\pi$ give rise to exact functors 
$$ k[V]_{\rho} [G] \column{\otimes}{k[G]} - : \mbox{\em Rep}[G] \rightarrow \mbox{ \em Mod}^G(V)$$ 
and 
$$k[G] \column{\otimes}{k[V]_{\rho} [G]} - :   \mbox{\em Mod}^G(V)  \rightarrow \mbox{\em Rep}[G]
$$
Both functors induce equivalences on $K$-theory spectra.  
\end{Theorem}  
We will also need to describe how induction of representations works in this context.  Let $H \subseteq G$ be an inclusion of finite groups, and let $\rho $ denote a $k$-linear representation of $H$, with $V_{\rho}$ as the representation space.  The induced representation $i_H^G(\rho)$ can be defined via the operation $$k[G]\column{\otimes}{k[H]} -: \mbox{Mod}(k[H])\rightarrow \mbox{Mod}(k[G])$$ where for a ring $A$ the category of finitely generated left $A$ -modules is denoted by $\mbox{Mod}(A)$. We denote the representation space of $i_H^G(\rho)$ by $V_{i_H^G(\rho)}$.  There is a natural $H$-equivariant inclusion $V_{\rho} \hookrightarrow V_{i_H^G(\rho)}$, and consequently an inclusion of rings 
$$k[V_{\rho}]_{\rho}[H] \hookrightarrow k[V_{i_H^G(\rho)}]_{{i_H^G(\rho)}}[G]
$$ 
\begin{Proposition} \label{inductionsquare}
Given $H,G$, and $\rho$ as above, the following diagram of functors commutes.  

$$ \begin{CD} 
\mbox{{\em Mod}}(k[H])  @>\phantom{\mbox{xxxxxxx}}{k[G]\column{{\otimes}}{k[H]} -}\phantom{\mbox{xxxxxxxxxx}} >> \mbox{\em Mod}(k[G]) \\
@V{k[V_{\rho}]_{\rho}[H] \column{\otimes}{k[H]} - }VV   @VV{k[V_{i^G_H(\rho)}]_{i_H^G(\rho)}[G]\column{\otimes}{k[G]}-}V \\
\mbox{\em Mod}(k[V_{\rho}]_{\rho}[H])  @>{k[V_{i^G_H(\rho)}]_{i_H^G(\rho)}[G]\column{\otimes}{k[V_{\rho}]_{\rho}[H]}-}>> \mbox{\em Mod}(k[V_{i^G_H(\rho)}]_{i_H^G(\rho)}[G])
\end{CD} 
$$

\end{Proposition} 
\subsection{Completions}
We recall from \cite{completion} that given a homomorphism of commutative ring spectra $f:A \rightarrow B$, and an $A$-module $M$, one may construct the {\em derived completion} $M^{\wedge}_f$.  Here are some  examples of this construction.

{\bf Example:} Any commutative ring may be thought of as a commutative  ring spectrum via the Eilenberg-MacLane construction.  If $f: A \rightarrow B$ is a surjective homomorphism of commutative Noetherian rings, and $M$ is a finitely generated $A$-module, then the derived completion $M^{\wedge}_f$ is the Eilenberg-MacLane construction on the usual completion $M^{\wedge}_I$, where $I$ is the kernel of $f$.  

{\bf Example} Let $f: S^0 \rightarrow \Bbb{H}_l$ be the mod $l$ Hurewicz map, and let $X$ be any spectrum (and therefore an $S^0$-module).  Then $X^{\wedge}_f$ is the usual Bousfield-Kan  notion of completion of the spectrum at the prime $l$ (see \cite{bousfield}).  

{\bf Example:}  Let $A = R[\Bbb{Z}_l]$  be the representation ring of the additive group of the $l$-adic integers, defined as the direct limit of the representation rings of the finite cyclic $l$-th power order groups.  Let $f: A \rightarrow \Bbb{F}_l$ be the augmentation followed by reduction mod $l$.  Then $A^{\wedge}_f$ is equivalent to the $l$-adic completion of the integral group ring of the simplicial group of points on the circle group.  In particular, the homotopy groups are $\cong \Bbb{Z}_l$ for dimensions 0 and 1, and are $\cong 0$ otherwise.

It has a number of useful properties which can be found in \cite{completion}, and we record some of them here. Let $f: A \rightarrow B$ be  a homomorphism of commutative ring spectra.

\begin{Proposition} \label{compcofib} Suppose that 
$$  M \rightarrow N \rightarrow P
$$
is a cofibration sequence of $A$-modules.  Then the naturally defined sequence 
$$ M^{\wedge}_f \rightarrow N ^{\wedge}_f \rightarrow P ^{\wedge}_f 
$$
is also a cofibration sequence.  
\end{Proposition} 

\begin{Proposition}\label{criterionone}  Suppose we are given a homomorphism $\lambda: M \rightarrow N$ of $A$-modules. Suppose further that the derived smash product map $id_B \wedge \lambda : B \column{\wedge}{A} M \rightarrow B \column{\wedge}{A}N$ is an equivalence of spectra.  Then the natural map $\lambda ^{\wedge}_f: M^{\wedge}_f \rightarrow N^{\wedge}_f$ is an equivalence of spectra.  In particular, if $ B \column{\wedge}{A}M \simeq *$, then $M^{\wedge}_{f} \simeq *$.  
\end{Proposition} 

{\bf Remark:} By the derived smash product, we mean the smash product construction applied to cofibrant replacements for $M$ and $N$.  If they are already cofibrant, then one can apply the smash product directly. 

\begin{Proposition} \label{algcriterion}Suppose that $A$ is a (-1)-connected  commutative ring spectrum, and that $M$ is a connective  $A$-module.  Suppose further that we have a homomorphism of commutative ring spectra $A \rightarrow \Bbb{H}_l$, and that for all $n$ and $i$,  we have that the groups 
$Tor^{\pi _0 A}_i(\pi_nM, \Bbb{F}_l)$ vanish.  Then we have $M \column{\wedge}{A} \Bbb{H}_l \simeq *$, and in particular  $M^{\wedge}_{\Bbb{H}_l} \simeq *$. 
\end{Proposition}
\begin{Proof} We have the K\"{u}nneth spectral sequence converging to $\pi _* (M \column{\wedge}{A} \Bbb{H}_l)$, whose $E_2 $-term is $Tor_{\pi_*A}(\pi_*M, \Bbb{F}_l)$. We also have the inclusion $\pi _0A \hookrightarrow \pi _*A$, as the inclusion of the degree zero component in the graded ring $ \pi _*A$.  We let $B_*$ denote the commutative graded ring $\Bbb{F}_l\column{\otimes}{\pi _0A} \pi _*A $.  We can now describe the functor $ \Bbb{F}_l\column{\otimes}{\pi _*A} -$ from graded $\pi _*A$-modules to graded $\Bbb{F}_l$-vector spaces as the composite of the functor $\Bbb{F}_l\column{\otimes}{\pi _0 A} -$ from graded $\pi _*A$-modules to graded $B_*$-modules with the functor $\Bbb{F}_l\column{\otimes}{B_*} -$ from the category of graded $B_*$-modules to graded $\Bbb{F}_l$-vector spaces.  The Grothendieck spectral sequence now shows that the vanishing of the derived functors of the first functor forces the vanishing of the derived functors of the composite. 
\end{Proof} 

\begin{Proposition}  Let $f: A \rightarrow B$ be a homomorphism of $(-1)$-connected commutative  ring spectra, with $\pi _0 (f) $ surjective,  and let $M$ be an $A$-module. Let $I$ denote the kernel of $\pi _0 (f)$.  Suppose that the $\pi _0A$-modules $\pi _j M$ are such that for each $j$, there is an integer $e_j$ so that $I^{e_j} \cdot \pi _j M = \{ 0 \}$.  Then the natural map $\eta : M  \rightarrow M^{\wedge}_f$ is an equivalence.
\end{Proposition}
\begin{Proof} Theorem 7.1 of \cite{completion} reduces us to the case where $A$ is  the ring $\pi _0A$, $B$ is the ring $\pi _0 B$,   and $M$ is a module for which $I^eM = \{ 0 \}$ for some $e$.  We then obtain the finite filtration 
$$ \{ 0 \} = I^e M \subseteq I^{e-1}M \subseteq \cdots \subseteq I^2M \subseteq IM \subseteq M
$$
each of whose quotients is annihilated by $I$.  It follows from (6) of Proposition 3.2 of \cite{completion} that $\eta : I^sM/I^{s+1}M \rightarrow (I^sM/I^{s+1}M )^{\wedge}_f$ is an equivalence.  A straightforward application of Proposition \ref{compcofib} now gives the result. 
\end{Proof} 

The following technical criterion will be useful in later sections. 

\begin{Proposition} \label{technical} Let 
$$ M_0 \stackrel{f_0}{\rightarrow} M_1 \stackrel{f_1}{\rightarrow}M_2 \stackrel{f_2}{\rightarrow}\cdots
$$
be a directed system of  module spectra over the commutative  ring spectrum $A$.  Suppose further that we are given a homomorphism $f : A \rightarrow B$ of commutative  ring spectra. Finally, suppose that the $Ind$-groups 
$$ \pi _t(B \column{\wedge}{A} M_0) \rightarrow   \pi _t(B \column{\wedge}{A} M_1) \rightarrow   \pi _t(B \column{\wedge}{A} M_2) \rightarrow \cdots 
$$
are $Ind$-trivial.  Then the $A$-module $M_{\infty} = \mbox{ \hspace{-.3cm}} \column{hocolim}{r} M_r$ satisfies $B \column{\wedge}{A}M_{\infty} \simeq *$, and therefore $(M_{\infty})^{\wedge}_f \simeq *$. \end{Proposition}
\begin{Proof}  Because smash products are colimit constructions, it is possible to commute the constructions $B \column{\wedge}{A} -$ and colim, which shows that $B \column{\wedge}{A} M_{\infty} \simeq *$.  
\end{Proof}

One situation in which one can  prove that the hypotheses of Proposition \ref{technical} are satisfied is treated in the following result. 

\begin{Proposition} \label{localize} Suppose $f:A \rightarrow B$ is a homomorphism of $(-1)$-connected commutative  ring spectra, which is surjective on $\pi _0$,  and suppose $L$ is a connective commutative $A$-algebra.
, equipped with a homomorphism of commutative $A$-algebras $\epsilon _L : L \rightarrow B$ so that the diagram
$$
\begin{diagram}\node{A} \arrow{e} \arrow{se,b}{f}    \node{L} \arrow{s,b}{\epsilon _L} \\
\node{}   \node{B}
\end{diagram}
$$
commutes. 
 Let $I_A \subseteq \pi _0 A$ denote the ideal $Ker(f)$.  
Consider now an Ind-$A$-module of the form 

\begin{equation} \label{system} 
 L \stackrel{\times \alpha _0}{\longrightarrow} L \stackrel{\times \alpha _1}{\longrightarrow} L \stackrel{\times \alpha _2}{\longrightarrow} \cdots
\end{equation} 
where $\alpha _i $ is an element of the ideal $I_A \cdot \pi _0 L$.  Then Ind-groups
\begin{equation} \label{indsystemone}
\pi _t (B \column{\wedge}{A} L )\stackrel{\pi _t(\times \alpha _0)}{\longrightarrow} \pi _t(B \column{\wedge}{A}L)\stackrel{\pi _t (\times \alpha _1)}{\longrightarrow}\pi _t (B \column{\wedge}{A}L )\stackrel{\pi _t (\times \alpha _2)}{\longrightarrow} \cdots
\end{equation}
are Ind-trivial, and therefore by   Proposition \ref{technical}, the homotopy colimit $L_{\infty} = \mbox{}\column{hocolim}{i} L$ satisfies $B \column{\wedge}{A} L_{\infty} \simeq *$.

\end{Proposition} 
\begin{Proof} From the theory of derived tensor products for modules over ring spectra, as developed in \cite{schwede}, there is a spectral sequence with $E_2$-term $Tor^{\pi _*A}_* (\pi _*B, \pi _ * L)$  converging to  $\pi _* (B \column{\wedge}{A} L)$.  We now have an inductive system of spectral sequences corresponding to diagram \ref{indsystemone} above, and  claim that maps in the system are identically zero.  To prove this, we must show that any $\lambda \in I_A \cdot \pi _0 L$ induces the zero map on all $Tor$-groups.  To see this, we write such a $\lambda $ as 
$$\lambda = \sum _s i_s \lambda _s
$$
where $i_s \in I_A$, and $\lambda _s \in \pi _0 L$.  If we can prove that each $i_s$ induces the zero map on $Tor$-groups, we will have verified the claim.  But this is clear, since multiplication by elements of $I_A$ induce the zero map on $\pi _* B$.  
 Further, the colimit of the spectral sequences converges to the colimit of the homotopy groups due to the strong convergence of the spectral sequences which is guaranteed by the connectivity hypothes on $A$ and $B$, which gives the result.   \end{Proof} 
 
 \begin{Corollary} \label{tootechnical} Suppose we are given $A,B,L, \mbox{ and } I_A,$ are as in Proposition \ref{localize},  and that $\epsilon _L : L \rightarrow B$ is a homomorphism of commutative ring spectra so that the diagram 
 
 $$
\begin{diagram}\node{A} \arrow{e} \arrow{se,b}{f}    \node{L} \arrow{s,b}{\epsilon _L} \\
\node{}   \node{B}
\end{diagram}
$$
commutes.  Let $I_L$ be the kernel of $ \pi _0 (\epsilon _L )$.  Suppose further that there is an integer $N$ so that $I_L^N \subseteq I_A \cdot \pi _0(L)$.   We suppose  we are given an Ind-$A$-module as in (\ref{system}) above, but with the assumption that $\alpha _i \in I_L$ rather than that $\alpha _i \in I_A$.  Then the Ind-groups in (\ref{indsystemone}) are trivial, and we again have that $B \column{\wedge}{A}L_{\infty} \simeq *$
 \end{Corollary} 
 \begin{Proof}  We note that any $N$-fold composite in the system (\ref{indsystemone}) has its image in $I_A \cdot \pi _0 (L)$, and a standard cofinality argument gives the result.  
 \end{Proof}

\section{Totally torsion free groups}

We recall that  a topological group $G$ is said to be {\em profinite } if it is Hausdorff, compact, and totally disconnected.  These properties imply that as a topological group, it can be expressed as an inverse limit of finite groups.  Let $l$ denote a prime.  A profinite $l$-group is a profinite group whose only finite quotients are finite $l$-groups.  For any profinite group $G$, we define the {\em commutator subgroup} of $G$, $[G,G]$ to be the closure of the subgroup generated by elements of the form $ghg^{-1}h^{-1}$ in $G$.   The subgroup $[G,G]$ is normal, and the quotient $G/[G,G]$ will be denoted by $G^{ab}$.  It is the universal profinite abelian quotient of $G$.  One also has the following. 

\begin{Proposition} \label{profiniteprops}
Any closed subgroup of a profinite group is itself a profinite group.  The quotient of a profinite group by a closed normal subgroup is also a profinite group.  In particular, any subgroup containing a closed subgroup of finite index is itself a closed subgroup.  
\end{Proposition} 
This is in $\S$2.1 of \cite{ribes}.

\noindent We say a profinite $l$-group $G$ is {\em totally torsion free}  if for every closed subgroup $K \subseteq G$ the quotient group $K^{ab}$ is torsion free.  Free abelian pro-$l$ groups and free pro-$l$ groups are totally torsion free,  and free products (in the category of profinite pro-$l$ groups) of totally torsion free pro-$l$ groups are totally torsion free.  However, the condition is quite restrictive, as the following example shows.  

{\bf Example:} The {\em $2$-adic Heisenberg group}, i.e. the group of $3 \times 3$ upper triangular matrices with entries in $\Bbb{Z}_2$ and ones on the diagonal is {\em not} TTF, although it is easily seen to be torsion free.  To see this, we note that the subgroup of matrices of the form 
 
 \[  \left[ \begin{array}{ccc}
 1 & 2k  &  l  \\
 0  &  1  &  m  \\
 0  &  0  &  1  
\end{array} 
\right]
\]
has abelianization isomorphic to 

$$ \Bbb{Z}_2 \oplus \Bbb{Z} _2 \oplus \Bbb{Z}/2 \Bbb{Z}
$$
The matrix 
\[ \left[ \begin{array}{ccc}
1  &  0  &  1  \\
0  &  1  &  0  \\
0  &  0  &  1
\end{array} \right] \]
is not a commutator, but its square is. 

\begin{Proposition} \label{totorfree}
Let $k$ be a field containing all roots of unity, with absolute Galois  group $G_k$.  Then $G_k$ is totally torsion free.\end{Proposition}
\begin{Proof} It clearly suffices to prove that $G_k^{ab}$ is torsion free for all fields $k$ containing all the roots of unity.  Kummer theory (see e.g. \cite{lang}, \S VI.8) tells us that $G_k^{ab}$ is the Galois group of the infinite extension $k^{ab} = \bigcup _mk((k^*)^{\frac{1}{m}})$, and therefore that any non-identity element $g \in G_k^{ab}$ acts non-trivially on an element of the form $\sqrt[m]{\kappa}$ for some $m$ and $\kappa \in k$.  Suppose now that $g^n = e$ in $G_k^{ab}$.   We have $g \cdot \sqrt[m]{\kappa} = \zeta \sqrt[m]{\kappa}$ for some $\zeta \in \mu _m (k)$.  Now choose any $n$-th root $\xi$ of $\sqrt[m]{\kappa}$.  We see that $g \cdot \xi = \eta \xi$ for some $mn$-th root of unity $\eta $ satisfying $\eta ^n = \zeta$.  Now, $g^n \cdot \xi = \eta ^n \xi = \zeta \xi \neq \xi$, which contradicts the assumption that $g$ is $n$-torsion.  
\end{Proof} 

{\bf Remark:} This structural property of Galois groups will allow us to identify completion constructions based on representation theory with $K$-spectra of certain schemes which appear to be the algebraic geometric version of the classifying space construction in topology.  It would be interesting to consider how this property relates with the Bloch-Kato conjecture.

\section{Construction of some continuous affine $G$-schemes}\label{schemeactions}

Let $G$ be a profinite $l$-group.  We will need to construct continuous affine $G$-actions so as to be able to assemble the scheme ${\cal E}G$ we require.    In this section, we will construct such actions   attached to affine representations over $\Bbb{Z}_l$  of  $G$.  By taking (infinite) products of schemes constructed in this way, we will obtain the scheme  ${\cal E}G$  discussed in the introduction when $G$ is totally torsion free.  

\begin{Definition} An {\em affine $\Bbb{Z}_l$-module} is a pair $(F,X)$, where $F$ is a free finitely generated $\Bbb{Z}_l$-module, and $X$ is a continuous simply transitive $F$-space.  Here $F$ is regarded as a topological abelian group, and the simple transitivity means that for every $x \in X$, the natural map $F \rightarrow X$ given by $f \rightarrow f \cdot x$ is a homeomorphism.  A {\em based affine $\Bbb{Z}_l$-module} is a triple $(F,X,B)$ so that the pair $(F,X)$ is an affine $\Bbb{Z}_l$-module and $B$ is a basis for the $\Bbb{Z}_l$-module $F$.  An {\em isomorphism}  of affine $\Bbb{Z}_l$-modules from $(F,X)$ to $(F^{\prime}, X^{\prime})$ is a pair $(\varphi ,\chi )$, where $\varphi :F \rightarrow F^{\prime}$ is a continuous isomorphism of $\Bbb{Z}_l$-modules and $\chi : X \rightarrow X^{\prime}$ is a homeomorphism, and where the relationship 
$$ \chi (f \cdot x) = \varphi (f) \cdot \chi (x)
$$
holds for all $f$ and $x$.  An isomorphism $\varphi$  of based affine $\Bbb{Z}_l$-modules is an isomorphism of the underlying affine $\Bbb{Z}_l$-modules which preserves $\Bbb{Z}_l$-multiples of basis elements, i.e. for all $z \in \Bbb{Z}_l$ and $b \in B$ there exist $z^{\prime} \in \Bbb{Z_l}$ and $b^{\prime} \in B$ so that $\varphi({zb} )= z^{\prime}b^{\prime}$.  If for all $z$ and $b$, there is an element $z^{\prime} \in \Bbb{Z}_l$ so that $\varphi (zb) = z^{\prime}b$, then the isomorphism is said to be a {\em translation}. 
\end{Definition} 
\begin{Proposition}
For $F$ a free,  finitely generated $\Bbb{Z}_l$-module  and a simply  transitive $G$-space $X$,  we have an isomorphism 
$$ \mbox{Aut} (F,X) \cong GL(F) \ltimes F
$$
For a based affine $\Bbb{Z}_l$-module $(F,X,B)$ we have 
$$\mbox{Aut}(F,X,B)\cong \Sigma _B \ltimes F
$$
where $\Sigma _B$ denotes the symmetric group of permutations of the set $B$.

\end{Proposition}
\begin{Proof} Clear from the definitions of the morphisms.  
\end{Proof} 

{\bf Remark:} The group $\mbox{Aut}(F,X,B)$ is a profinite group.  For every $n > 0$, we have the quotient of $X$ by the equivalence relation $\simeq _n$ given by 
$x \simeq _n  y$ if and only if there is an element $f \in l^n F$ so that $f\cdot x = y$.  These equivalence relations are respected by the action of $\mbox{Aut}(F,X,B)$, and the quotients are finite sets.  

Let $(F,X,B)$ be a based affine $\Bbb{Z}_l$-module, and let $\frak{A} = \mbox{Aut}(F,X,B)$.  Also, let $\frak{T} \subseteq \frak{A}$ denote the subgroups of translations.  We now  construct a family $\{ \rho _s \}$ of $k$-linear representations of $\frak{A}$.  Let $T^{\infty}$ denote the inverse limit $\column{lim}{\leftarrow} \mu _{l^n}$, and choose a continuous identification $\alpha :\Bbb{Z}_l \stackrel{\sim}{\rightarrow} T^{\infty}$. Let $\pi _s : T^{\infty} \rightarrow \mu_{l^s}$ be the projection.   Fixing an element $b \in B$,   we define a subgroup 
$\frak{A}^{\prime}  \subseteq \frak{A} $ to be $p^{-1}(\Sigma_B(b))$, where $p: \frak{A} \rightarrow \Sigma _B$ is the natural projection, and where $\Sigma_B(b) \subseteq \Sigma _B$ denotes the stabilizer of the element $b$.   We define a one dimensional character $\sigma _s $ of $\frak{A}^{\prime}$ by setting $\sigma _s(b)  = \pi _s  \alpha (1)$, $\sigma _s (b^{\prime}) = 1 $ for $b^{\prime} \neq b$, and $\sigma _s (\Sigma_B(b)) = \{ 1 \}$.   We define $\rho _s$ to be the induced representation
$$ \mbox{Ind}_{\frak{A}^{\prime} }^{\frak{A}} (\sigma _s)
$$
Of course, each $\rho _s$ now gives an action of $\frak{A}$ on the $n$-dimensional affine space over $k$ by algebraic automorphisms. We denote this variety with $\frak{A}$-action by $\Bbb{A}^n(k,\rho _s)$.  Note that the affine space also is equipped with a fixed choice of coordinates, since it is an induced representation.  The coordinates correspond to the cosets $\Sigma _B/\Sigma _B(b)$, which in turn correspond to the basis $B$ of $F$.  We define a self map $\theta: \Bbb{A}^n(k) \rightarrow \Bbb{A}^n (k)$ by requiring that 
$$ \theta (\sum  x_b b) = \sum  x_b^l b
$$
\begin{Proposition}
The map $\theta$ is $\frak{A}$-equivariant as a map from $\Bbb{A}^n(k, \rho _{s+1} )$ to $\Bbb{A}^n(k, \rho _{s} )$
\end{Proposition} 
\begin{Proof} Immediate.
\end{Proof}

We now consider the  pro-scheme  

$$  \cdots \stackrel{\theta}{ \rightarrow} \Bbb{A}^n(k,\rho _{s+1})  \stackrel{\theta}{ \rightarrow} \Bbb{A}^n(k,\rho _{s}) \stackrel{\theta}{ \rightarrow}  \Bbb{A}^n(k,\rho _{s-1}) \rightarrow \cdots \stackrel{\theta}{ \rightarrow} \Bbb{A}^n(k,\rho _{1})\stackrel{\theta}{ \rightarrow}  \Bbb{A}^n(k,\rho _{0})
$$
It is equipped with the action by $\frak{A}$, and we will denote its inverse limit by ${\cal A}_{(F,X,B)}$.  This inverse limit exists because the system consists of affine varieties, and is itself the affine scheme associated to the ring $\bigcup_n k[t^{\frac{1}{l^n}}]$. The $\frak{A}$-action is a continuous affine $\frak{A}$-action.   We also construct the related continuous affine $\frak{A}$-action ${\cal T}_{(F,X,B)}$, which is obtained by deleting the coordinate hyperplanes from each $\Bbb{A}^n(k,\rho _s)$ for each $s$.  We also define ${\cal A}^{s}_{(F,X,B)}$ and ${\cal T}^{s}_{(F,X,B)}$  to be the representation $\Bbb{A}^n(k,\rho _s)$ and the result of removing the coordinate hyperplanes from $\Bbb{A}^n(k,\rho _s)$, respectively.  

{\bf Remark:} The construction of ${\cal A}_{(F,X,B)}$ and ${\cal T}_{(F,X,B)}$   depends on the choice of $ b \in B$. It is easy to verify that varying this choice does not affect ${\cal A}_{(F,X,B)}$ or ${\cal T}_{(F,X,B)}$  up to isomorphism.

\begin{Proposition}\label{translationfree}
The restriction of the $\frak{A}$ action on ${\cal T}_{(F,X,B)}$ to $\frak{T}$ is an etale  $\frak{T}$-action.
\end{Proposition} 
\begin{Proof} Consider the direct system of rings 
$$ k[t_1^{\pm 1}, \ldots , t_n^{\pm 1}] \subseteq k[t_1^{\pm\frac{1}{l} }, \ldots , t_n^{\pm\frac{1}{l} }] \subseteq k[t_1^{\pm\frac{1}{l^2} }, \ldots , t_n ^{\pm\frac{1}{l ^2} }] \subseteq \cdots 
$$ 
 where $n$ is the rank of $F$, and let ${\cal L}_n$ denote its colimit.  This system is identified with the system of rings one obtains by applying the affine coordinate ring functor to the system defining  ${\cal T}_{(F,X,B)}$.  Further, ${\cal L}_n$ is an infinite Galois extension of $k[t_1^{\pm 1}, \ldots, t_n ^{\pm 1}]$, with Galois group $\Gamma = \Bbb{Z}_l^n$.  The restriction of the $\frak{A}$-action on ${\cal T}_{(F,X,B)}$ constructed above using the linear representations $\rho _s$, when restricted to $\frak{T}$, clearly provide an identification of  $\frak{T} \cong F$ with $\Gamma$, and the result now follows.  
\end{Proof}

\begin{Definition}
For any pro-$l$ group $G$, an {\em affine $l$-adic representation} (or $l$-BAR) is a continuous homomorphism from $G$ to the automorphism group of a based affine $\Bbb{Z}_l$-module.  For any affine $l$-adic representation $ \eta$ of $G$ on an affine $\Bbb{Z}_l$ module $(F,X,B)$, we obtain actions of $G$ on the schemes ${\cal A}_{(F,X,B)}$ and ${\cal T}_{(F,X,B)}$ by pulling back the action of the automorphism group of $(F,X,B)$ on these schemes.  We will denote these $G$-schemes by ${\cal A}_{\eta}$ and ${\cal T}_{\eta}$.  We also write ${\cal A}^s_{\eta}$ and ${\cal T}_{\eta} ^s$ for the pullback along $\eta$ of the actions ${\cal A}_{(F,X,B)}^s$ and ${\cal T}_{(F,X,B)}^s$. 
\end{Definition} 

There is a useful description of ${\cal T}_{\eta}^s$ which we will need. 

\begin{Proposition}  \label{usefulprop} The scheme ${\cal T}_{\eta}^s$ can be identified with the quotient ${\cal T}_{\eta}/l^s \cdot \frak{T}$, where $\frak{T}$ denotes the full translation subgroup of $\frak{A}$.    $G$ acts on this quotient, since the image of $G$ under $\eta$ normalizes all the subgroups $l^s\cdot \frak{T}$.  
\end{Proposition}
\begin{Proof}  This is just the ring theoretic result that 
$$  \bigcup _n k[t_1^{\pm \frac{1}{l^n}}, \ldots , t_k^{\pm \frac{1}{l^n}}]^{l^s\cdot \frak{T}} = k[t_1^{\pm \frac{1}{l^s}}, \ldots , t_k^{\pm \frac{1}{l^s}}]
$$
\end{Proof} 

\begin{Proposition}  Let $\eta$ be any $l$-BAR.  Then the continuous affine $G$-actions ${\cal A}_{\eta} \times Spec(L) $   and   ${\cal T}_{\eta} \times Spec(L)$ are both coherent, where $L$ is any field equipped with a continuous action by $G$.  As always, it is assumed that $L$ is equipped with the discrete topology.    
\end{Proposition}
\begin{Proof} It suffices to prove that the rings 
 $$ \bigcup _n L[ x_1^{ \frac{1}{l^n}}, \ldots , x_s^{ \frac{1}{l^n}}]$$ and $$ \bigcup _n L[ x_1^{\pm \frac{1}{l^n}}, \ldots , x_s^{\pm \frac{1}{l^n}}]$$ are coherent rings for any field $L$.  This follows immediately from Proposition \ref{coherent} since the inclusions $$ L[ x_1^{ \frac{1}{l^n}}, \ldots , x_s^{\frac{1}{l^n}}] \hookrightarrow L[ x_1^{\frac{1}{l^{n+1}}}, \ldots , x_s^{ \frac{1}{l^{n+1}}}]
$$  and $$ L[ x_1^{\pm \frac{1}{l^n}}, \ldots , x_s^{\pm \frac{1}{l^n}}] \hookrightarrow L[ x_1^{\pm \frac{1}{l^{n+1}}}, \ldots , x_s^{\pm \frac{1}{l^{n+1}}}]
$$ are all flat.  \end{Proof}

 Given a possibly infinite set of $l$-BAR's $\{ \eta _i\}$, the schemes $\prod {\cal A}_{\eta _i} $ and $\prod {\cal T}_{\eta _i} $ both are also continuous affine $G$-schemes. 
 \begin{Corollary}  Both $\prod {\cal A}_{\eta _i} $ and $\prod {\cal T}_{\eta _i} $ are coherent. 
\end{Corollary}
\begin{Proof}  It will suffice to prove that the rings  $$ \bigcup _n L[ x_1^{ \frac{1}{l^n}}, \ldots , x_s^{ \frac{1}{l^n}}, \ldots ]$$ and $$ \bigcup _n L[ x_1^{\pm \frac{1}{l^n}}, \ldots , x_s^{\pm \frac{1}{l^n}},\ldots ]$$ are coherent.  This follows easily from Proposition \ref{coherent}  since the inclusions 
$$ \bigcup_n L[ x_1^{ \frac{1}{l^n}}, \ldots , x_{s}^{\frac{1}{l^n}}] \hookrightarrow \bigcup _nL[ x_1^{\frac{1}{l^{n}}}, \ldots , x_{s+1}^{ \frac{1}{l^{n}}}]
$$  and $$ \bigcup _n L[ x_1^{\pm \frac{1}{l^n}}, \ldots , x_s^{\pm \frac{1}{l^n}}] \hookrightarrow \bigcup_nL[ x_1^{\pm \frac{1}{l^{n}}}, \ldots , x_{s+1}^{\pm \frac{1}{l^{n}}}]
$$ are all flat. 
\end{Proof} 

We now have the following definition. 

\begin{Definition} Let $G$ be a profinite group, and suppose $K$ is any closed subgroup of finite index.  We say an affine $l$-adic representation $\eta$  of $G$  is {\em $K$-principal} if the stabilizer in $G$ of any $x \in X_{\eta}$ is equal to $K$. Note that it follows that $K$ is normal.   Clearly, the condition is equivalent to the condition that $X$ is a free $G/K$-set.     \end{Definition}

The following lemma is an $l$-adic analogue of the results of Auslander and Kuranishi concerning the structure of fundamental groups of flat manifolds \cite{auslander}.

\begin{Lemma} \label{torsion free} Let $K \subseteq G$ be a closed normal subgroup of a profinite group $G$.  Then a  based affine $l$-adic representation $\eta$ of a $G$ group is $K$-principal if and only if the image of $\eta$ in the automorphism group of $X_{\eta}$ is torsion free.  
\end{Lemma} 
\begin{Proof}  We let $F$ be a free finitely generated $\Bbb{Z}_l$-module with basis $B$, so that $$(F_{\eta}, X_{\eta}, B_{\eta} ) \cong (F,F,B)$$ where $F$ acts on itself by translation.  The automorphism group of $(F,F,B)$ is a semidirect product $\Sigma_B \ltimes F$.  The group $\Sigma _B$ acts on $F$ by permuting the basis $B$, and we describe the action of an element $(\sigma, f)$, where $\sigma \in \Sigma _B$ and $f\in F$ by $(\sigma , f) \cdot f^{\prime} = \sigma \cdot f^{\prime} + f$.  Given an element $(\sigma , f) $ of $\Sigma _B \ltimes F$, the action of $(\sigma , f)$ on $F$ has a fixed point if and only if 
$f = v - \sigma v$ for some $v \in F$.  This condition is in turn equivalent to the requirement that any $\Bbb{Z}_l$-linear function which satisfies $f\cdot \sigma = f$ vanishes on $f$.  Functions satisfying this condition form a free $\Bbb{Z}_l$-module of rank equal to the number of cycles in the cycle decomposition of the permutation $\sigma$, with the characteristic functions on orbits forming a basis.  On the other hand, $(\sigma , f)$ has finite order $n$ if and only if $\sigma ^n = e$ and 
$$f+ \sigma\cdot f + \sigma ^2 \cdot f + \cdots + \sigma ^{n-2} \cdot f + \sigma ^{n-1} \cdot f = 0
$$
Let $\{ \frak{o}_i \} _i$ denote the collections of orbits under action of $\sigma$, and let $\pi _i$ denote the projection on the set of basis elements in $\frak{o}_i$.  The summands corresponding to the orbits are invariant under $\sigma$, and we obtain a $\sigma$-invariant direct sum decomposition
$$ F \cong \bigoplus _i F_i
$$
where $F_i$ is the span of the basis elements belonging to $\frak{o}_i$.   The above equation now demonstrates that 
$$ \sum_{s = 0}^{n-1} \sigma^s (\pi _i (f)) = 0
$$
It follows that all the $\sigma$-invariant $\Bbb{Z}_l$-linear functions on $F$ vanish on $f$.  So we have proved that for any element $(\sigma , f) \in Aut(F,F,B)$, the condition that $(\sigma , f)$ as a fixed point under its action on $F$ is equivalent to its having finite order.  The result now follows. 
\end{Proof} 

We will  write ${\cal T}_{\eta}/G$ for the orbit scheme of ${\cal T}_{\eta}$.  

We now have the following result. 
\begin{Proposition}\label{translation}  Let $G$ be a profinite group, with $N \subseteq G$ a closed normal subgroup.  Let $\eta$ denote a based affine $l$-adic representation of $G$.  If  $\eta$ is $N$-principal, then ${\cal T}_{\eta}$ is an etale  $G/N$-scheme. \end{Proposition}
\begin{Proof} Let $(F_{\eta}, X_{\eta}, B_{\eta})$ be the based $l$-adic representation space, so $\eta$ is a continuous homomorphism $G \rightarrow \frak{A}$, where $\frak{A}$ denotes the automorphism group of  $(F_{\eta}, X_{\eta}, B_{\eta})$.  As before, let $\frak{T}$ denote the translation subgroup of $\frak{A}$.  The subgroup $\eta ^{-1}(\frak{T})$ is now a closed subgroup of finite index in $G$, containing $N$, and therefore projects to a closed subgroup $T$ of finite index in $G/K$.  Let $\frak{X}$ denote the set of all morphisms from schemes of the form $Spec(F)$, where $F$ is an algebraically closed field, to the inverse limit ${\cal T}_{\eta}$.  In order to prove the freeness of the $G/N$-scheme, we must prove the freeness of he $G/N$ action on the set $\frak{X}$.  We know from Proposition \ref{translationfree} that the restriction of the action on $\frak{X}$ to the subgroup of finite index $T$ is free, and that the group $G/N$ is torsion free.  In order to prove the freeness of the full $G/N$-action, we suppose that there is are elements $\gamma \in G/N$, with $\gamma \neq e$,  and $x \in \frak{X}$, so that   $\gamma \cdot x = x$.  It would follow that for any $n$, the element $\gamma ^n$ would also fix $x$.  But for sufficiently large $n$, $\gamma ^n \in T$, and further that $\gamma^n \neq e$, because of the torsion freeness proved in Proposition \ref{torsion free}.  The freeness of the $N$ -action now contradicts the existence of $\gamma$.  
\end{Proof} 

We actually need a slightly stronger result, namely that for a fixed $N$-principal $l$-BAR, there is an $s \geq 0$ so that the scheme ${\cal T}^s_{\eta}$ is acted on freely by the quotient group $G/N\cdot \eta ^{-1}(l^{s\pr} \frak{T})$ for all $s^{\prime} \geq s$.  This will require a few algebraic facts.  

Consider any extension of groups 
$$  0 \rightarrow A \rightarrow H \rightarrow Q \rightarrow 0
$$
where $Q$ is a finite cyclic group of order $l$, and $A$ is abelian.  The extension gives an action of $Q$ on $A$ by automorphisms.  Let $N(A) \subseteq A$ denote the subgroup of elements of the form $ a + \sigma \cdot a + \sigma^2 \cdot a + \cdots \sigma^{p-1}a$, where $\sigma \in Q$ is a generator.  Of course the group is independent of the choice of generator.   We define a function $c:Q \rightarrow A/N(A)$  as follows.  For any element $\sigma \in Q$, select a lift $\overline{\sigma} \in H$, i.e. so that $\pi \overline{\sigma} = \sigma$, and define $c(\sigma) = \overline{\sigma}^l$.  It is easy to check that $c(\sigma)$ is independent of the choice of $\overline{\sigma}$.  

\begin{Lemma} \label{smallfact} An element $\sigma \in Q$ lifts to an element of order $l$ if and only if $c(\sigma) = 0$.  
\end{Lemma} 
\begin{Proof}  The ``only if" part is clear, since if $\sigma$ lifts to an element of order $l$, $c(\sigma)$ is by definition equal to zero.  Suppose on the other hand that $c(\sigma) = 0$, which means that there is a lift $\overline{\sigma}$ and an element $\alpha \in A$ so that 
$$\overline{\sigma} ^l = \alpha + \sigma \alpha + \cdots \sigma ^{l-1} \alpha
$$
If we now set $\overline{\sigma}\pr = \overline{\sigma} \cdot (- \alpha)$, then $(\overline{\sigma} \pr )^l = e$. 
\end{Proof} 

\begin{Lemma} \label{simplefact}  Let $\Gamma$ denote a pro-$l$ group which fits in an exact sequence 
$$ \{ e \}\rightarrow L  \rightarrow \Gamma \stackrel{\pi}{\rightarrow} Q \rightarrow \{ e \}
$$
where $Q$ is a finite $l$-group, and $L \cong \Bbb{Z}_l^n$.  Let $\Gamma \pr \subseteq \Gamma$ be a subgroup so that the composite 
$$  \Gamma \pr \rightarrow \Gamma \rightarrow Q
$$
is surjective.  Suppose further that $\Gamma \pr$ is torsion free, and that we are given a continuous $\Gamma$-set $X$ on which $\Gamma \pr$ and $L$ both act freely.  Then there is an $s \geq 0$ so that $\Gamma \pr / \Gamma \pr \cap l^{s\pr}\cdot L$ acts freely on the set $X/l^{s\pr}\cdot L$ for all $s \pr \geq s$.   The number $s$ depends only on the groups involved, not the $\Gamma$-set $X$.  
\end{Lemma} 
\begin{Proof}  We first claim that it suffices to prove the result for $Q$ a cyclic finite $l$-group.  To see this, we suppose we are able to prove the result for each subgroup  $ \Gamma\pr_C = (\pi | \Gamma\pr)^{-1}(C) \subseteq \Gamma \pr$.  This means that for every cyclic subgroup $C \subseteq Q$, there is an $s(C) \geq 0$ so that $\Gamma_C\pr / \Gamma_C \pr \cap l^{s\pr}L$ acts freely on $X/l^{s\pr}L$ for all $s\pr \geq s(C)$ .  Since there are only finitely many cyclic subgroups of $Q$, we may set $S = \column{max}{C} s(C)$.  It now follows that $\Gamma \pr /\Gamma \pr \cap l^{s \pr}L$ acts freely on $X/l^{s \pr}L$ whenever $s \pr \geq S$,  since all the groups 
$\Gamma_C\pr / \Gamma_C \pr \cap l^{s\pr}L$ do, and since these the union of these subgroups is all of $\Gamma\pr / \Gamma \pr \cap l^{s \pr} L$.  From this point on, we assume that $Q$ is cyclic of prime power order.  The next reduction is that it will suffice to prove the result for $Q$ a cyclic group of order $l$.  To see this,  we first assume that the result holds for all groups $\Gamma\pr$ with $Q$ of order $l$, and  let $\Gamma ^{\prime \prime}\subseteq \Gamma \pr $ denote the subgroup $\pi ^{-1} (Q_l)$, where $Q_l$ is the unique subgroup of order $l$.  Suppose $\gamma \in \Gamma \pr - \Gamma^{\prime \prime}$, and let $l^k$ denote the order $\pi (\gamma) \in Q$.  We then have that $\gamma^{l^{k-1}} \in Q_l$, and $\gamma^{l^{k-1}} \neq e$.  By our assumption, there is an $s \geq 0$ so that the group $\Gamma^{\prime \prime} /  \Gamma^{\prime \prime} \cap l^s L$ acts freely on $X / l^{s \pr} \cdot LX$ for all $s \pr \geq s$.  This means that $\gamma$ also acts freely, since a power of it acts freely.  This give the reduction to the case where $Q$ has order $l$. To handle this case, we note that we have the function $c:Q \rightarrow L \cap \Gamma \pr$ defined above.   By Lemma \ref{smallfact} above, the torsion freeness of $\Gamma^{\prime} $ means that $c(\sigma) \neq 0$ in the group $L \cap \Gamma \pr/N(L) \cap \Gamma \pr $.  But it is clear that $L\cap \Gamma \pr /N(L) \cap \Gamma \pr$ is identified with the inverse limit  
$$  \column{lim}{\column{\longleftarrow}{s}} L \cap \Gamma \pr/(N(L) \cap \Gamma \pr + l^s L \cap \Gamma \pr)
$$ 
so  for sufficiently large $s$, we have $\overline{\sigma}^l \neq 0 $ in $L \cap \Gamma \pr/ l^s L \cap \Gamma \pr$, for every possible lift $\overline{\sigma}$ of $\sigma$.  Since the $L$ action on $X$ is assumed to be free,  it follows  that $L \cap \Gamma \pr/ l^s L \cap \Gamma \pr$ acts freely on $X /l^s\cdot  L$ for sufficiently large $s$.  We now have that every lift of $\sigma$ acts freely for sufficiently large $s$.  It is clear also that the same calculation applies to all non-identity powers  of $\sigma$, which together with the fact that $L/l^s\cdot L $ acts freely on $X /l^s\cdot  L$, as always for sufficiently large $s$,  gives the result. 
\end{Proof} 
\begin{Proposition}  \label{finitestageetale} Let $ \eta$ be a $N$-principal $l$-BAR for some closed normal subgroup of $G$.  Then there exists an $s \geq 0$ so that for every $s \pr \geq s$, the action of the finite group   $G /(N \cdot \eta^{-1}(l^{s\pr}\frak{T}))$  on ${\cal T}^{s\pr}_{\eta}$ is etale.  
\end{Proposition} 
\begin{Proof} The $G/N$-action on  ${\cal T}_{\eta}$ is etale.  Let $\Omega$ denote an algebraically closed field and  let $i: Spec \Omega \rightarrow {\cal T}_{\eta}/G$ be a morphism.  Then we have the set of all maps $\Omega \rightarrow {\cal T}_{\eta}$ for which the diagram 
$$  \begin{diagram}
 \node{Spec(\Omega)} \arrow{e} \arrow{se,t}{i} \node{{\cal T}_{\eta}} \arrow{s} \\
\node{}  \node{{\cal T}_{\eta}/G}
\end{diagram}
$$
commutes, which we have denoted by $\frak{H}(A[{\cal T}_{\eta}], \Omega, i)$.  The etaleness of the $G/N$-action on ${\cal T}_{\eta}$ means that the $G/N$-action on $\frak{H}(A[{\cal T}_{\eta}], \Omega, i)$ is free.  
   The group $G/N$ is identified by the representation with   a subgroup of 
$\Sigma_n(l) \ltimes \Bbb{Z}_l^n$, where $\Sigma_n(l)$ denotes an $l$-Sylow subgroups,and is torsion free, so satisfies the hypotheses of Lemma \ref{simplefact} above.  Note that the hypothesis that $L$ acts freely in Lemma \ref{simplefact} is satisfied due to Proposition \ref{translationfree}.  The group $\Gamma$ will be $I \ltimes \Bbb{Z}_l^n$, where $I$ is the image of the composite $$G/K \rightarrow \Sigma_n(l) \times \Bbb{Z}_l^n \rightarrow \Sigma_n (l)$$
 and the affine representation $\eta$ identifies $G/N$ with a subgroup of $\Gamma$.  It now follows that there is an $s$ so that for all $s \pr \geq s$, the group $G /(N \cdot \eta^{-1}(l^{s\pr}\frak{T}))$ acts freely on the orbit set $$\frak{H}(A[{\cal T}_{\eta}], \Omega, i)/ (l^{s\pr}\frak{T})$$ But, because of the description of ${\cal T}_{\eta}^s$ given in Proposition \ref{usefulprop},  this orbit set is just the set $\frak{H}(A[{\cal T}_{\eta}^s], \Omega, i)$,  which gives the result.   Note that we are using the self evident fact that if a group $G$ acts freely on a set $X$, and $N$ is a normal subgroup of $G$, then the naturally defined action of $G/N$ on $X/N$ is also free.  
\end{Proof}

  We record a useful property of $K$-principal actions.  

\begin{Proposition}\label{principalcriterion}  Let $\eta _1$ and $\eta _2$ be affine $l$-adic representations of a profinite group $G$. Then we  have the evident product affine $l$-adic representation $\eta _1 \times \eta _2$.  Let $\{e\} \subseteq K^{\prime} \subseteq K \subseteq G$ be closed subgroups.  Suppose that $\eta _1$ is $K$-principal, and that $\eta_2|K$ is $K^{\prime}$-principal.  Then $\eta _1 \times \eta _2$ is $K^{\prime}$-principal.  In particular, if $\eta _i$ is $K_i$-principal, for $K_i \subseteq G$, then $\eta _1 \times \eta _2$ is $K_1 \cap K_2$-principal.    \end{Proposition} 
\begin{Proof} Clear.  
\end{Proof}

We also consider the notion of induction of affine based  $l$-adic representations.  Let $\eta = (F_{\eta}, X_{\eta}, B_{\eta})$ be a based affine $l$-adic representation of $K$, where $K$ is a closed subgroup of finite index in a profinite group $G$.  Then we define an induced based affine $l$-adic representation $i_K^G(\eta)
= (F_{i_K^G(\eta)}, X_{i_K^G(\eta)}, B_{i_K^G(\eta)})$ as follows.   The sets $F_{i_K^G(\eta)}$ and $X_{i_K^G(\eta)}$ are defined as the sets of equivariant functions $F^K(G, F_{\eta})$ and $F^K(G, X_{\eta})$, respectively, where $K$ acts by left multiplication on $G$ and via $\eta$ on $F_{\eta}$ and $X_{\eta}$.  The group $G$ acts by right multiplication on $G$, yielding a left action on the two sets.  The free action of $F^K(G, F_{\eta})$ on  $F^K(G, X_{\eta})$ is given by pointwise multiplication.  To define $B_{i_K^G(\eta)}$, we first define elements $f_b \in F_{i_K^G(\eta)}$ for each $b \in B_{\eta}$ by setting $f_b(e) = b$.  The remainder of $B_{i_K^G(\eta)}$  is obtained by translating the elements using the $G$ action.  It is clear that this set will give a basis for $F_{i_K^G(\eta)}$.  

We say an affine based $l$-adic representation $\eta$ of a profinite group $G$ is a {\em translation representation} if for every $g \in G$, the automorphism $\eta (g)$ is a translation.

\begin{Proposition}  \label{induction} Let $G$ be a  profinite group, and let $K \lhd G$ be a normal subgroup of finite index.  Suppose we are given a one-dimensional based affine $l$-adic representation $\eta$ of $K$.  Then we have the following. 

\begin{enumerate}
\item{The restriction of $i_K^G(\eta)$ to $K$ is a translation representation. }
\item{The restriction of $i_K^G(\eta)$ to $K$ is $K^{\prime}$-principal, where $K^{\prime}$ is defined as follows.  For each element $\gamma \in G/K$, we denote by $\eta ^{\gamma}$ the result of composing $\eta$ with the conjugation action of a coset representative of $\gamma$.  Then $K^{\prime}$ is given by 

$$ K^{\prime} = \bigcap _{\gamma  \in G/K} Ker(\eta ^{\gamma} )
$$
}
\end{enumerate} 
\end{Proposition} 

\section{Construction of ${\cal E}G$} \label{universal} 

In this section, we will use products of the actions constructed in the last section to build pro-objects in the category of continuous $G$-actions of a profinite $l$-group which will satisfy the requirements formulated in the introduction.   Throughout this section, we will allow ourselves to refer to a based affine $l$-adic representation as an $l$-BAR, in the interest of brevity. 
We begin with  a technical lemma. 

\begin{Lemma} \label{principalfinal} Let $G$ be a totally torsion free pro-$l$-group. Let $G_0 \subseteq G$ be any closed subgroup of finite index.  Then there exists a closed normal subgroup $K \subseteq G_0$ and a $K$-principal 
$l$-BAR  $\eta$ of $G$.  \end{Lemma} 
\begin{Proof}  
We first  observe that it suffices to prove the result for $G_0$ normal in $G$, since there is a always a normal subgroup of finite index contained in $G_0$.  
 
 To construct an $l$-BAR   which is $K$-principal for some $K \subseteq G_0$, we first note that the set of subgroups of $G$ containing $G_0$ consists entirely of closed subgroups, by Proposition \ref{profiniteprops}.  We now   assume by induction that we are able to  construct such $l$-BARs  for all closed subgroups of index less than the index of $G_0$ in $G$.  We  construct such $l$-BAR's for every normal subgroup containing $G_0$ properly.  There are obviously only finitely many such subgroups.    We consider the finite list $\{ N_i \}_{i= 1}^s$ of all normal subgroups of $G$ which contain $G_0$ properly (they are necessarily closed), and for each one we construct a $K_i$-principal $l$-BAR $\eta_i$ of $G$, where $K_i \subseteq N_i$.   The $l$-BAR  $\eta_1 \times \eta _2 \times \cdots \times \eta _s$ is now $K_1 \cap \cdots \cap K_s$-principal by repeated application of Proposition \ref{principalcriterion}.  We now consider the subgroup $G_0  \cdot (K_1 \cap  \cdots \cap  K_s) $.   This subgroup contains $G_0 $, and therefore corresponds to a subgroup of the quotient $G/G_0$, contained in $(N_1/G_0 ) \cap (N_2/G_0)  \cap \cdots \cap (N_s/G_0) $.  This group can be characterized as the intersection $\sigma  (G/N)$ of all non-trivial normal subgroups in $G/G_0$.  Since $G_0 $ is a closed normal subgroup of the pro-$l$ group $G$, we know that $G/G_0 $ is an $l$-group, and hence its center $Z(G/G_0 )$ is non-trivial. Therefore, $\sigma (G/G_0 ) \subseteq Z(G/G_0 )$ is a central subgroup, which implies that every subgroup of $\sigma (G/G_0 )$ is normal in $G/G_0 $.  By the definition of $\sigma (G/G_0 )$, this means that $\sigma (G/G_0 )$ has no non-trivial proper subgroups, and the only such $l$-groups are $\{ e \}$ and $\Bbb{Z}/l \Bbb{Z}$. 
 Consequently, either 
$G_0  \cdot (K_1 \cap  \cdots \cap  K_s) /G_0$ is trivial, in which case $K_1 \cap  \cdots \cap  K_s \subseteq G_0$, and we are done, or it is cyclic of order $l$.  In this case, we are given a homomorphism from $ \overline{G}_0  = G_0  \cdot(K_1 \cap \cdots \cap K_s)$ to $\Bbb{Z}/l\Bbb{Z}$, and therefore a character $\chi _1$ of   
$\overline{G}_0 $.  Since the group $G$ is totally torsion free, we find that the group  $\overline{G}_0 ^{ab}$ is torsion free, and therefore that the character group is divisible.  Construct a sequence of characters $\chi _i$ so that $\chi _i ^l = \chi _{i-1}$, with $\chi _1$ the previously defined character.  This sequence constructs a one-dimensional $l$-BAR   of $\overline{G}_0$, which we now induce up to an $l$-BAR  ${\eta}$ of dimension equal to the index of $G_0 $ in $G$.  We claim that the $l$-BAR  $$\eta \times \eta _1 \times \cdots \times \eta _s$$ is $Ker(\eta) \cap (K_1 \cap K_2 \cap \cdots \cap K_s)$-principal.  To see this we note that  the stabilizer of any point in  $X_{\eta _1} \times \cdots \times X_{\eta _s}$ is equal to  $K_1 \cap \cdots \cap K_s \subseteq \overline{G}_0$, and therefore by Proposition \ref{principalcriterion} it suffices to prove that the restriction of the $G$-action on $X_{\eta}$ to $K_1 \cap \cdots \cap K_s$ is $Ker(\eta) \cap K_1 \cap \cdots \cap K_s$-principal.   Proposition \ref{induction} asserts that the restriction to $G_0$ of $\eta$ is $Ker(\eta)$-principal. The result now follows from Proposition \ref{principalcriterion}.\end{Proof} 

We need a second fact concerning representations of a profinite $l$-group.  

\begin{Proposition} \label{repall}   Let $\sigma$ denote any irreducible continuous representation of a totally torsion free profinite $l$-group $G$, over an algebraically closed field of characteristic prime to $l$.    Then $\sigma$ is isomorphic to ${\cal T}^s_{\eta}$  for some $l$-$BAR$ $\eta$ and some $s$. 
\end{Proposition}
\begin{Proof}  The continuity of the $\sigma$ gives that $\sigma$ factors through a finite quotient $G/N$ of $G$, where $N$ is a closed normal subgroup of finite index in $G$.  By Blichfeldt's theorem, $\sigma $ (regarded as a representation of $G/N$)  is monomial, and is therefore of the form $i_K^G(\chi )$ for some closed subgroup of finite index $K \subseteq G/N$, and some one-dimensional character $\chi$ of $K$.  Let $\overline{K}$ denote the inverse image of $K$ in $G$. The group $K$ is also totally torsion free, and therefore any homomorphism from $K$ to a finite cyclic $l$-group can be lifted to a continuous homomorphism $\hat{\chi }$  from $K$  to $\Bbb{Z}_l$.  This homomorphism is a one-dimensional $l$-BAR, and it is clear from the construction that $\chi$ is equal to $\rho _s \compcirc \hat{\chi }$.  It now readily follows that the induced (from $\overline{K}$ to $G$)  $l$-BAR is the required $\eta$.  \end{Proof} 

We now construct the continuous affine $G$-scheme  ${\cal E}G $.
\begin{Definition} \label{cof}  Let $G$ be any profinite group, which is acting continuously on an affine scheme $X = Spec(A)$.  We say $X$ is {\em cofinite} if there is a closed normal  subgroup $N \subseteq G$  of finite index so that the $G$ action on $X$ factors through the quotient $G/N$.  For any cofinite affine $G$-scheme $X$, we define the subgroup ${\cal N}(X) \subseteq G$ to be the maximal closed normal subgroup so that the $G$-action on $X$ factors through $G/{\cal N} (X)$.  It is clear that this defines ${\cal N}(X)$ uniquely, and moreover that if we have a map $X \rightarrow Y$ of cofinite affine $G$-schemes, then ${\cal N}(X) \subseteq {\cal N} (Y)$.  
We say a cofinite affine $G$-scheme $X$ is etale if the group action of $G/{\cal N}(X)$ on $X$ is a free action in the sense of Proposition \ref{free}. We  define $\frak{C}^{cof}(G)$ to be the category whose objects are coherent  cofinite affine $G$-schemes, and whose morphisms are the faithfully flat equivariant maps of such $G$-schemes.   The full subcategory of etale actions will be denoted by $\frak{C}^{cof}_{et}(G)$.  \end{Definition}  
\begin{Proposition} \label{sorite} Suppose that we have a $G$-equivariant map $f: X \rightarrow Y$ of cofinite etale affine $G$-schemes, which is faithfully flat.  Then for any closed subgroup  $K$ of finite index in $G$, the natural map $X/K \rightarrow Y/K$ is also faithfully flat.  
\end{Proposition}
\begin{Proof} Since the etale property of an action  clearly descends to closed subgroups, it suffices to prove this for the case $K = G$.  We now have the diagram 
$$
\begin{diagram}
\node{X} \arrow{e,t}{f}  \arrow{s,t} {\pi _X}  \node{Y} \arrow{s,b}{\pi _Y} \\
\node{X/G} \arrow{e,t}{f/G} \node{Y/G}
\end{diagram} 
$$
The morphism  $\pi _Y \compcirc f$ is faithfully flat, since $f$ is and $\pi _Y$ is since it is etale and surjective. $\pi _X$ is faithfully is also faithfully flat by the same reasoning.  It follows that $f/G$ is faithfully flat.  
\end{Proof}

 For each closed normal subgroup of finite index $N \lhd G$ we apply Lemma \ref{principalfinal} to construct an $l$-BAR $\eta _N$ which is $N^{\prime} $-principal for some closed normal subgroup $N^{\prime}\subseteq N$.  Let $\frak{N}^G$ denote the  partially ordered set of closed normal subgroups of $G$.  For each $(s,N) \in \Bbb{N} \times \frak{N}^G$, we let $\theta^{\prime } (s,N)$ denote the (finite) product 
$$  \theta^{\prime}  (s,N)  = \prod _{\overline{N} \supseteq N, \overline{N} \in \frak{N} } {\cal T}_{\eta _{\overline{N}}}^s
$$
The construction $\theta\pr$ clearly gives a functor from the category $\Bbb{N} \times \frak{N}^G$ to $\frak{C}^{cof}(G)$, where the effect on morphisms consists on  bonding maps ${\cal T}_{\eta}^{s+1} \rightarrow {\cal T}_{\eta}^s$ together with product projections. We will need to build a more complicated pro-object.  Let $\frak{R}^G$ denote the set of all isomorphism classes of finite dimensional continuous $k$-linear representations of $G$.  For each $N \in \frak{N}^G$, we let $\frak{R}^G(N) \subseteq \frak{R}^G$ denote the subset of all  representations which factor through the quotient $G/N$.  For each $\alpha \in \frak{R}^G$, we select an $l$-BAR $\eta _{\alpha}$ so that for some $s$, $\rho _s \compcirc \eta _{\alpha}$ belongs to the isomorphism class $\alpha$.  This is possible because of Proposition \ref{repall}.  We define 
$$ \theta ^{\prime \prime }(s,N) = \prod _{\alpha \in \frak{R}^G(N)} {\cal T}_{\eta _{\alpha}}^s
$$
$\theta^{\prime\prime} $ defines a functor from $\Bbb{N} \times \frak{N}^G$ to $\frak{C}^{cof}(G)$, again with structure maps arising from those in the pro objects ${\cal T}_{\eta}$ and product projections.    We let ${\cal E}G$ denote the inverse limit of   $\theta = \theta^{\prime}  \times \theta ^{\prime \prime}$  in the category of  affine $G$-schemes.  ${\cal E}G$ is thus an affine scheme whose affine coordinate ring is the colimit of the rings $A(\theta (s,N))$ over the category $\Bbb{N} \times \frak{N}^G$, where $A(X)$ denotes the affine coordinate ring of $X$.  

\begin{Proposition}  For any pair $N \subseteq N \pr \subseteq G$ of closed normal subgroups of finite index, and any $s \leq s \pr$, the map of schemes $ \theta (s \pr , N \pr)  \rightarrow \theta (s,N)$ is faithfully flat.  It follows that ${\cal E}G$ is coherent. 
\end{Proposition} 
\begin{Proof} On affine coordinate rings, each map is of the form 
$$ A = k[x_1^{\pm \frac{1}{l^s}}, \ldots, x_t^{\pm \frac{1}{l^s}}] \hookrightarrow k[x_1^{\pm \frac{1}{l^{s\pr}}}, \ldots, x_{t \pr}^{\pm \frac{1}{l^{s \pr}}}] = B
$$
where $t \pr \geq t$.  In this case, it is an easy computation that $B$ is a free module over $A$.  Both rings are Noetherian, therefore coherent, and the result now follows from Proposition \ref{coherent}. 
\end{Proof} 

\begin{Proposition}  The $G$-action on the affine coordinate ring $A[{\cal E}G]$ is continuous, when $A[{\cal E}G]$ is equipped with the discrete topology. 
\end{Proposition}
\begin{Proof}  Each ring $A[\theta(s,N)]$ is acted on by $G$ via an action which factors through a finite quotient of $G$.  Since $A[{\cal E}G]$ is the colimit of the rings $A[\theta(s,N)]$, the result follows. 
\end{Proof} 

We will now prove that the $G$-action on ${\cal E}G$ is etale.  Let $\overline{\theta}$ denote the functor $(s,N) \rightarrow \theta(s,N)/N$.  $\overline{\theta}(s,N)$ is equipped with a $G/N$-action.  
\begin{Lemma} The natural map $\column{lim}{\leftarrow} \theta \rightarrow \column{lim}{\leftarrow} \overline{\theta}$ is an isomorphism of affine schemes.
\end{Lemma} 
\begin{Proof}  We prove that the map induces an isomorphism on affine coordinate rings.   The map $A[\column{lim}{\leftarrow} \overline{\theta}] \rightarrow A[\column{lim}{\leftarrow} \theta ]$ is clearly injective, since for every $(s,N)$, the map $A[\overline{\theta}(s,N)] \rightarrow A[ \theta(s,N)]$ can be identified with the inclusion $A[ \theta(s,N)]^N \rightarrow A[ \theta(s,N)]$.   Given an element $x \in A[\column{lim}{\leftarrow} \theta]$, it is in the image of $A[\theta(s,N)]$ for some $s$ and $N$.  This means that it is fixed by a normal subgroup of finite index $N \pr \subseteq N$, and consequently that it is in the image of $A[\theta(s,N\pr)]^{N \pr} = A[\overline{\theta}(s,N \pr)]$.  This proves surjectivity.  
\end{Proof}

\begin{Lemma}  Let $N \subseteq G$ be a closed subgroup of finite index.  Then $\column{lim}{\column{\leftarrow}{s}} \overline{\theta}(s,N)$ an etale $G/N$-action for all $N$.  
\end{Lemma} 
\begin{Proof}  Proposition \ref{translation} shows that $\column{lim}{\column{\leftarrow}{s}}{\theta}(s,N)$ is an etale action of $G/N\pr$ for some closed normal subgroup $N \pr \subseteq N$.  It follows that the action of $G/N$ on 
$$ \column{lim}{\column{\leftarrow}{s}} {\theta}(s,N)/N = \column{lim}{\column{\leftarrow}{s}} \overline{\theta}(s,N)
$$
is etale.  
\end{Proof} 
\begin{Corollary} \label{etalequantitative}   For every closed normal subgroup $N \subseteq G$ of finite index,  the $G/N$-action on $\overline{\theta}(s \pr, N)$ is etale for sufficiently large $s$.  
\end{Corollary}
\begin{Proof}  We let $N \pr \subseteq N \subseteq G$ be the closed normal subgroup of $G$ so that $\eta _N$ is $N \pr$-principal, and therefore that  $\column{lim}{\column{\leftarrow}{s}}{\theta}(s,N)$ is an etale $G/N\pr$-action.    It follows from Proposition \ref{finitestageetale} that there is an $s \geq 0$ so that the action of $G/  \eta ^{-1}( l^{s \pr}\frak{T})$ on $ \column{lim}{\column{\leftarrow}{s}}{\theta}(s,N)/(\eta ^{-1} (l^{s \pr}\frak{T})) $ is etale for all $s\pr \geq s$.  Because $N $ is of finite index, there is also a $t \geq 0$ so that $\eta ^{-1}(l^{t \pr} \frak{T}) \subseteq N $ for all $t\pr \geq t$.  Letting $S = max(s,t)$, we have that for every $s\pr \geq S$, the $G/N$-action on $\overline{\theta}(s \pr, N)$ is defined as the natural action of $G/N$ on the $N /\eta ^{-1}( l^{s \pr}\frak{T})$-orbits of the action of the $ G/\eta ^{-1}( l^{s \pr}\frak{T})$-action on  $\column{lim}{\column{\leftarrow}{s}}{\theta}(s,N)/(\eta ^{-1} (l^{s \pr}\frak{T})) $.  Since the $ G/\eta ^{-1}( l^{s \pr}\frak{T})$-action is etale, so is the $G/N$-action.  
\end{Proof}  

\begin{Proposition} ${\cal E}G$ is etale.
\end{Proposition} 
\begin{Proof} We must show that for any closed subgroup of finite index $N$, the $G/N$ action on $({\cal E}G)/N$ is etale. We have that $({\cal E}G)/N$ is the inverse limit of the functor $(s,N) \rightarrow \theta(s,N)/N$.  This is  a diagram, parametrized by a directed set, with values in the category of affine schemes with $G/N$-action.  For every closed normal subgroup of finite index $N \pr \subseteq G$, let $s(N\pr) \geq 0$ be an integer so that the $G/N$-action  on $\theta(s,N\pr)/N$ is etale for all $s \geq s(N\pr)$.  The subset  $$\hat{\frak{N} }= \{(s,N\pr) \mbox{ such that } s \geq s(N \pr) \}$$ contained in $\Bbb{N} \times \frak{N}^G$ is clearly final, and so ${\cal E}G/N$ can be described as the inverse limit of the restriction of the functor $N\pr \rightarrow \theta(s,N\pr ) /N$ to the subset $ \hat{\frak{N}}$.  Proposition \ref{inverseetale} now gives the desired result.  \end{Proof}

\begin{Proposition}  Let $W$ denote any continuous affine $G$-scheme, and assume that $A[W]$ is Noetherian approximable.  Then ${\cal E}G \times W$ is coherent, and the $G$-action is etale.  The natural map 
$$  K({\cal E}G \column{\times}{G} W) \rightarrow K^G({\cal E}G \times W) 
$$
therefore exists and is an equivalence.
\end{Proposition} 
\begin{Proof} ${\cal E}G \times W$ clearly is an etale action, using the base change property for etale actions.   The affine coordinate ring of ${\cal E}G \times W$ can be written as the colimit of rings of the form $A[W][x_1^{\pm \frac{1}{l^s}}, \ldots, x_t^{\pm \frac{1}{l^s}}]$, which are coherent by Proposition \ref{noetherianapproximable} since $A[W]$ is Noetherian approximable.  Further, each of the inclusions in the colimit system is of the form 
$$A[W][x_1^{\pm \frac{1}{l^s}}, \ldots, x_t^{\pm \frac{1}{l^s}}] \hookrightarrow A[W][x_1^{\pm \frac{1}{l^{s\pr}}}, \ldots, x_{t \pr}^{\pm \frac{1}{l^{s \pr}}}] 
$$
and are therefore clearly faithfully flat. This demonstrates the coherence of $A[{\cal E}G \times W]$.   The result now follows by Proposition \ref{orbitequiv}.   \end{Proof} 

{\bf Remark:} The conditions on $\theta^{\prime}$ provide a freeness property of ${\cal E}G$.  Extending by creating the product with $\theta^{\prime \prime}$ enlarges the pro-scheme so as to trivialize the action of $R[G]$ on the equivariant $K$-groups  of ${\cal E}G$.

We prove  the above mentioned  technical result  on the action of $R[G]$ on $K^G_*({\cal E}G)$.  which will be useful in analyzing completions of $K^G ({\cal E}G)$. 

\begin{Proposition} \label{trivialaction} Let $G$ be a totally torsion free $l$-profinite group.  Then the action of $\pi _0 K^G(k) \cong R[G]$ satisfies $I_G \cdot \pi _* K^G({\cal E}G) = \{ 0 \}$.  
\end{Proposition} 
\begin{Proof}  It clearly suffices  show that $j(I_G) = \{ 0 \}$.  This means that we will need to show that for any representation $\rho $ of $G$, we have that $[\rho ] - dim(\rho)$ vanishes in $K^G_0({\cal E}G)$.  Blichfeldt's theorem (see \cite{serre}) tells us that  any irreducible representation of $G$ is monomial, i.e it is isomorphic to an induced representation $i_{G_0}^G(\chi)$ where $G_0$ is a subgroup of index $dim(\rho)$ and $\chi$ is a one dimensional character of $G_0$.  Because of the totally torsion free property of $G$, $\chi$ is of the form $ \rho _s \compcirc \eta_{\chi}$ where $\eta _{\chi}$ is a one-dimensional $l$-BAR of $G_0$, and the $\rho _s$'s are defined in Section \ref{schemeactions}.   Inducing $\eta _{\chi}$ up to $G$ we obtain an $l$-BAR  $\eta$ of $G$, and it is clear from the definition that $\rho$ can be identified with $\rho _s \compcirc \eta$.  The conclusion is that every irreducible representation of $G$ is of the form $ \rho _s \compcirc \eta$ for some $l$-BAR of $G$, and some $s$.  Next, we observe that any $l$-BAR is a factor of a principal one.  This is an immediate consequence of Lemma \ref{principalfinal}.  For any $l$-BAR  $\eta$ we have the associated permutation representation $\overline{\eta}$.  We will  first prove that $[\rho _s \compcirc \eta] - [ \rho _s \compcirc \overline{\eta}] = 0$ in $K_0^G({\cal E}G)$. Fix $s \geq 0$.  Given an $l$-BAR $\eta$,  we let $K_{\eta}$ denote the kernel of $\rho _s \compcirc \eta$.  Similarly, let $\eta ^{\prime} $ be a $K_{\eta^{\prime}}$-principal  $l$-BAR containing $\eta$ as a factor, where $K_{\eta ^{\prime}} \subseteq K_{\eta}$.  We now have a factorization 
$$R[G/K_{\eta}] \cong K_0^{G/K_{\eta}}(k) \rightarrow K_0^{G/K_{\eta}} ({\cal T}^s_{\eta}) \stackrel{j}{\longrightarrow} K_0^{G/K_{\eta^{\prime}}}({\cal T}^s_{\eta \pr})\rightarrow K_0^G({\cal E}G)
$$
since ${\cal T}_{\eta}^s$ is one of the factors occurring in the inverse system defining  ${\cal T}_{\eta}$, and therefore one of the factors occurring in ${\cal E}G$. The map $j$ is the composite
$$ K_0^{G/K_{\eta}}{\cal T}_{\eta}^s \rightarrow K_0^{G/K_{\eta^{\prime} }} {\cal T}_{\eta}^s\rightarrow K_0^{G/K_{\eta ^{\prime}}}{\cal T}_{\eta \pr}^s
$$
where the left hand map is the inflation map associated to the projection $G/K_{\eta^{\prime}} \rightarrow G/K_{\eta}$, and the right hand map is the map induced by the equivariant  projection ${\cal T}_{\eta \pr}^s \rightarrow {\cal T}_{\eta}^s$.
 Letting $\rho = \rho_s \compcirc \eta$ and $\overline{\rho} = \rho _s \compcirc \overline{\eta}$ it will suffice to prove that the image of $[\rho] - [ \overline{\rho}] = 0$ in $K_0^G\Bbb{T}(k, \rho) $ vanishes.  Let $N$ be the kernel of  $\rho $, a normal subgroup of finite index in $G$. Note that $\rho$ and $\overline{\rho} $ may be regarded as elements in $K^{G/N}_0 \Bbb{A}(k, \rho)$.  We observe that there is a factorization 
$$ R[G/N] \cong K^{G/N}(k) \rightarrow K^{G/N}_0 \Bbb{A}(k, \rho ) \rightarrow K_0^G\Bbb{T}(k, \rho )\mbox{  } (= K_0^G{\cal T}_{\eta}^s)
$$We now select a basis $\{ x_1, \ldots , x_n \}$ for the representation space of $\rho$, where $n = dim(\rho)$, so that $\rho$ is monomial with respect to this basis, and we let $\overline{\rho} $ denote the associated permutation representation.  
 The group $K^{G/N}_0 \Bbb{A}(k, \rho)$ is the group completion of the commutative monoid of isomorphism classes of finitely generated left modules over the skew polynomial  ring 
 $A_{\rho} = k[x_1, \ldots , x_n]_{\rho} [G/N]$, by Proposition \ref{veryformal}.   We  construct a short exact sequence of left $A_{\rho}$-modules as follows.  First, let $G_0 \subseteq G$ be the stabilizer of the line spanned by the element $x_1$, so $\rho$  is induced from a one-dimensional representation of $G_0/N$, so that $\rho$ is induced from a one-dimensional representation $\chi$ of $G_0$.    We let $F$ be the module $k[x_1]_{\chi}[G_0/N] \column{\otimes}{k[G_0/N]}\epsilon$, where $\epsilon $ denotes the trivial representation of $G_0/N$. The module $F$ is free of rank one over $k[x_1]$, and we may write $F = k[x_1]b$, where $b = 1 \otimes 1$ is the  basis element which is acted on trivially by $G_0/N$.    We let $k$ denote the left $k[x_1]_{\chi}[G_0/N] $-module which as a $k$-module is cyclic of rank one, and on which all elements of $G_0/N$ act by the identity and such that $x_1 \cdot k = 0$.  There is a natural homomorphism $\pi: F \rightarrow k$, which is the identification $ F/xF \cong k$.  Let $F^{\prime}$ denote the kernel of $\pi$.  It is the free cyclic $k[x]$-module generated by $xb$, and from this description  it is clear that $F^{\prime}$ is isomorphic to $k[x_1]_{\chi}[G_0/N] \column{\otimes}{k[G_0/N]}\chi$.  The conclusion is that the equation $[k] = [\epsilon] - [\chi]$ holds in $K_0^{G_0/N}\Bbb{A}(k, \chi) \cong R[G_0/N]$, where we are interpreting this group as the projective class group of finitely generated modules over $k[x_1]_{\chi}[G_0/N]$.  Performing induction from $G_0/N$, and using Proposition \ref{inductionsquare}, we find that the equation 
 $$ [M] = [i_{G_0/N}^{G/N}(\epsilon) ] -  [i_{G_0/N}^{G/N}(\chi )] 
 $$ holds in  $K_0^{G/N}\Bbb{A}(k, \rho) \cong R[G/N]$, where 
 $$  M = k[x_1, \ldots , x_n] _{\rho}
  [G/N] \column{\otimes}{K[x_1] _{\chi} [G_0/N]} k
 $$
 We observe that $i_{G_0/N}^{G/N}(\epsilon)$ is the permutation representation $\overline{\rho} $ of $G$.  It is also readily checked that the element $x_1x_2, \cdots x_n$ annihilates the module $M$, and therefore that it becomes trivial on tensoring up to $k[x_1^{\pm1}, \ldots , x_n ^{\pm1}]$, which means that the equation 
 $$0= [i_{G_0/N}^{G/N}(\epsilon) ] -  [i_{G_0/N}^{G/N}(\rho] $$
 holds in 
 $K_0^{G/N}\Bbb{T}(k, \rho ) \cong R[G/N]$.  This gives the fact that $[\rho] -[\overline{\rho}] = 0 $ holds in $K_0({\cal E}G)$ for every irreducible representation  $\rho$ of $G$.   Since every permutation representation of a finite group is decomposable (it contains a copy of the trivial representation of $G$), this means that for every irreducible representation $\rho$ of $G$, the class $[\rho] \in K_0({\cal E}G)$ is equal to a sum of images of representations of smaller dimension.  Note that the special case $n=1$ of this result shows that one-dimensional representations $\chi$ satisfy the condition that  $ [\chi] - dim(\chi)$ vanishes in $K^G ({\cal E}G)$.  An induction on $n$ now shows that $[\rho]$ is equal to $dim(\rho)$, which was to be shown.   \end{Proof} 
 
 \begin{Corollary} \label{trivialactiontwo}  Let $W$ be any affine scheme with continuous $G$-action.  Then the action of $R[G]$ satisfies $I_G \cdot K^G({\cal E}G \times W)  = \{ 0 \}$. 
 \end{Corollary}
 \begin{Proof} This follows since $ K^G({\cal E}G \times W)  = \{ 0 \}$ is a ring spectrum, and we have that $I_G$ maps trivially into $\pi _0 K^G({\cal E}G \times W)  $ since there is a factorization 
 $$ R[G] \rightarrow K^G({\cal E}G)  \rightarrow K^G({\cal E}G \times W) 
 $$
 and we already know from Proposition \ref{trivialaction} that the map $R[G]\rightarrow K^G({\cal E}G)$ is trivial.  
 \end{Proof}
 
 Now we consider the diagram of commutative  ring spectra
 $$ \begin{diagram}
 \node{S^0} \arrow{e} \arrow{s,t}{r}   \node{K^G(k)} \arrow{s,b}{\epsilon}  \\
 \node{\Bbb{H}_l} \arrow{e,t} {=}\node{\Bbb{H}_l}
 \end{diagram} 
 $$
 where $r$ denotes the ring map induced by the fundamental class in the sphere spectrum, and $\epsilon$ is the augmentation.  
 The functorality of the derived completion construction means that the diagram induces a map 
 $$\beta : K^G({\cal E}G\times W) ^{\wedge}_l = K^G({\cal E}G \times W) ^{\wedge}_r  \rightarrow K^G({\cal E}G \times W) ^{\wedge}_{\epsilon} $$
\begin{Corollary} \label{nonequiv} The map $\beta $ is an equivalence of spectra.  
\end{Corollary} 
\begin{Proof}
 The algebraic-to-geometric spectral sequence (Theorem 7.1 of \cite{completion}) and the preceding Proposition show that it suffices to prove that the map on derived completions $\hat{\beta}: M^{\wedge}_l \rightarrow M^{\wedge}_{\epsilon}$ induced by the diagram 

 $$ \begin{diagram}
 \node{\Bbb{Z} } \arrow{e} \arrow{s,t}{r}   \node{R[G]} \arrow{s,b}{\epsilon}  \\
 \node{\Bbb{F}_l} \arrow{e,t} {=}\node{\Bbb{F}_l}
 \end{diagram} 
 $$
 is an equivalence for modules $M$ with trivial action by $I_G \subseteq R[G]$, or equivalently which are obtained by restricting scalars along the homomorphism $\epsilon : R[G] \rightarrow \Bbb{Z}$.  In order to prove this, we regard the ring $R[G]$ as a commutative  ring spectrum, and let $\tilde{\Bbb{Z}}$ denote a cofibrant replacement for the 
 $R[G]$-algebra $\Bbb{Z}$, with algebra structure given by $\epsilon$, in the model category of commutative  ring spectra discussed in \cite{shipley}.  We have the natural structure map $R[G] \rightarrow \tilde{\Bbb{Z}}$, as well as a homomorphism $\tilde{\Bbb{Z}} \rightarrow R[G]$ of commutative  ring spectra obtained by extending the homomorphism of rings $\Bbb{Z} \rightarrow R[G]$ to the cofibrant replacement $\tilde{\Bbb{Z}}$.  We also have a  homomorphism $\tilde{\Bbb{Z}} \rightarrow \Bbb{H}_l$ which is compatible with the mod-$l$ reduction of the augmentation on $R[G]$.    We now have two monads $S$ and $T$ on the category of $\tilde{\Bbb{Z}}$-module spectra, where $S(M) = \Bbb{H}_l \column{\wedge}{\tilde{\Bbb{Z}}} M$ and $T(M) = \Bbb{H}_l \column{\wedge}{R[G]} M$.  Further, there are natural transformations $T \rightarrow S$ and $S \rightarrow T$ induced by the ring homomorphisms $ \tilde{\Bbb{Z}} \rightarrow R[G] \rightarrow \tilde{\Bbb{Z}}$.  It follows that the natural transformations of triples induce equivalence on the total spectra of the cosimplicial  objects of the triples $S$ and $T$, by Theorem 2.15 of \cite{completion}.  These total spectra attached to $S$ and $T$ are the derived completions over $\tilde{\Bbb{Z}}$  and $R[G]$ respectively.   The fact that weak equivalences of commutative ring spectra induce equivalences of derived completions now gives the result.  
 \end{Proof}
 
 We also need a similar result for $K^G(\overline{F})$, where $G$ is the absolute Galois group of a field $F$, containing $k$ as a subfield.  
 
 \begin{Proposition}
 With $G,F,\overline{F}$ as above, the action of $R[G]$ on $K^G_*(\overline{F})$ is trivial in the sense that $I_G \cdot K^G_*(\overline{F}) = \{ 0 \}$.  
 \end{Proposition} 
 \begin{Proof}  It is proved in \cite{grothendieck} that the category of finitely generated $G$-twisted $\overline{F}$-modules is equivalent to the category of finite dimensional vector spaces over $F$, and that the equivalence can be taken by applying the functor $V \rightarrow \overline{F} \column{\otimes}{F}V$.   There is a natural twisted $G$-action on $\overline{F} \column{\otimes}{F}V$.  Since $\pi _0 K(F)$ is a free group of rank one, given only by the dimension, it is therefore clear that for any finite dimensional  $k$-linear representation $(W,\rho)$ of $G$, $\overline{F} \column{\otimes}{F} W$ is isomorphic to $\overline{F} \column{\otimes}{F} W^0$, where $W^0$ denotes $W$ equipped with trivial action.  This gives the result.  
 \end{Proof}
 \begin{Corollary}  \label{fineequiv} There is a natural equivalence $K^G(\overline{F})^{\wedge}_{\epsilon} \cong K(F)^{\wedge}_l$.
 \end{Corollary}
 \begin{Proof}  The proof is identical to that given in Proposition \ref{nonequiv} above. That proof uses only the triviality of the $I_G$-action on $\pi _0$, and so applies here as well.  
 \end{Proof}
 
 We conclude this section by proving that the equivariant $K$-theory of ${\cal E}G \times W$ can be computed as a particular kind of  colimit  of equivariant $K$-theory spectra of smooth schemes of finite type over $k$.  This will be important to us because the rigidity results we prove in \cite{royg} will apply only to such systems.  
 
 \begin{Definition}
 Let $G$ be a profinite group, and as usual let $\frak{N}^G$ denote the directed set of closed normal subgroups of finite index in $G$.  Let $F: \frak{N}^G \rightarrow \frak{C}^{cof}_{et}(G)$ denote any functor ($\frak{C}^{cof}_{et}(G)$ was defined in Definition \ref{cof} above).    We say $F$ is {\em adapted} if
 ${\cal N}(F(N)) = N$ for all  closed normal $N\subseteq G$.  Let $ \rho$ be any continuous homomorphism  from $G$ to the absolute Galois group of a field $K$. Let $L$ denote the absolute closure of $K$.  We define  spectrum valued functors $K(F,L)$ and  $K^G(F,L)$ associated to this data by 
 $$  K(F,L) = \column{hocolim}{N \in \frak{N}^G} K(F(N) \column{\times}{G/N} L^N)
 $$
 and 
 $$  K^G(F,L) = \column{hocolim}{N \in \frak{N}^G} K^{G/N}(F(N) \times  L^N)
 $$
 These colimits are naturally equivalent by Proposition \ref{free}.  $F$ is said to be {\em smooth} (respectively { \em of finite type}) if each $F(N)$ is smooth (respectively of finite type). 
 \end{Definition}  
 
 We are now able to express the $K$-theory of ${\cal E}G$ in these terms.   
 
  \begin{Proposition} \label{roycompare} Let $G$ be a totally torsion free profinite $l$-group, and let $\rho $ be a homomorphism into the absolute Galois group of a field $K$.  Let $L$ denote the algebraic  closure of $K$. Then there is a smooth  adapted functor $F(G,\rho): \frak{N}^G \rightarrow  \frak{C}^{cof}_{et}(G)$ of finite type, and an equivalence $ K^G(F(G, \rho), L) \rightarrow K^G({\cal E}G \times Spec(L))$.  This equivalence is natural for inclusions of fields with continuous $G$-action.  
  \end{Proposition}
  \begin{Proof}  For every closed normal subgroup of finite index $N \subseteq G$, Corollary \ref{etalequantitative} asserts that there is an integer $s(N)$ so that the action of $G/N$ on $\overline{\theta}(s \pr, N)$ for every $s \geq s(N)$, is etale.  We now define a function $j: \frak{N}^G \rightarrow \Bbb{N} \times \frak{N}^G$ by
  $$ j(N) = \column{max}{N \pr \subseteq N } s(N)
  $$
  It is clear from the definition that $J: \frak{N}^G \rightarrow \Bbb{N} \times \frak{N}^G$, defined by $J(N) = (j(N), N)$ is a map of directed sets, that $\overline{\theta} \compcirc J$ is an adapted functor,  and further that the image of $J$ is a final subset of $\frak{N}^G$.  It follows that we have the functors $K^G(\overline{\theta} \compcirc J, L)$ for any $\rho$ and $L$ as in the statement of the Proposition.   Moreover, we have an equivalence of colimits from $K^G(\overline{\theta} \compcirc J, L)$ to $K^G({\cal E}G \times Spec(L))$.      
  Finally, each of the schemes $F(N)$ are of finite type and smooth,  since they are  quotients by a finite etale group action on a torus.      \end{Proof} 
  
  {\bf Remark:} Note that the cases we are interested in are $K = L = k$ and $K =F $, $ L=\overline{F}$.
  
We write $F^G_L$ for $l$-adic completion of the the functor $N \rightarrow \theta (\nu (N))$ given  in the proof of the Proposition above.      We will prove that the maps $K^G(k) \rightarrow K^G({\cal E}G)$ and $K^G(Spec(\overline{F})) \rightarrow K^G ({\cal E}G \times Spec(\overline{F}))$ induce equivalences on derived completions along the homomorphism $\epsilon: K^G(k) \rightarrow \Bbb{H}_l$, and further that the completions of $K^G({\cal E}G)$ and $K^G ({\cal E}G \times Spec(\overline{F}))$ are just  $l$-adic completions.  Consider the diagram  
  $$
  \begin{diagram} \node{K^G(k)^{\wedge}_{\epsilon}} \arrow{e} \arrow{s} \node{K^G({\cal E}G)^{\wedge}_{\epsilon} \simeq  K^G({\cal E}G)^{\wedge}_l} \arrow{s}  \node{\column{hocolim}{N \in \frak{N}^G} F^G_k} \arrow{w} \arrow{s}  \\
   \node{K^G(\overline{F})^{\wedge}_{\epsilon}} \arrow{e} \node{K^G({\cal E}G\times Spec(\overline{F}))^{\wedge}_{\epsilon} \simeq  K^G({\cal E}G \times Spec(\overline{F}))^{\wedge}_l}  \node{\column{hocolim}{N \in \frak{N}^G}F^G_{\overline{F}}} \arrow{w} 
  \end{diagram} 
  $$
  
  We will be proving that the left hand horizontal arrows are equivalences.  The right hand horizontal arrows are equivalences by Proposition \ref{roycompare} above.   It follows that in order to prove that the left hand vertical arrow is an equivalence, it will suffice to prove that the right hand vertical arrow is one.  This is what will be proved in \cite{royg}.

\section{Algebraic properties of certain representation rings}

In this section, we analyze some particular properties of the representation rings of certain semidirect products of  symmetric groups with profinite groups of the form $\Bbb{Z}_l^n$.  The results of this section are the key technical tools which allow us to prove our main results.  

\begin{Theorem} \label{general} Let $k$ be a field of characteristic $l$, and let $\epsilon: A \rightarrow k$ be an augmented $k$-algebra with augmentation ideal $I_A$. Let $G \subseteq G^{\prime} \subseteq \Sigma _n$ denote subgroups of $l$-power order.   Let $ \displaystyle B = \mathop{\otimes} ^{n}_{k} A$, equipped with the action of $\Sigma _n$, and therefore of $G$ and $G^{\prime}$. The algebra $B$ is given the tensor product augmentation $\epsilon _B$, and   the algebras $B^G$ and $B^{G^{\prime}}$ both obtain augmentations by restriction of $\epsilon_B$, and we denote the restrictions by $\epsilon _{B^G}$ and $\epsilon _{B^{G^{\prime}}}$ respectively.  Further, we let $I_{B^G}$ and $I_{B^{G^{\prime}}}$ denote the corresponding augmentation ideals.  Finally, suppose that every element in $A$ has an $l$-th root.  Then the equality of ideals 
$$ I_{B^G} = I_{B^{G^{\prime}}}\cdot B^{G}
$$
holds. 
\end{Theorem}
\begin{Proof} We first prove that the $k$-algebra $B^{G}$ also has the property that every element admits an $l$-th root.  To see this, we let $\frak{A}$ denote a $k$-basis for $A$.  Then we form the product $\frak{B} = \frak{A}^n$, note that it forms a $k$-basis for $B$, and let $G$ act on $\frak{B}$ via the inclusion $G \hookrightarrow \Sigma _n$.  For every orbit $\frak{o} = \{ \beta _1, \ldots , \beta _s \}$,  we let $t(\frak{o})$ denote the element $\beta_1 +  \cdots +\beta _s \in B^{G}$.  It is clear that the elements $t(\frak{o})$, as $\frak{o}$ ranges over the orbits of the $G$ action on $\frak{B}$, form a $k$-basis for $B^G$, and therefore that it suffices to prove that each element $t(\frak{o})$ admits an $l$-th root in $B^G$.  We suppose that we are given an orbit $\frak{o}$, with orbit representative $\beta \in \frak{o}$.  Let $G_\beta$ denote the stabilizer of $\beta$.  The group $G_{\beta}$ now acts on the set $\{ 1, \ldots , n \}$ via its inclusion into $\Sigma _n$, and we let 
$$ \{ 1 , \ldots , n \} =  S_1 \cup S_2 \cup \cdots  \cup S_l
$$
denote the orbit decomposition of this action.  Since $G_{\beta}$ stabilizes $\beta$, it is clear that there are elements $\{ \alpha _j \}_{j = 1}^l$, with $\alpha \in \frak{A}$, so that 
$$ \beta = \alpha _{j(1) }\otimes \alpha _{j(2)} \otimes \cdots \otimes \alpha _{j(n)}
$$
where $j(i)$ denotes the integer so that $i \in S_{j(i)}$.  For each $j = 1, \ldots, l$, we select an element $\overline{\alpha}_j$ so that  $(\overline{\alpha}_j)^l = \alpha _j$.  Also, select left coset representatives $\{ g_1, \ldots , g_k \}$ for the collection of left cosets $G/G_{\beta}$. The element 
$$ \sum _{m=1}^k g_m( \overline{\alpha}_{j(1)} \otimes \overline{\alpha}_{j(2)} \otimes \cdots \otimes \overline{\alpha}_{j(l)})
$$
is now clearly an element in $B^G$, and it is equally clearly an $l$-th root of $t(\frak{o})$, which gives us the claim.  We now proceed to the proof of the theorem.

 We consider first the case where $G \lhd G^{\prime}$ and where the index of $G$ in $G^{\prime}$ is $l$.   Consider any element $b \in I_{B^G}$.  We let $\overline{b}$ denote any $l$-th root of $b$ in $B^G$.  It is immediate  that $\overline{b} \in I_{B^G}$ as well.  Now construct the  polynomial 
$$ \varphi (X) = \prod _{\gamma  \in G^{\prime}/G} (X - \overline{b}^{\gamma})
$$
Note that the notation $\overline{b}^{\gamma}$ is well defined for $\gamma \in G^{\prime}/G$ because $\overline{b} \in B^G$.  
From the definition of $\varphi$, it is clear that $\varphi (\overline{b}) = 0$, and that $\varphi$  is of the form 
$$\varphi (X) = X^l + \sum _{i=1}^l (-1)^i c_i X^{l-i}
$$ 
where $c_i$ is the $i$-th elementary symmetric function in the elements $\overline{b}^{\gamma}$, as $\gamma $ varies over all the elements of $G^{\prime}/G$.  Note that $c_i \in  I_{B^{G^{\prime} }}$. 
Since $\overline{b}^l = b$, it follows that 
$$ b =   \sum _{i=1}^l (-1)^{i+1} c_i \overline{b}^{l-i}
$$
which gives the result in this case.  For the general case, we may form a sequence of subgroups 
$$  G  = G_t \subseteq G_{t-1} \subseteq G_{t-2} \subseteq \cdots \subseteq G_{1} \subseteq G_0 = G^{\prime}
$$
where each inclusion is the inclusion of a normal subgroup of index $l$.  This is an immediate consequence of the fact that $G^{\prime}$ is an $l$-group.  The result for the general case can now be obtained from the case of a normal subgroup of index $l$ via an induction on $t$.  We wish to prove that $I_{B^{G_0}} B^{G_t} = I_{B^{G_t}}$, and by induction we suppose that we have proved that $I_{B^{G_{0}}}B^{G_{t-1}} =  I_{B^{G_{t-1}}}$ and $I_{B^{G_{t-1}}}B^{G_t} = I_{G_t}$.    We   obtain 
$$ I_{B^{G_t}} =  I_{B^{G_{t-1}}} B^{G_t}  = I_{B^{G_{0}}} B^{G_{t-1}} B^{G_t} =  I_{B^{G_{0}}} B^{G_t}
$$
 \end{Proof} 

We now  obtain the result which allows us to apply Corollary \ref{tootechnical}. 

\begin{Theorem}  \label{maindiagram} Let $k$ be a field of characteristic $l$, and $A$ an augmented  $k$-algebra in which every element has an $l$-th root.  Let $ G \subseteq G^{\prime} \subseteq \Sigma _n$ be an inclusion of $l$-groups, and let $ \displaystyle B = \mathop{\otimes} ^{n}_{k} A$ be equipped with the natural action of $\Sigma_n$ and therefore of $G$ and $G^{\prime}$. Also equip $B$ with the tensor product augmentation.  Suppose we have a diagram of augmented (to $k$) rings
$$
\begin{diagram}
\node{R^{\prime}} \arrow{s,t}{f} \arrow{e,t} {\pi ^{\prime}} \node{B^{G^{\prime}}} \arrow{s,b}{i} \\
\node{R} \arrow{e,t}{\pi} \node{B^{G}}
\end{diagram} 
$$
in which $\pi $ and $\pi ^{\prime}$ are surjective, and so that the kernel $\frak{I}_R$ of $\pi$ satisfies $(\frak{I}_R)^N = \{ 0 \}$ for some positive integer $N$.  Then there is an integer $M$ so that $(I_R)^M \subseteq  I_{R^{\prime}}\cdot  R$, where $I_R$ and $I_{R^{\prime}}$ denote the augmentation ideals of $R$ and $R^{\prime}$ respectively.

\end{Theorem} 
\begin{Proof}  Let $r \in I_R$.  Then $\pi (r) \in I_{B^G}$, and consequently, by Proposition \ref{general}, $\pi (r) \in I_{B^{G \pr}} B^G$.  We can therefore write 
$$ \pi (r) = \sum_s i(y_s)z_s
$$
where $y_s \in  I_{B^{G \pr}}$ and $z_s \in B^G$.  Since $\pi$ and  $\pi \pr$ are surjective, we can write $y_s = \pi \pr (\overline{y}_s)$ and $z_s = \pi (\overline{z}_s)$ for each $s$.  We let 
$$ \overline{r} = \sum _s f(\overline{y}_s) \overline{z}_s
$$
It is clear from the construction that $\pi(r - \overline{r}) = 0$, and of course that $\overline{r} \in 
I_{R \pr} R$.   We now let $M$ be any power of $l$ which is greater than $N$.  Then 
we have 
$$  0 = (r -\overline{r})^M = r^M - \overline{r}^M
$$
so $r^M = \overline{r}^M$, which gives that $r^M \in I_{R \pr} R$. 
\end{Proof}


We will use this result to study certain representation rings.  A first step is to study the behavior of the mod-$l$ representation rings of certain wreath products.  Recall that for any groups $G$ and $H$, with $G$ finite, the wreath product $G \wr H$ is the semidirect product $G \ltimes H^{\# (G)}$, where the $G$ action on $H^{\# (G)}$  is by permutation of coordinates.  For any group $G$, we will write $\frak{R}_l [G]$ for $R[G] \column{\otimes}{\Bbb{Z}} \Bbb{F}_l$.  

\begin{Proposition} \label{cycliccase}  Let $C_l$ denote the cyclic group of order $l$, where $l$ is a prime.  Then for any profinite group $K$, we construct $C_l \wr K$, with the inclusion $i: K^l \hookrightarrow C_l \wr K$.  The then following results hold.  
\begin{enumerate}
\item{  The restriction map $i^* : \frak{R}_l [C_l \wr K] \rightarrow \frak{R}_l [ K^l] $ is surjective onto the invariant subring $ \frak{R}_l [ K^l] ^{C_l}$.  }
\item{The kernel $\frak{I}$ of $i^*$ satisfies $\frak{I}^l = 0$. }
\end{enumerate} 
\end{Proposition} 
\begin{Proof}  Consider first the case where $K$ is finite.  Then the results of section 4.3 of \cite{james} show that an additive basis of $\frak{R}_l [C_l \wr K] $ can be constructed as follows.  Let $\mbox{Irr}(K)$ denote the set of irreducible representations of $K$.  Then we have the permutation action of $C_l$ on $\mbox{Irr}(K)^l$, which can then be decomposed into orbits, which either consist of one element or of $l$ elements.  Letting the orbit set be denoted by $\frak{O}$, we write $\frak{O} = \frak{O}_{(1)} \coprod \frak{O}_{(l)}$ for this decomposition.  A basis for $\frak{R}_l [C_l \wr K] $ is now in bijective correspondence with the set $\frak{O}_{(l)} \coprod \frak{O}_{(1)} \times \mbox{Irr}(C_l)$.  The correspondence is given by assigning to elements $\{ \rho _1, \ldots , \rho _l \}$  of $\frak{O}_{(l)}$ the induced representation $i_{K^l}^{C_l \wr K}(\rho _1)$, and on $ \frak{O}_{(l)} \times \mbox{Irr}(C_l)$ by assigning to elements $(\rho, e)$ the representation whose underlying $K^l$ representation is $\otimes ^l \rho$, and whose action by $C$ is given by permuting the tensor product, and to  $(\{\rho \} , \chi )$ the representation  $(\rho, e)\otimes \chi$.  The restriction map $i^*$ is now clearly surjective onto the invariant subring $\frak{R}_l [ K^l] $, which is the first result.  The kernel ideal is  clearly spanned by differences of the form $\langle \rho, \chi \rangle - \langle \rho, \chi _0 \rangle $, where $\chi _0$ denotes the trivial character of $C_l$.  Such elements lie in
$I_{\frak{R}_l[C_l]}\cdot \frak{R}_l[C_l \wr K]$, where  $\frak{R}_l[C_l \wr K]$ is an $\frak{R}_l[C_l]$-algebra via the homomorphism induced by the group homomorphism  $C_l \wr K \rightarrow C_l$.  Since the augmentation ideal in $\frak{R}_l [C_l ]$ has trivial $l$-th power, the second result follows.  The extension to profinite groups is a straightforward direct limit argument.  
\end{Proof} 

Consider an $l$-Sylow subgroup $\Bbb{S}(l,n)$ of $\Sigma _{l^n}$.  It is well-known that this group can be described as an iterated wreath product 
$$ \Bbb{S}(l,n) = \wr^{i} C_l
$$
where $C_l$ denotes a cyclic group of order $l$.  We also  consider the semidirect product 
$\frak{S}(l,n) = \Bbb{S}(l,n)\ltimes \Bbb{Z}_l^{l^n}$.  

\begin{Proposition} \label{sylow}The image of the homomorphism $i^* : \frak{R}_l[ \frak{S}(l,n)] \rightarrow \frak{R}_l[\Bbb{Z}_l^{l^n}]$ is  the invariant subring $\frak{R}_l[\Bbb{Z}_l^{l^n}]^{\Bbb{S}(l,n)}$.  If we let $\frak{I}$ denote the kernel of the homomorphism, then $\frak{I}^{l^n} = \{ 0 \}$.  
\end{Proposition} 
\begin{Proof} We first note that  $\frak{S}(l,n) \cong C_l \wr \frak{S}(l,n-1)$.  We now use this observation to perform an induction over $n$, as follows.  The $n = 1 $ case is Proposition \ref{cycliccase}.  For the case $n \geq 1$, the observation above show us the restriction map  $$
\frak{R}_l[\frak{S}(l,n)] \rightarrow \frak{R}_l[\frak{S}(l,n-1)^l] \cong \otimes^l \frak{R}_l[\frak{S}(l,n-1)]
$$
surjects on to the invariant subring $(\otimes^l \frak{R}_l[\frak{S}(l,n-1)])^{C_l}$.  We observe the general fact that if $f:V \rightarrow W$ is a homomorphism of $\Bbb{F}_l$-vector spaces, then the map 
$$ \otimes ^l f:( \otimes ^l V)^{C_l} \rightarrow (\otimes ^l W) ^{C_l}
$$
is also surjective, since the map $f$ admits a right inverse. By the inductive hypothesis,
the map 
$$ \frak{R}_l[\frak{S}(l,n-1)] \rightarrow \frak{R}_l[\Bbb{Z}^{l^{n-1}}_l]
$$ 
surjects onto the ring of invariants under the $\frak{S}(l,n-1)$ action on $\frak{R}_l[\Bbb{Z}^{l^{n-1}}_l]$.  What the induction now shows is that the image of the restriction map 
$$ \frak{R}_l[\frak{S}(l,n)] \rightarrow \frak{R}_l[\Bbb{Z}^{l^{n}}_l]
$$
surjects onto the $C_l$ invariant part of $(\otimes^l (\frak{R}_l[\Bbb{Z}^{l^{n-1}}_l])^{C_l}$.  Now it is easy to check that this invariant part is exactly identified with the $\frak{S}(l,n)$-invariant part of 
$\frak{R}_l[\Bbb{Z}_l^{l^n}]$ under the identification 
$$  \frak{R}_l[\Bbb{Z}_l^{l^n}] \cong \otimes ^l \frak{R}_l[\Bbb{Z}_l^{l^{n-1}}]
$$
which concludes the proof, which proves the first assertion.  The second is an easy induction using Proposition \ref{cycliccase}. 
\end{Proof}

We will need to work with some composite order finite groups as well.  We will need some terminology.

\begin{Definition}  For any finite group $G$, we will write $\hat{\frak{R}}_l[G]$ for the completion of the ring $\frak{R}_l[G]$ at its augmentation ideal.  For any profinite group $G$ and finite quotient $Q$ of $G$, we define 
$$ \hat{\frak{R}}_l[G;Q] = \hat{\frak{R}}_l[Q] \column{\otimes}{\frak{R}_l[Q]} \frak{R}_l[G]
$$
\end{Definition} 

We collect some elementary properties. 

\begin{Proposition} \label{manyproperties} Let $G$ be  a finite group. 
\begin{enumerate}
\item{If $G$ is an $l$-group, then the  natural homomorphism $\frak{R}_l[G] \rightarrow \hat{\frak{R}}_l[G]$ is an isomorphism. More generally, 
let $\Gamma$ denote a profinite group, with a quotient $G$ which is a finite $l$-group, then the natural homomorphism $\frak{R}_l[G] \rightarrow \hat{\frak{R}}_l[\Gamma; G]$ is an isomorphism.}
\item{For any subgroup $H \subseteq G$, the induction homomorphism $i_H^G:R[H] \rightarrow R[G]$ naturally extends to a homomorphism $\hat{i}_H^G: \hat{\frak{R}}_l[H] \rightarrow \hat{\frak{R}}_l[G]$.}
\item{We let $G_l$ denote the $l$-Sylow subgroup of $G$.  Then the restriction map $\hat{\frak{R}}_l[G] \rightarrow \hat{\frak{R}}_l[G_l]$  is an injection. More generally, let $\Gamma$ denote a profinite group, with a finite quotient $G$.  Let $G_l \subseteq G$ be a $l$-Sylow subgroup of $G$.  Then the natural map $\hat{\frak{R}}_l[\Gamma;G] \rightarrow \hat{\frak{R}}_l[\Gamma_l, G_l]$ is an injection, where $\Gamma _l$ denotes the inverse image of $G_l$ in $\Gamma$  under the projection $\Gamma \rightarrow G$.  }
\item{For any finite group $G$, the kernel of the (surjective) homomorphism 
$$ R[G] \rightarrow \hat{\frak{R}}_l[G]
$$
is the ideal 
$$  (l) +  ((l)+I_G)^{\infty}
$$
where for any ideal $J^{\infty}$ denotes the intersection $\bigcap _n J^n$.  More generally, given a profinite group $\Gamma$ and a finite quotient $G$, the kernel of the homomorphism 
$$ R[\Gamma] \rightarrow \hat{\frak{R}}_l[\Gamma;G]
$$
is the ideal $( (l) +  ((l)+I_G)^{\infty})R[\Gamma]$. }
\item{For any finite $l$-group $G$, the ideal $((l)+I_G)^{\infty} \subseteq R[G]$ is the zero ideal.  }
\end{enumerate} 
\end{Proposition}
\begin{Proof} Part (1) follows from  Proposition 3.5 of \cite{segal}. For part (2), we first quote Corollary 3.9 of \cite{segal}, which asserts that the $I_H$ topology and $I_G$ topology on $R[H]$ are the same.  It is easy to check that this extend directly to $\frak{R}_l$.  To prove that we have a natural extension of $i_H^G$, we need to check that for every $n$, there is a $d(n)$ so that $i_H^G(I_H^{n(d)}) \subseteq I_G^{n}$. The observation above now shows that it suffices to prove that for every $n$, there is some number $N(n)$ so that $i_H^G(I_G^{N(n)}\frak{R}_l[H]) \subseteq I_G^{n}\frak{R}_l[G]$.  But Frobenius induction has the property that $i_H^G(I_G^{n}\frak{R}_l[H]) \subseteq I_G^{n}\frak{R}_l[G]$, so we are done.  For part (3), we find that the composite 
$$ \hat{\frak{R}}_l[G] \rightarrow \hat{\frak{R}}_l[H]  \stackrel{\hat{i}_H^G}{\rightarrow} \hat{\frak{R}}_l[G] 
$$
has the property that $1$ is carried to an element $\lambda$  whose augmentation is equal to $\#(G/G_l)$.  This number is prime to $l$, and so is a unit in $\Bbb{F}_l$.  Since the augmentation ideal in $\hat{\frak{R}}_l[G]$ is nilpotent, it follows that $\lambda $ is a unit, so $\hat{\frak{R}}_l[G]$ is in fact a summand (as $\hat{\frak{R}}_l[G]$-modules) of $\hat{\frak{R}}_l[H]$. Part (4) follows immediately from the definitions.  Part (5) Is a direct consequence of Part (1).  
\end{Proof} 

\begin{Corollary} \label{symmcomplete} 
The ring homomorphism  $\hat{\frak{R}}_l[\Sigma _n \ltimes \Bbb{Z}_l^n; 
\Sigma _n ] \rightarrow \frak{R}_l[ \Bbb{Z}_l^n ]$ has image the invariant subring $ \frak{R}_l[ \Bbb{Z}_l^n ]^{\Sigma _n}$.  Further, its kernel $\frak{I}$ satisfies $\frak{I}^M = \{ 0 \}$, where $M$ is the smallest power of $l$ which is greater than $n$.  
\end{Corollary}
\begin{Proof} Let $\Sigma _n(l) \subseteq \Sigma _n$ denote an $l$-Sylow subgroup.  Letting $\{ \sigma _1, \ldots , \sigma _t \}$ denote a complete set of coset representatives for $\Sigma _n(l)$ in $\Sigma _n$, we have a homomorphism $\tau: \frak{R}_l[ \Bbb{Z}_l^n ]^{\Sigma _n(l)} \rightarrow \frak{R}_l[ \Bbb{Z}_l^n ]^{\Sigma _n}$ defined by 
$$  \tau (x) = \sum _i \sigma _i(x)
$$
The composite $$\frak{R}_l[ \Bbb{Z}_l^n ]^{\Sigma _n} \stackrel{i}{\rightarrow} \frak{R}_l[ \Bbb{Z}_l^n ]^{\Sigma _n(l) } \stackrel{\tau}{\rightarrow} \frak{R}_l[ \Bbb{Z}_l^n ]^{\Sigma _n}$$
is equal to multiplication by the index of $\Sigma _n(l)$ in $\Sigma _n$, from which it follows that $\tau$ is a surjection onto a summand.  It is now easy to check that there is a commutative diagram 
$$\begin{diagram}\node{\hat{\frak{R}}_l[\Sigma _n(l) \ltimes \Bbb{Z}_l^n; \Sigma _n(l)] } \arrow{s}\arrow{e} \node{ \frak{R}_l[ \Bbb{Z}_l^n ]^{\Sigma _n (l)}} \arrow{s,t} {\tau}\\
\node{\hat{\frak{R}}_l[\Sigma _n \ltimes \Bbb{Z}_l^n; \Sigma _n] } \arrow{e} \node { \frak{R}_l[ \Bbb{Z}_l^n ]^{\Sigma _n }}
\end{diagram}
$$
where the left hand vertical arrow is induction from $\Sigma_n(l) \ltimes \Bbb{Z}_l^n$ to $\Sigma _n \ltimes \Bbb{Z}_l^n$.  The upper horizontal arrow is surjective by Proposition \ref{sylow}, and we have observed that $\tau$ is surjective.  It follows that the lower horizontal arrow is also surjective, which is the first assertion in the Corollary.  The second statement follows directly from (3) in Proposition \ref{manyproperties}. 
%
\end{Proof} 

\begin{Corollary} \label{fintech} Let $n$ be a positive integer, and let $s$ and $t$ be positive integers so that $s+t= n$.    The the homomorphism $\hat{\frak{R}}_l[(\Sigma _s \times \Sigma_t )\ltimes \Bbb{Z}_l^n] \rightarrow \frak{R}_l [ \Bbb{Z}_l ^n]$ has image equal to the invariant subring $ \frak{R}_l [ \Bbb{Z}_l ^n]^{\Sigma _s \times \Sigma _t}$, and there is an integer $M$ so that the $M$-th power of the kernel ideal is $= \{ 0 \}$.  
\end{Corollary} 
\begin{Proof}  Immediate from Corollary \ref{symmcomplete}, since this example is simply a tensor product of two of the examples to which it applies.  
\end{Proof} 

\begin{Corollary} \label{reptech} Let $n, s$, and $t$ be as in the preceding corollary.  Then $\hat{\frak{R}}_l[(\Sigma _s \times \Sigma_t)\ltimes \Bbb{Z}_l^n;\Sigma_s \times \Sigma _t] $ becomes an  $\hat{\frak{R}}_l[\Sigma _n \ltimes \Bbb{Z}_l^n;\Sigma _n ] $-algebra under the ring homomorphism induced by the inclusion $(\Sigma _s \times \Sigma_t)\ltimes \Bbb{Z}_l^n \rightarrow \Sigma _n \ltimes \Bbb{Z}_l^n$.  We let $\hat{I}[-] \hookrightarrow \hat{\frak{R}}_l[-]$ denote the augmentation ideal.  Then there is an integer $M$  so that 

$$ \hat{I} [(\Sigma _s \times \Sigma_t )\ltimes \Bbb{Z}_l^n;\Sigma_s \times \Sigma _t]^M \subseteq \hat{I}[\Sigma _n \ltimes \Bbb{Z}_l^n ;\Sigma _n] \cdot \hat{\frak{R}}_l[(\Sigma _s \times \Sigma_t )\ltimes \Bbb{Z}_l^n;\Sigma _s \times \Sigma _t] 
$$
\end{Corollary} 
\begin{Proof}  This is a direct application of Theorem \ref{maindiagram}, with $R = \hat{\frak{R}}_l[(\Sigma _s \times \Sigma_t )\ltimes \Bbb{Z}_l^n;\Sigma _s \times \Sigma _t] $, $R^{\prime} = \hat{\frak{R}}_l[\Sigma _n \ltimes \Bbb{Z}_l^n ;\Sigma _n] $, $A = \frak{R}_l[\Bbb{Z}_l]$, $f$ is the homomorphism induced by the inclusion 
$(\Sigma _s \times \Sigma_t )\ltimes \Bbb{Z}_l^n \hookrightarrow \Sigma _n \ltimes \Bbb{Z}_l^n $, and where $\pi $ and $\pi ^{\prime} $ are the homomorphisms induced by the inclusions $\Bbb{Z}_l ^n \hookrightarrow (\Sigma _s \times \Sigma_t )\ltimes \Bbb{Z}_l^n$ and $\Bbb{Z}_l ^n \hookrightarrow \Sigma _n\ltimes \Bbb{Z}_l^n$ respectively.  That the hypothesis of Theorem \ref{maindiagram} is satisfied is guaranteed by  Corollary \ref{fintech}. 
\end{Proof} 
\begin{Corollary} \label{reptech2} Fix $n,s$, and $t$ as in the preceding corollary.  Then there is an integer $M$ so that 

$$I_{(\Sigma _s \times \Sigma _t) \ltimes \Bbb{Z}_l^n}^M \subseteq 
((l) + I_{\Sigma _n \ltimes \Bbb{Z}_l^n}  + ((l) + I_{\Sigma _s \times \Sigma _t})^{\infty}) R[\Sigma_s \times \Sigma _t \ltimes \Bbb{Z}_l^n]
$$

\end{Corollary}
\begin{Proof} Immediate consequence of Part (4) of Proposition   \ref{manyproperties}, together with Corollary \ref{reptech}.
\end{Proof}

\section{Equivariant $K$-theory of ${\cal A}_{\eta}$ and ${\cal T}_{\eta}$} \label{equivariantone}

We wish to study the map $K^{G} {\cal A}_{\eta} \rightarrow K^G {\cal T}_{\eta}$, with an eye to proving that it induces an equivalence on completions along the homomorphism of commutative ring spectra $\frak{K}^G(k) \rightarrow \Bbb{H}_l$.  More generally, we will want  to study the map   $K^{G}( {\cal A}_{\eta} \times W)\rightarrow K^G( {\cal T}_{\eta}\times W)$, and prove that it induces an equivalence on completions, where $W$ is a Noetherian approximable Affine scheme with continuous $G$-action.  We will carry out the analysis in the case $W = Spec(k)$, and then indicate any necessary modifications in the general case.   In order to do this, we will create a filtration by abelian subcategories on the category of all finitely generated modules over affine space whose support consists of coordinate hyperplanes (including the full affine space).  Of course, this filtration will induce a filtration of the $K$-theory spectra associated to the categories, via application of the localization sequence.  This section is entirely devoted to the analysis of the subquotients in this filtration as modules over the {\em weak ring spectra} $K(k)$ and $K^G(k)$.  We emphasize that all calculations in this section involve only modules over {\bf weak ring spectra}. These calculations will then be used in the next section to draw conclusions about the derived completions of modules over $\frak{K}^G(k)$, via comparisons such as the ones performed in \cite{repassembly}, section 4.

Let $\Bbb{A} = \Bbb{A}(k) = Spec(k[x_1, x_2, \ldots , x_n ])$ be the $n$-dimensional affine space over $k$, with explicit choice of coordinates $x_i$.  Of course, Theorem 7 of \cite{quillen} tells us that $K(\Bbb{A}) \cong K(k)$.   For any $1 \leq s \leq n$, we let $V(s)$ denote the hyperplane defined by $x_s = 0$.  More generally, let $S \subseteq \{ x_1, \ldots , x_n \}$, and define
$$V(S) = \bigcap _{s \in S} V(s)
$$
For any $i$, we will let $V_i$ denote the union 
$$  V_i = \bigcup _{\# (S) \geq i}V(S)
$$
Note that $V_{i+1} \subseteq V_i$.  We will also denote by $\Bbb{T}$ the scheme $\Bbb{A}- \cup _s V_s$, which is $Spec(k[x_1^{\pm 1}, \ldots , x_n ^{\pm 1} ]$.  
 We will refer to subvarieties of the form $V(S)$ as {\em coordinate subspaces}, and the corresponding ideals as {\em coordinate ideals}.  
 Note that the case $S = \emptyset$ is included, with the full $\Bbb{A}$ as corresponding variety.
 \begin{Definition}
Let  $\mbox{Mod}_i(\Bbb{A} )$ to be the subcategory of finitely generated modules whose support is contained in $V_i$.  For a set $S \subseteq \{ 1, \ldots , n \}$ of cardinality $i$, we also define $\mbox{Mod}_i(\Bbb{A} ;S) \subseteq \mbox{Mod}_i(\Bbb{A})$ to be the subcategory of modules supported on the subvariety $V(S) \cup V_{i+1}$.  We note that all the categories $\mbox{Mod}_i(\Bbb{A})$ and $ \mbox{Mod}_i(\Bbb{A};S)$ are abelian categories, and all inclusions are inclusions of Serre subcategories.  
Let $\Bbb{A}_S$ denote the subscheme  $\Bbb{A} - \bigcup _{s \notin S} V(s)$,  and $\overline{\Bbb{A}}_S \subseteq \Bbb{A}_S$ denote the subscheme $\{ (x_1, \ldots ,x_n) | x_i =0 \mbox{ for } i \in S \}$.  We also let $\mbox{\em Nil}(\Bbb{A};S)$ denote the category of coherent $\Bbb{A}_S$-modules on which the generators $\{ x_i \}_{i \in   S}$ act nilpotently. \end{Definition}


We clearly have an increasing sequence of subcategories
$$\{ 0 \} = \mbox{Mod} _{n+1}(\Bbb{A} ) \subseteq \mbox{Mod}_n(\Bbb{A} ) \subseteq \mbox{Mod}_{n-1}(\Bbb{A} ) \subseteq \cdots \subseteq \mbox{Mod}_1(\Bbb{A} ) \subseteq \mbox{Mod}_0(\Bbb{A} ) = \mbox{Mod}(\Bbb{A} )
$$

\begin{Proposition} \label{idtwo} There is a natural functor  from the  quotient abelian category  $\mbox{\em Mod}(\Bbb{A})/\mbox{\em Mod}_1(\Bbb{A})$  to the category $\mbox{\em Mod}(\Bbb{T})$, which induces an equivalence of $K$-theory spectra. 
\end{Proposition} 
\begin{Proof} This is simply the localization sequence associated to the removal of $V_1 $ from $V_0 \ \Bbb{A}$, which is $\Bbb{T}$. 
\end{Proof} 

We let $\Phi(\Bbb{A}, \Bbb{T})$ denote the homotopy fiber of the map of $K$-theory spectra 
$ K(\Bbb{A}) \rightarrow K(\Bbb{T})$.  It follows readily from Proposition \ref{idtwo} above that $\Phi(\Bbb{A}, \Bbb{T})$ is equivalent to the $K$-theory spectrum $K(\mbox{Mod}_1)$.  
\begin{Proposition} \label{filtration} There is a filtration of $\Phi = \Phi(\Bbb{A}, \Bbb{T})$ by spectra $\Phi _i =  \Phi _i (\Bbb{A}, \Bbb{T})$, with 
$$ * \simeq \Phi_{n+1} \subseteq \Phi _n \subseteq \Phi _{n-1} \subseteq \cdots \subseteq \Phi _2 \subseteq \Phi _1 = \Phi 
$$
and so that $\Phi _i/\Phi _{i+1} \cong K(\mbox{\em Mod}_i(\Bbb{A})/\mbox{\em Mod}_{i+1} (\Bbb{A}))$. \end{Proposition} 
\begin{Proof} The spectra $\Phi _i$ can be taken to be the homotopy fibers of the maps of spectra 
$$ K(\Bbb{A}) \rightarrow K(\mbox{Mod}(\Bbb{A})/\mbox{Mod}_i(\Bbb{A}))$$ 
\end{Proof} 

 Each inclusion $\mbox{Mod}_{i+1}(\Bbb{A} ) \hookrightarrow \mbox{Mod}_i(\Bbb{A} )$ is an  inclusion of an abelian subcategory, and it is clear from the discussion of quotient abelian categories in \cite{swan}  that the quotient $\mbox{Mod}_i (\Bbb{A} )/ \mbox{Mod}_{i+1}(\Bbb{A} ) $ can be analyzed as follows.  For each $S$ with $\# (S) = i$, we may consider the inclusion $i_S : \mbox{Mod}_i(\Bbb{A} ;S) \hookrightarrow  \mbox{Mod}_i(\Bbb{A} )$, and therefore the inclusion of subquotients $j_S: \mbox{Mod}_i(\Bbb{A} ;S) / \mbox{Mod}_{i+1}(\Bbb{A} ) \hookrightarrow  \mbox{Mod}_i (\Bbb{A} )/ \mbox{Mod}_{i+1}(\Bbb{A} ) $.  We may then sum over all $S$ of cardinality $i$  and apply $K$-theory to obtain 

\begin{equation} \label{special}\bigvee_{\#(S) = i} j_S: \bigvee_{\#(S) = i} K(\mbox{Mod}_i(\Bbb{A} ;S)  / \mbox{Mod}_{i+1}(\Bbb{A} ))  \longrightarrow  K(\mbox{Mod}_i(\Bbb{A} ) / \mbox{Mod}_{i+1}(\Bbb{A} ))
\end{equation}

This functor is an equivalence since as in the proof of Theorem 5.4 of \cite{quillen} it can be checked that the underlying functor is an equivalence of categories.  Further, there is a natural equivalence of abelian categories 
$$ \mbox{Mod}_i(\Bbb{A} ;S)  / \mbox{Mod}_{i+1}(\Bbb{A} ) \simeq \mbox{Nil}(\Bbb{A};S)
$$

 and hence  equivalences of $K$-theory spectra of these categories 
$$ K(\mbox{Mod}_i(\Bbb{A} ;S)  / \mbox{Mod}_{i+1}(\Bbb{A} )) \simeq K(\mbox{Nil}(\Bbb{A};S))
 \simeq K(\mbox{Mod}(\overline{\Bbb{A}}_{{S}}))
$$
The second of the equivalences  is proved as in \cite{quillen} by the devissage principle.  

We will also consider the situation of a product of a fixed Noetherian affine scheme $W$ with $\Bbb{A}$ and $\Bbb{T}$.  In this case, we define the subcategory $\mbox{Mod}_i(\Bbb{A},W) \subseteq \mbox{Mod} (\Bbb{A} \times W)$ to consist of all $\Bbb{A} \times W$-modules whose support is contained in  $\pi ^{-1}(V_i)$, 
where $\pi: \Bbb{A} \times W \rightarrow \Bbb{A}$ is the projection.  There are obvious analogues $\mbox{Mod}_i(\Bbb{A},W;S)$ and $\mbox{Nil} (\Bbb{A},W;S)$ of $\mbox{Mod}_i(A;S)$ and $\mbox{Nil}(\Bbb{A};S)$.







It is further clear from Section \ref{ktheory}  that all the spectra $$K(\mbox{Mod}_{i}(\Bbb{A} ,W)),K(\mbox{Mod}_i(\Bbb{A} ,W;S)), \mbox{ and }K(\mbox{Nil}(\Bbb{A},W;S))$$ are module spectra over the weak ring spectrum $K(\Bbb{A}\times W)$,  that all the inclusion and restriction functors induce maps of $K(\Bbb{A}\times W)$-module spectra, and therefore that the quotient spectra $K (\mbox{Mod}_{i}(\Bbb{A},W ))/K(\mbox{Mod}_{i+1}(\Bbb{A},W ))$ are also $K(\Bbb{A} \times W)$-module spectra.  It is also clear that $K(\mbox{Nil}(\Bbb{A},W;S))$ is naturally a $K(\Bbb{A}_S\times W)$-module spectrum, in a way which extends the $K(\Bbb{A}\times W)$-module structure via the map of  weak ring spectra $K(\Bbb{A}\times W) \rightarrow K(\Bbb{A}_S \times W)$.  Moreover, it is clear that the decomposition of spectra

$$ K(\mbox{Mod}_i(\Bbb{A} ,W)) / K(\mbox{Mod}_{i+1}(\Bbb{A},W )) \cong \bigvee_{\# (S) = i }  K(\mbox{Nil}(\Bbb{A},W;S))
\cong \bigvee_{\# (S) = i }  K(\mbox{Mod}(\overline{\Bbb{A}}_{{S}} \times W))
$$

given  by the homotopy inverse to the map $\bigvee_{\# (S) = i } K(j _s)$ from \ref{special} above  is also an equivalence of $K(\Bbb{A} \times W)$-module spectra, where the $K(\Bbb{A} \times W)$ action on the summand $K(\mbox{Nil}(\Bbb{A},W;S))$ and $K(\mbox{Mod}(\overline{\Bbb{A}}_{{S}} \times W))$ are   obtained by restriction of scalars along the map of  weak ring spectra $K(\Bbb{A}\times W) \rightarrow K(\Bbb{A}_S \times W)$.

We next begin the analysis of how this filtration works equivariantly.  Let $\rho$ be any monomial $k$-linear representation  of a finite group $G$, and let $\{ x_1, \ldots , x_n \}$ denote any basis for the representation space $V_{\rho}$ in which the representation is monomial.   Let $\overline{\rho}: G \rightarrow \Sigma _n$ denote the corresponding permutation representation.  As usual, we  write $\Bbb{A} = \Bbb{A}(k,\rho)$ for the associated affine variety, and therefore on the associated module category.    It is clear that the filtration of $\mbox{Mod}(\Bbb{A})$ by the subcategories $\mbox{Mod}_i(\Bbb{A})$ is invariant under the action defined by $\rho$, and we therefore easily obtain a filtration of $\mbox{Mod}^G(\Bbb{A})$ by Serre subcategories $\mbox{Mod}_i^G(\Bbb{A})$.  The equivariant version of Proposition \ref{filtration} holds. Let $\Phi^G = \Phi^G(\Bbb{A}, \Bbb{T})$ denote the homotopy fiber of the map of spectra 
$$  K^G(\Bbb{A}) \rightarrow K^G( \Bbb{T})
$$

\begin{Proposition}\label{filtrationone} There is a filtration of $\Phi^G $ by spectra $\Phi _i^G =  \Phi _i^G (\Bbb{A}, \Bbb{T})$, with 
$$ * \simeq \Phi_{n+1}^G \subseteq \Phi _n^G \subseteq \Phi _{n-1}^G \subseteq \cdots \subseteq \Phi _2^G \subseteq \Phi _1^G = \Phi ^G
$$
and so that $\Phi^G _i/\Phi ^G _{i+1} \cong K^G(\mbox{\em Mod}_i(\Bbb{A})/\mbox{\em Mod}_{i+1} (\Bbb{A}))$.
\end{Proposition}

We will have several categories  of $\Bbb{A}$-modules which are closed under the action of $G$ or of subgroups of $G$.  For each such category, we can then denote its equivariant version with the superscript $G$, and define  its equivariant $K$-theory as the $K$-theory of the equivariant version of the category.  For instance, we will write $K^G(\mbox{Mod}_i(\Bbb{A}))$ for $K(\mbox{Mod}^G_i(\Bbb{A}))$.   We now identify the subquotients $\mbox{Mod}_i^G(\Bbb{A})/\mbox{Mod}_{i+1}^G(\Bbb{A})$. Let ${\cal P}_{i}$ denote the set of subsets of $\{ 1, \ldots , n \}$ of cardinality $i$, and let 
$$ {\cal P}_i =  \frak{o}_1 \cup \cdots \cup \frak{o}_t
$$
be an orbit decomposition, with orbit representatives $S_j \in \frak{o}_j$. Note that for each $j$, the subcategory $\mbox{Mod}_i(\Bbb{A}; \frak{o}_j)$ is invariant under the action of $G$.  Then an argument entirely analogous to the non-equivariant one above shows that 

$$  \mbox{Mod}_i^G(\Bbb{A})/\mbox{Mod}_{i+1}^G(\Bbb{A}) \cong \bigvee _{j = 1}^t \mbox{Mod}^G_i(\Bbb{A};\frak{o}_j) / \mbox{Mod}^G_{i+1}(\Bbb{A})
$$
Moreover, it is now easy to see that each summand is described by 
$$  \mbox{Mod}^G_i(\Bbb{A};\frak{o}_j) / \mbox{Mod}^G_{i+1}(\Bbb{A}) \cong  \mbox{Mod}^{G_j}_i(\Bbb{A};S_j) / \mbox{Mod}^{G_j}_{i+1}(\Bbb{A})
$$
where $G_j$ is the stabilizer of the orbit representative $S_j$, and that further we have an equivalence of abelian categories
\begin{equation} \label{gamma} \mbox{Mod}^{G_j}_i (\Bbb{A}, S_j ) / \mbox{Mod}_{i+1}^{G_j}(\Bbb{A}) \cong \mbox{Nil}^{G_j}(\Bbb{A};{{S}_j})
\end{equation}
Also, devissage gives us an equivalence of spectra
\begin{equation} \label{delta} K(\mbox{Nil}^{G_j}(\Bbb{A};S_j))  \cong K^{G_j}(\mbox{Mod}(\overline{\Bbb{A}}_{{S}_j}))
\end{equation}

 and therefore a decomposition of $K$-theory spectra 

\begin{equation}\label{alpha}  K(\mbox{Mod}^G_i(\Bbb{A} ))/K(\mbox{Mod}^G_{i+1}(\Bbb{A}))  \cong 
\bigvee _{j=1}^t K^{G_j}(\mbox{Nil} (\Bbb{A};{{S}_l})) \cong  \bigvee _{j=1}^t K^{G_j}( \overline{\Bbb{A}} _{{S}_j})
\end{equation}

The corresponding decomposition for the case of a product with a fixed affine Noetherian scheme with $G$-action $W$ 

\begin{equation}\label{alphaone}  K(\mbox{Mod}^G_i(\Bbb{A},W ))/K(\mbox{Mod}^G_{i+1}(\Bbb{A},W))  \cong 
\bigvee _{j=1}^t K^{G_j}(\mbox{Nil} (\Bbb{A},W;S_l)) \cong  \bigvee _{j=1}^t K^{G_j}( \overline{\Bbb{A}} _{{S}_j}\times W)
\end{equation}
follows similarly

We now let $G$ be a profinite group, and let $\eta$ denote a fixed $l$-BAR of $G$.  We also let $W$ denote a continuous Noetherian affine $G$-scheme.  In practice, we will be concerned only with the case where $W = Spec(k)$, and the $G$-action is trivial,  or the case where $G= G_F$ is the absolute  Galois group of an extension $F$ of $k$, and $W = Spec(\overline{F})$.  We wish to analyze $K^G({\cal A}_{\eta} \times W)$, $K^G({\cal T}_{\eta}\times W)$, and the map $K^G({\cal A}_{\eta}\times W)\rightarrow K^G({\cal T}_{\eta}\times W)$.  We obtain $k$-linear representations $\rho _s \cdot \eta$ of $G$, and the corresponding varieties ${\cal A}_{\eta}^s$ and ${\cal T}_{\eta}^s$.  
We also denote the affine coordinate ring of ${\cal A}_{\eta}^s$ by $A_s$. We also  defined the $l$-th power map $\theta :  {\cal A}_{\eta}^{s+1} \rightarrow  {\cal A}_{\eta}^s$, which is equivariant.
 Because the map $\theta$ respects the coordinate axes, it is clear that it carries coordinate ideals into coordinate ideals.  From the definition of the subcategories $\mbox{Mod}^{G}_i$, it follows readily that the functor $A_{s+1} \column{\otimes}{A_s} - $ carries $\mbox{Mod}_i^{G}({\cal A}_{\eta}^s) \subseteq \mbox{Mod} ^{G} ({\cal A}_{\eta}^s)$ to the subcategory $\mbox{Mod}_i^{G}({\cal A}_{\eta}^{s+1}) \subseteq \mbox{Mod} ^{G}({\cal A}_{\eta}^{s+1})$.  Note that these tensoring functors are exact because ${A}_{s+1}$ is always a flat $A_s$-algebra.  This permits us to define categories $\mbox{Mod}_i^G({\cal A}_{\eta} )$ as the colimits of the systems 
$$\cdots \rightarrow \mbox{Mod}_{i} ^{G}({\cal A}_{\eta}^{s-1}) \rightarrow \mbox{Mod}_i^{G}({\cal {A}}_{\eta}^s)\rightarrow \mbox{Mod}_{i}^{G}({\cal A}_{\eta}^{s+1}) \rightarrow \cdots 
$$
A more intrinsic way to describe the subcategories of $\mbox{Mod}_i({\cal A_{\eta}})$ described above is as follows.  For any $s \in \{ 1, \ldots , n \}$, we let $V(s)$ denote the subscheme defined by the ideal $(x_s, x_s^{\frac{1}{l}} , x_s^{\frac{1}{l^2}}, \ldots )$.  We can now construct $V(S)$ for $S \subseteq \{ 1, \ldots , n \}$ as before, and the corresponding module categories of finitely presented modules supported on the corresponding subschemes as well.

We  note that  the analyses above apply equally well to a situation where we take products with a fixed affine Noetherian  $G$-scheme.  The arguments from above are adapted in a straightforward way.  

 We now have the analogue of the Proposition \ref{filtrationone}. Let $\Phi ^G({\cal A}_{\eta}, {\cal T}_{\eta})$ denote the homotopy fiber of the map $$
K^G({\cal A}_{\eta} ) \rightarrow K^G({\cal T}_{\eta} )
$$

\begin{Proposition} \label{newfiltration} There is a filtration of $\Phi^G({\cal A}_{\eta}, {\cal T}_{\eta}) $ by spectra $\Phi _i^G({\cal A}_{\eta}, {\cal T}_{\eta}) $, with 
$$ * \simeq \Phi_{n+1}^G({\cal A}_{\eta}, {\cal T}_{\eta}) \subseteq \Phi _n^G({\cal A}_{\eta}, {\cal T}_{\eta}) \subseteq \Phi _{n-1}^G({\cal A}_{\eta}, {\cal T}_{\eta}) \subseteq \cdots \subseteq \Phi _2^G ({\cal A}_{\eta}, {\cal T}_{\eta})\subseteq \Phi _1^G({\cal A}_{\eta}, {\cal T}_{\eta}) = \Phi ^G({\cal A}_{\eta}, {\cal T}_{\eta})
$$
and so that $$\Phi^G _i({\cal A}_{\eta}, {\cal T}_{\eta})/\Phi ^G _{i+1} ({\cal A}_{\eta}, {\cal T}_{\eta})\cong K^G(\mbox{\em Mod}_i({\cal A}_{\eta} )/\mbox{\em Mod}_{i+1} ({\cal A}_{\eta} ))$$
\end{Proposition} 

Further,    for any orbit $\frak{o}_j$  under the action of $G$ on  ${\cal P}_i$, $A_{s+1} \column{\otimes}{A_s} - $ carries the subcategory $\mbox{Mod}^G_i({\cal A}^s_{\eta}, \frak{o}_j)$ into the subcategory 
$\mbox{Mod}^{G}_i({\cal A}_{\eta}^{s+1}, \frak{o}_j)$.  It follows that the decompositions \ref{alphaone}  are respected by the map  $\theta$. We can therefore define subcategories 
$$ \mbox{Mod}_i^G({\cal A}_{\eta}; \frak{o}) \subseteq  \mbox{Mod}_i^G({\cal A}_{\eta})
$$ as the colimits of the system of categories
$$\cdots \rightarrow \mbox{Mod}_{i} ^{G}({\cal A}_{\eta}^{s-1};\frak{o}) \rightarrow \mbox{Mod}_i^{G}({\cal A}^s_{\eta};\frak{o} )\rightarrow \mbox{Mod}_{i}^{G}({\cal A}_{\eta}^{s+1};\frak{o}) \rightarrow \cdots 
$$
and we obtain a decomposition 
\begin{equation} \label{bardecomp} K^G(\mbox{Mod}_i({\cal A}_{\eta} )/\mbox{Mod}_{i+1}({\cal A}_{\eta} )) \cong \bigvee _{j = 1, \ldots t} 
K^G(\mbox{Mod}_i
({\cal A}_{\eta} ;
\frak{o_j})/\mbox{Mod}_{i+1}({\cal A}_{\eta} ))
\end{equation}
where $\{ \frak{o}_i, \ldots , \frak{o}_t \}$ is an orbit decomposition of the action of $G$ on the subsets of cardinality $i$ of $\{ 1, \ldots , n \}$.  


We will now focus on analyzing the effect of the map
$$  K(\mbox{Mod}^{G}_i({\cal A}_{\eta}^s;\frak{o}_j) / 
\mbox{Mod}^{G}_{i+1}({\cal A}_{\eta}^s)) \rightarrow K(\mbox{Mod}^{G}_i({\cal A}_{\eta}^{s+1};\frak{o}_j) / \mbox{Mod}^{G}_{i+1}({\cal A}_{\eta}^{s+1}))
$$
induced by $\theta$ on the individual factors, so as to understand the groups $$\pi _* K^G(\mbox{Mod}_i
({\cal A}_{\eta} ;
\frak{o_j})/\mbox{Mod}_{i+1}({\cal A}_{\eta} ))$$  We let $S_j$ denote a representative for the orbit $\frak{o}_j$, and let $G_j \subseteq G$ denote the stabilizer of $S_j$.  
We also have the identifications in \ref{gamma} and \ref{delta} above, and we wish to identify the effect of $\theta$ under this identification.  The category $\mbox{Nil}^{G_j}({\cal A}_{\eta}^s;{{S}_j})$ is described as follows.  We suppose that an orbit representative for the set $S_j$ is the set $\{1, \ldots , i \}$.  If not, perform a reordering of the indexing.  The permutation representation of $G_j$ preserves the decomposition 
\begin{equation} \label{sets} \{ 1, \ldots , i \} \coprod \{ i+1, \ldots , n \}
\end{equation} 
and yields decompositions $\rho _s \cdot \eta \simeq \sigma _s \oplus \tau _s$ and $\rho _{s+1}\cdot \eta \simeq \sigma _{s+1} \oplus \tau _{s+1}$, and corresponding decompositions of schemes with $G$-action \begin{equation}\label{splitone}{\cal A}_{\eta}^s = {\cal B}_{\eta}^s \times {\cal C}_{\eta}^s  = \Bbb{A}(k, \sigma _s) \times  \Bbb{A}(k, \tau _s)\end{equation}  for all $s$.  We also remove the coordinate hyperplanes corresponding to the elements of the complement of $S_j$, and obtain the variety
$$ \overline{{\cal A}}^s_{\eta} = {\cal{B}}^s_{\eta} \times \overline{{\cal C}}^s_{\eta} = \Bbb{A}(k, \sigma _s) \times \Bbb{T}(k, \tau _s)
$$
Moreover, the power map $\theta = \theta_{\rho}$ preserves the decomposition, so that we have a  commutative diagram 
$$ \begin{diagram} \node{\overline{{\cal A}}_{\eta}^{s+1}} \arrow{s,t}{\theta_{\rho}} \arrow{e,t}{\sim}  \node{{\cal B}_{\eta}^{s+1}\times \overline{{\cal C}}_{\eta}^{s+1}} \arrow{s,b}{\theta_{\sigma}\times \theta_{\tau}} \\
 \node{\overline{{\cal A}}_{\eta}^s} \arrow{e,t} {\sim} \node{{\cal B}_{\eta}^s \times \overline{{\cal C}}_{\eta}^s} 
\end{diagram} 
$$
of $G_j$-equivariant schemes.  The category 
$\mbox{Nil}^{G_j}({\cal A}_{\eta}^s;{{S}_j})$ corresponds to the category of all finitely generated $\overline{{\cal A}}_{\eta}^s$-modules  
  whose support is  concentrated on the subscheme $\{ 0 \} \times \overline{{\cal C}}_{\eta}^s$, and similarly for $\overline{{\cal A}}_{\eta}^{s+1}$. We let $A_s = k[x_1, \ldots , x_i, x_{i+1}^{\pm 1}, \ldots , x_n ^{\pm 1}], B_s = k[x_1, \ldots , x_i]$, and $C_s = k[x_{i+1}^{\pm 1} , \ldots , x_n^{\pm 1}] $ denote the corresponding affine coordinate rings, and similarly for $s+1$.  
\begin{Proposition} \label{products} Since the product 
decompositions \ref{splitone} are equivariant, we have the following. 
\begin{enumerate}
\item{ The projections $\overline{\cal A}_{\eta}^s \rightarrow \overline{\cal C}_{\eta}^s$ and $\overline{\cal A}_{\eta}^{s+1} \rightarrow \overline{\cal C}_{\eta}^{s+1}$ induce equivalences $K^{G_j}(\overline{\cal A}_{\eta}^s ) \cong K^{G_j}(\overline{\cal C}_{\eta}^s)$ and $K^{G_j}(\overline{\cal A}_{\eta}^{s+1} ) \cong K^{G_j}(\overline{\cal C}_{\eta}^{s+1})$ of weak ring spectra. }
\item{There are weak equivalences $$\varepsilon_s:K^{G_j}(\overline{\cal C}_{\eta}^s) \xrightarrow{\sim} K(\mbox{\em Nil}^{G_j}({\cal A}_{\eta}^s ;S_j))$$ and $$\varepsilon_{s+1}:K^{G_j}(\overline{\cal C}_{\eta}^{s+1})\xrightarrow{\sim} K(\mbox{\em Nil}^{G_j}({\cal A}_{\eta}^{s+1} ;S_j))$$ of modules over the weak ring spectra $ K^{G_j}({\cal C}_{\eta}^s)$ and $K^{G_j}({\cal C}_{\eta}^s)$, which is  given on $\pi _0$ by $K_0(\varepsilon _s )(1)  = [k \column{\otimes}{k} C_s]$, where $k$ is the quotient of $B_s$ by the maximal ideal of the origin, where the group action on $k$ is by the identity, and on $C_s$ is given by the already prescribed action.  The description of $\varepsilon _{s+1}$ is identical. }
\end{enumerate} 
\end{Proposition} 
\begin{Proof} The first result follows from the homotopy property in \cite{quillen}.  The second two are straightforward devissage arguments.  
\end{Proof} 

Next, we note that by the homotopy property of $K$-theory, we have equivalences $K^{G_j}(Spec(k))\rightarrow K^{G_j}({\cal B}_{\eta}^s) $ and $K^{G_j}(Spec(k))\rightarrow K^{G_j}({\cal B}_{\eta}^{s+1})$ induced by the projections ${\cal B}_{\eta}^s \rightarrow Spec(k)$ and ${\cal B}_{\eta}^{s+1} \rightarrow Spec(k)$, respectively.  These equivalences, taken together with the map $ K^{G_j}({\cal B}_{\eta}^s) \rightarrow K^{G_j}({\cal B}_{\eta}^{s+1})$ induced by $\theta _{\sigma}$, give a map $\overline{\theta} _{\sigma}: K^{G_j}(Spec(k) )\rightarrow K^{G_j}(Spec(k))$, which we wish to identify.  In order to do this, we consider a monomial $k$-linear representation $\nu$ of a finite group $G$, with the representation space $V_{\nu}$ of $G$ equipped with a basis $B_{\nu}$ in which the representation is monomial. Given any vector space $V$ equipped with a basis $B$, we may construct the symmetric algebra (over $k$) $S[V]$, and factor out the ideal generated by the set $\{ b^n \}_{b \in B}$.  We denote this vector space by $T_n(V,B)$.  Note that the construction only depends on the set of lines spanned by the individual basis elements, and can hence be defined as a functor on the category of vector spaces equipped with a family of independent and spanning lines.  This means that in our context, we can construct from the representation $\nu$ a new representation $T_n(\nu)$, which has $T_n(V_{\nu}, B_{\nu})$ as its representation space.  We can now describe the map $\overline{\theta}_{\sigma}$ as follows.  
\begin{Proposition} \label{characterize} The map $\overline{\theta}_{\sigma}: K^{G_j}(k) \rightarrow K^{G_j}(k)$ is characterized up to homotopy   by the following two requirements. 
\begin{enumerate}
\item{$\overline{\theta}_{\sigma}$ is a map of $K^{G_j}(k)$-modules.}
\item{$\pi _0(\overline{\theta}_{\sigma})(1) = T_l(\sigma _{s+1})$}
\end{enumerate} 
\end{Proposition}
\begin{Proof} The first statement is immediate  since the map $ K^{G_j}({\cal B}_{\eta}^s) \rightarrow K^{G_j}({\cal B}_{\eta}^{s+1})$ is a map of $K^{G_j}(k)$-modules.  For the second, we pass to affine coordinate rings.  We write $${\cal B}_{\eta}^{s+1} = Spec(B_{s+1}) = Spec(k[x_1, \ldots , x_i ]) $$ with the $G_j$-action given by an action on the  vector space spanned by $\{ x_1, \ldots , x_i \}$, which is monomial with respect to this basis.  We have $${\cal B}_{\eta}^s = Spec(B_s) = Spec(k[x_1^l, \ldots , x_i^l])$$ and the map $\theta _{\sigma}: {\cal B}_{\eta}^{s+1} \rightarrow {\cal B}_{\eta}^s$ is given by the evident inclusion $B_s \hookrightarrow B_{s+1}$.  The restriction of the $G_j$-action to $B_s$ factors through the homomorphism $G_j \rightarrow G_j$.  We are interested in the behavior of the functor $B_{s+1} \column{\otimes}{B_s}-$ restricted to the category $\mbox{Nil}^{G_j^s}({\cal B}_{\eta}^s)$, the category of modules whose support is the origin in the affine space ${\cal B}_{\eta}^s$.  A devissage argument shows that a generator for $K_0(\mbox{Nil}^{G_j}({\cal B}_{\eta}^s))$ is the field $k$, regarded as a module over $B_s$ via the trivial (zero) action of each of the generators $x_t^l$ and as a $G_j$-representation by declaring the action to be the identity action.  We now consider $B_{s+1} \column{\otimes}{B_s} k$.  It can clearly be identified with $B_{s+1}/(x_1^l, \ldots , x_i ^l)$, which as a $G_j$-representation is isomorphic to $T_l(\sigma _{s+1})$.  It now follows from the homotopy property that this element in $K_0^{G_j}({\cal B}_{\eta}^{s+1})$ is equal to the element $[T_l(\sigma_{s+1})]$.  
\end{Proof}

\begin{Corollary} \label{characterizationlemma} There is a commutative diagram 
$$\begin{diagram}
\node{K_i^{G_j}(k)} \arrow{e} \arrow{s,t}{T_l(\sigma _{s+1})\cdot} \node{K_i^{G_j}(\overline{\cal C}_{\eta}^s)} \arrow{e,t} {\varepsilon _{s}}\arrow{s,t}{ T_l(\sigma _{s+1}) \cdot K^{G_j}_i(\theta _{\tau}) } \node{ K_i(\mbox{\em Nil}^{G_j}({\cal A}_{\eta}^s ;S_j)) } \arrow{s,b} {K_i(\mbox{\em \tiny Nil}^{G_j}( \theta _{\rho};S_j))}\\
\node{K_i^{G_j}(k)} \arrow{e} \node{K_i^{G_j}(\overline{\cal C}_{\eta}^{s+1} )} \arrow{e,t} {\varepsilon _{s+1} }
\node{ K_i(\mbox{\em Nil}^{G_j}({\cal A}_{\eta}^{s+1} ;S_j))}
\end{diagram} 
$$ for all $i$,
where $\varepsilon _s$ and  $\varepsilon _{s+1}$ are isomorphisms.  
\end{Corollary} 
\begin{Proof} This is a simple product calculation, noting that ${\cal A}^s \cong {\cal B}_{\eta}^s \times {\cal C}_{\eta}^s$.  
\end{Proof}

\begin{Corollary} \label{computation} We let $i$ be an integer between $1$ and $n$, that $\frak{o}_j$ is an orbit of the action of $G$ on the collection of subsets of cardinality $i$ in $\{ 1, \ldots , n \}$ with orbit representative $S_j$, and that $G_j$ is the stabilizer of $S_j$. Then we have the following. 
\begin{enumerate}
\item{ The graded $R[G]$-module  $K_*^G(\mbox{ \em Mod}_i({\cal A}_{\eta}; \frak{o}_j)/\mbox{\em Mod }_i ({\cal A}_{\eta} ))$ is isomorphic to the colimit of the system 
\begin{equation} \label{colimitfinal}K_*^{G_j}(\overline{{\cal C}}^0_{\eta})  \stackrel{T_l(\sigma _1)\cdot K^{G_j}_*(\theta _{\tau}) }{\longrightarrow} K_*^{G_j}(\overline{{\cal C}}^1_{\eta})\stackrel{T_l(\sigma _2)\cdot K^{G_j}_*(\theta _{\tau}) }{\longrightarrow}K_*^{G_j}(\overline{{\cal C}}^2_{\eta}) \stackrel{T_l(\sigma _3)\cdot K^{G_j}_*(\theta _{\tau}) }{\longrightarrow}\cdots 
\end{equation}
where $\sigma _s$ is the $i$-dimensional representation of $G$ on the basis elements corresponding to the subset $S_j$ on the affine space $\Bbb{A}_s$. It follows that it is isomorphic to the colimit of the system 
\begin{equation} \label{colimitfinalone}K_*^{G_j}(\overline{{\cal C}}_{\eta})  \stackrel{T_l(\sigma _1) }{\longrightarrow} K_*^{G_j}(\overline{{\cal C}}_{\eta})\stackrel{T_l(\sigma _2) }{\longrightarrow}K_*^{G_j}(\overline{{\cal C}}_{\eta}) \stackrel{T_l(\sigma _3)  }{\longrightarrow}K_*^{G_j}(\overline{{\cal C}}_{\eta}) \stackrel{T_l(\sigma _4) }{\longrightarrow}\cdots\end{equation}
 }
\item{Let $\overline{G}_j$ denote the image of $G_j$ under $\eta$.  Then all of the elements $T_l(\sigma _s) \in R[G_j]$ are in the image of the ring homomorphism $R[\overline{G}_j] \rightarrow R[G_j]$.}
\end{enumerate}
\end{Corollary}
\begin{Proof} For part (1), we observe that the colimit (\ref{colimitfinal}) above can be written  as a two parameter colimit, of the following form. 
$$\begin{diagram}
\node{K_*^{G_j}(\overline{{\cal C}}^0_{\eta}) } \arrow{e,t}{T_l(\sigma_1)} \arrow{s,t}{K^{G_j}_*(\theta _{\tau})} \node{K_*^{G_j}(\overline{{\cal C}}^0_{\eta}) } \arrow{s,t}{K^{G_j}_*(\theta _{\tau})} \arrow{e,t}{T_l(\sigma_2)} \node{K_*^{G_j}(\overline{{\cal C}}^0_{\eta}) } \arrow{s,t}{K^{G_j}_*(\theta _{\tau})} \arrow{e,t}{T_l(\sigma_3)} \node{\cdots} \\
\node{K_*^{G_j}(\overline{{\cal C}}^1_{\eta}) } \arrow{e,t}{T_l(\sigma_1)} \arrow{s,t}{K^{G_j}_*(\theta _{\tau})}  \node{K_*^{G_j}(\overline{{\cal C}}^1_{\eta}) } \arrow{s,t}{K^{G_j}_*(\theta _{\tau})} \arrow{e,t}{T_l(\sigma_2)} \node{K_*^{G_j}(\overline{{\cal C}}^1_{\eta}) } \arrow{s,t}{K^{G_j}_*(\theta _{\tau})} \arrow{e,t}{T_l(\sigma_3)} \node{\cdots} \\
\node{K_*^{G_j}(\overline{{\cal C}}^2_{\eta}) } \arrow{e,t}{T_l(\sigma_1)} \arrow{s,t}{K^{G_j}_*(\theta _{\tau})}  \node{K_*^{G_j}(\overline{{\cal C}}^2_{\eta}) } \arrow{s,t}{K^{G_j}_*(\theta _{\tau})} \arrow{e,t}{T_l(\sigma_2)} \node{K_*^{G_j}(\overline{{\cal C}}^2_{\eta}) } \arrow{s,t}{K^{G_j}_*(\theta _{\tau})} \arrow{e,t}{T_l(\sigma_3)} \node{\cdots} \\
\node{\vdots}   \node{\vdots}   \node{\vdots}   \node{} \\
\end{diagram} 
$$
This colimit may be computed in iterative fashion, first in the vertical direction and then in the horizontal.  The diagram (\ref{colimitfinalone}) above is obtained from this two parameter diagram by performing the vertical colimits.  For part (2) above, it is clear that the representations $\sigma _s$ all factor through $\overline{G}_j$ by definition.  
\end{Proof}

\section{Completed equivariant $K$-theory of  ${\cal T}_{\eta}$}

We suppose that we are given an $l$-BAR $\eta$, and construct the spectra ${\cal A}_{\eta} $ and ${\cal T}_{\eta}$.  Also let $W$ be a Noetherian approximable continuous affine $G$-scheme, as in the previous section.   We wish to prove that the map 
$$\theta: \frak{K}^G (k) \cong \frak{K}^G({\cal A}_{\eta} \times W) \rightarrow  \frak{K}^G({\cal T}_{\eta}\times W)
$$ induces an equivalence on  derived completions at the homomorphism
$$ \epsilon: \frak{K}^G(k) \rightarrow \Bbb{H}_l
$$
where $\Bbb{H}_l$ denotes the mod $l$ Eilenberg spectrum.  
In order to do this, it will suffice to prove that the map 
$$  \Bbb{H}_l \cong \frak{K}^G({\cal A}_{\eta} \times W) \column{\wedge}{\frak{K}^G(k)} \Bbb{H}_l \xrightarrow{ \theta \wedge id} \frak{K}^G({\cal T}_{\eta}\times W) \column{\wedge}{\frak{K}^G(k)} \Bbb{H}_l
$$ is an equivalence of spectra by Proposition \ref{criterionone}.  
This analysis  will require the filtrations discussed in Section \ref{equivariantone}.  The computations there were exclusively using the weak ring spectra and modules constructed from Waldhausen's $S_.$-construction.  In order to apply those results, we have to construct the $\frak{K}^G$-counterparts to the spectra $K^G(\mbox{Mod}_i(-))$ and $K^G(\mbox{Mod}^G_i(-)/\mbox{Mod}^G_{i+1}(-))$. To begin, we suppose  that we are given a monomial representation $l$-BAR of $G$ , and the corresponding actions of $G$ on ${\cal A}_{\eta}^s$ and ${\cal T}_{\eta}^s$.  In Section \ref{equivariantone}, we defined the abelian subcategories $\mbox{Mod}^G_i({\cal A}_{\eta}\times W)$ of the category $\mbox{Mod}^G({\cal A}_{\eta}\times W)$. We now  define abelian categories $Q_i$ to be the quotient abelian categories $ \mbox{Mod}^G({\cal A}_{\eta} \times W)/\mbox{Mod}^G_i({\cal A}_{\eta} \times W) $.  Each of the categories $Q_i$ admits a strictly associative and coherently commutative tensor produce, and so is therefore suitable input to the Elmendorf-Mandell machine.  Consequently, we may construct commutative $\frak{K}^G(k)$-algebras  $\frak{K}^G(Q_i)$, and we have natural homomorphisms $\frak{K}^G(Q_i) \rightarrow \frak{K}^G(Q_{i-1})$ of such algebras.  We now define the module $\frak{M}^G_i({\cal A}_{\eta}\times W )$ over $\frak{K}^G(k)$ to be the homotopy fiber of the  map
$$ \frak{K}^G(k)({\cal A}_{\eta} \times W) \rightarrow \frak{K}^G(Q_i)
$$
The localization theorem in Waldhausen theory shows that $\frak{M}^G_i({\cal A}_{\eta}\times W )$ has the same homotopy type as $K(\mbox{Mod}_{i}^G({\cal A}_{\eta} \times W ))$ as a weak module spectrum over $K^G({\cal A}_{\eta}\times W )$.  We now have natural homomorphisms of $\frak{K}^G(k)$-modules $\frak{M}^G_{i+1}({\cal A}_{\eta}\times W ) \rightarrow \frak{M}^G_{i} ({\cal A}_{\eta} \times W)$, and we can construct their cofibers $\frak{M}^G_i({\cal A}_{\eta} \times W)/\frak{M}^G_{i+1}({\cal A}_{\eta}\times W )$ to obtain $\frak{K}^G$ counterparts to the subquotients $K^G(\mbox{Mod}_i^G({\cal A}_{\eta} \times W )/\mbox{Mod}_{i+1}^G({\cal A}_{\eta} \times W))$.  We also want to construct $\frak{K}^G$-counterparts $\frak{M}_i^G({\cal A}_{\eta} \times W ;\frak{o})$ to the spectra  $K^G(\mbox{Mod}^G_i({\cal A}_{\eta}\times W ;\frak{o}))$  where $\frak{o}$ is any orbit under the action of $G$ on the collection of subsets of cardinality $i$ of $\{ 1, \ldots , n \}$.  In order to do this, we proceed as above by letting $Q_i(\frak{o})$ denote the quotient abelian category $ \mbox{Mod}^G({\cal A}_{\eta}\times W) /\mbox{Mod}^G_i({\cal A}_{\eta}\times W ;\frak{o}) $, and defining $\frak{M}_i^G({\cal A}_{\eta} \times W ;\frak{o})$ to be the homotopy fiber of the map $\frak{K}^G(k) \rightarrow \frak{K}^G(Q_i(\frak{o}))$.  We may apply the $K^G$ case to conclude that there is a decomposition
\begin{equation}\label{fancydecomp} \frak{M}_i^G({\cal A}_{\eta} \times W)/\frak{M}_{i+1}^G({\cal A}_{\eta}\times W ) \cong \bigvee _{\frak{o}} \frak{M}_i^G({\cal A}_{\eta} \times W ;\frak{o}) / \frak{M}^G_{i+1}({\cal A}_{\eta}\times W ) 
\end{equation}

We now have the following Proposition concerning these modules. 

\begin{Proposition} \label{fibers} There are natural maps $\frak{K}^G({\cal A}_{\eta} \times W) \rightarrow Q_{\mbox{dim}(\eta) +1} \rightarrow Q_1 \rightarrow \frak{K}^G({\cal T}_{\eta}\times W )$. 
\begin{enumerate}
\item{The map $\frak{K}^G({\cal A}_{\eta}\times W ) \rightarrow Q_{\mbox{dim}(\eta)+1}$ is an equivalence of spectra, so we have $\frak{M}^G_{\dim (\eta) + 1} ({\cal A}_{\eta} \times W)\simeq	 *$.}
\item{The map  $Q_1 \rightarrow \frak{K}^G({\cal T}_{\eta} \times W )$ is an equivalence of spectra, so $\frak{M}_1^G({\cal A}_{\eta} \times W)$ is  equivalent to the homotopy fiber of the map $\frak{K}^G({\cal A}_{\eta} \times W) \rightarrow \frak{K}^G({\cal T}_{\eta} \times W )$.}
\end{enumerate} 
\end{Proposition} 
\begin{Proof} For item (1), we note that the category $\mbox{Mod}^G_{\mbox{dim}(\eta)+1}({\cal A}_{\eta}\times W )$ consists only of the zero module, from which (1) follows.  For (2), $Q_0$ consists of all modules on ${\cal A}_{\eta} \times W $ for which the support consists of coordinate subspaces, including the zero subspace.  The localization is equivalent to the category of free modules over ${\cal T}_{\eta} \times W$.  
\end{Proof} 


\begin{Proposition} In order to prove that the map 
$$  \Bbb{H}_l \cong \frak{K}^G({\cal A}_{\eta}\times W ) \column{\wedge}{\frak{K}^G(k)} \Bbb{H}_l \xrightarrow{ \theta \wedge id} \frak{K}^G({\cal T}_{\eta}\times W) \column{\wedge}{\frak{K}^G(k)} \Bbb{H}_l
$$ is  an equivalence of spectra, it will suffice to prove that 
$$ \frak{M}^G_i({\cal A}_{\eta} \times W ;\frak{o} )/  \frak{M}_{i+1}^G({\cal A}_{\eta} \times W) \column{\wedge}{\frak{K}^G(k)} \Bbb{H}_l \cong *
$$
for all $i =1, \ldots , n$ and all orbits $\frak{o}$ of the action of $G$ on the subsets of cardinality $i$ of $\{ 1, \ldots , n \}$.  
\end{Proposition} 
\begin{Proof} Since by Proposition \ref{fibers}, $\frak{M}^G_1({\cal A}_{\eta}  \times W )$ is equivalent to the homotopy fiber of the map 
$$\frak{K}^G({\cal A}_{\eta} \times W) \rightarrow \frak{K}^G({\cal T}_{\eta} \times W )$$ it will clearly suffice to prove that 
$$ \frak{M}^G_1({\cal A}_{\eta} \times W  ) \column{\wedge}{\frak{K}^G(k)} \Bbb{H}_l \simeq *
$$
But we have a diagram of module spectra
$$ * \simeq \frak{M}^G_{dim(\eta)+1} ({\cal A}_{\eta} \times W )  \rightarrow \frak{M}^G_{dim(\eta)}({\cal A}_{\eta} \times W ) \rightarrow \cdots \rightarrow   \frak{M}^G_2({\cal A}_{\eta} \times W ) \rightarrow  \frak{M}^G_1({\cal A}_{\eta} \times W  ) 
$$
from which it is clear that it suffices to prove that
$$ \frak{M}^G_i({\cal A}_{\eta}  \times W) / \frak{M}^G_{i+1} ({\cal A}_{\eta} \times W ) \column{\wedge}{\frak{K}^G(k)} \Bbb{H}_l \simeq *
$$
for $i =1,2,\ldots , dim(\eta)$.  The result now follows easily from the decomposition (\ref{fancydecomp}) above.
\end{Proof} 
\begin{Corollary} To prove that 
$$  \Bbb{H}_l \cong \frak{K}^G({\cal A}_{\eta} \times W) \column{\wedge}{\frak{K}^G(k)} \Bbb{H}_l \xrightarrow{ \theta \wedge id} \frak{K}^G({\cal T}_{\eta}\times W) \column{\wedge}{\frak{K}^G(k)} \Bbb{H}_l
$$ is  an equivalence, it will suffice to prove that the groups 
$$ Tor^{R[G]}_j(\pi _k (\frak{M}^G_i({\cal A}_{\eta} \times W   ; \frak{o})/  \frak{M}_{i+1}^G({\cal A}_{\eta} \times W )) , \Bbb{F}_l)
$$
vanish for all $i=1, \ldots , n$, and all $ j, k$, and $\frak{o}$.  

\end{Corollary}
\begin{Proof} Follows immediately from Proposition \ref{algcriterion}. 
\end{Proof} 

Let $\eta $ denote an $n$-dimensional  $l$-BAR of $G$, and consider an orbit $\frak{o}$ under the action of $G$ on the subsets of order $i$ of $\{ 1, \ldots , n \}$.  Let $H$ be the stabilizer of one orbit representative in $\frak{o}$.  We now have a commutative diagram 
\begin{equation}\label{gooddiagram} \begin{diagram}  \node{H} \arrow{s,t}{ \iota}\arrow{e,t}{j} \node{\Sigma_i \times \Sigma _{n-i} \ltimes \Bbb{Z}_l^n} \arrow{s} \\
\node{G} \arrow{e,t}{k} \node{ \Sigma _n \ltimes \Bbb{Z}_l^n} 
\end{diagram} 
\end{equation}
and corresponding diagram of representation rings
\begin{equation}\label{betterdiagram}
\begin{diagram}
\node{\hat{\frak{R}}_l[\Sigma _n \ltimes \Bbb{Z}_l^n;\Sigma _n]} \arrow{s} \arrow{e}\node{\hat{\frak{R}}_l[G;Q_G]} \arrow{s} \node{\frak{R}_l[G]} \arrow{s} \arrow{w,t}{\sim}\\
\node{\hat{\frak{R}}_l[ \Sigma _i \times \Sigma _{n-i} \ltimes \Bbb{Z}_l^n;\Sigma _i \times \Sigma _{n-i}]} \arrow{e} \node{\hat{\frak{R}}_l[H,Q_H]} \node{\frak{R}_l[H]} \arrow{w,t}{\sim}
\end{diagram}
\end{equation}
where the quotients $Q_G$ and $Q_H$ are the images of the composites 
$$ G \stackrel{k}{\longrightarrow} \Sigma _n \ltimes \Bbb{Z}_l^n \rightarrow \Sigma _n
$$
and 
$$ H \stackrel{j}{\longrightarrow} \Sigma _i \times \Sigma _{n-i} \ltimes \Bbb{Z}_l^n \rightarrow \Sigma _i \times \Sigma _{n-i}
$$
The right hand horizontal maps are isomorphisms since the groups $Q_G$ and $Q_H$ are finite $l$-groups, by Proposition \ref{manyproperties}.  
We draw the following consequence of  Proposition \ref{computation}.  

\begin{Proposition}  The $R[G]$ -module $\pi _k (\frak{M}^G_i({\cal A}_{\eta} \times W   ; \frak{o})/  \frak{M}_{i+1}^G({\cal A}_{\eta} \times W )) $ is of the form 
$$  M \stackrel{\times r_0}{\longrightarrow}  M \stackrel{\times r_1}{\longrightarrow} M \stackrel{\times r_2}{\longrightarrow} M \stackrel{\times r_3} {\longrightarrow} \cdots
$$
where
\begin{enumerate}
 \item{$M$ is a module over $R[H]$.}
\item{Each $r_j$ is a representation of dimension divisible by $l$}
 \item{Each $r_j$ is the pullback along $j$ of a representation $\rho _i$  of $\Sigma_i \times \Sigma _{n-i} \ltimes \Bbb{Z}_l^n$. }
 \end{enumerate}
\end{Proposition} 
\begin{Proof}  The first follows from Part (1) of Corollary \ref{computation}, with $M = K^{G_j}_i(\overline{\cal C}_{\eta})$.  The second  follows because the construction $T_l$ used in Proposition \ref{computation} always yields dimension a multiple of $l$.  The last statement follows since the representations $\sigma _s$ are pulled back along a homomorphism to $\Sigma _i \ltimes \Bbb{Z}_l^i$, which is a subgroup of  $\Sigma _i \times \Sigma _{n-i} \ltimes \Bbb{Z}_l^n$.  
\end{Proof} 
\begin{Corollary} There exists an $M$ so that for all $i$, the product $r_i \cdot r_{i+1} \cdot \cdots \cdot r_{i+M-1} \cdot r_{i+M}$ lies in the submodule $(l)+I_G\cdot R[H]$. 
\end{Corollary} 
\begin{Proof}  We note that each $\rho_i$ has dimension divisible by $l$, and it readily follows that 
$$\rho _i \in (l) + I_{\Sigma_i \times \Sigma _{n-i} \ltimes \Bbb{Z}_l^n} $$ It follows from Corollary \ref{reptech2} that there is an $M$ so that 
$$ \rho _i \cdot \rho _{i+1} \cdot \cdots \cdot \rho _{i+M -1} \cdot \rho _{i+M} \in ((l) + I_{\Sigma _n \ltimes \Bbb{Z}_l^n} + ((l)+I_{\Sigma _i \times \Sigma _{n-i} })^{\infty})R[\Sigma _i \times \Sigma _{n-i} \ltimes \Bbb{Z}_l^n ]
$$ 
Referring to the diagram \ref{betterdiagram} above, we have 
$$j^*(\rho _i \cdot \rho _{i+1} \cdot \cdots \cdot \rho _{i+M -1} \cdot \rho _{i+M}) = r_i \cdot r_{i+1} \cdot \cdots \cdot r_{i+M-1} \cdot r_{i+M}$$
and therefore that 

$$r_i\cdot r_{i+1} \cdot \cdots r_{i+M-1} \cdot r_{i+M} \in 
((l) + I_G + ((l) + I_{Q_H})^{\infty})R[H]
$$
But by Proposition \ref{manyproperties}, Part (5), the ideal $((l)+I_{Q_H})^{\infty}$ is equal to the zero ideal, since $Q_H$ an $l$-group. 
It now follows that $ r_i \cdot r_{i+1} \cdot \cdots \cdot r_{i+M-1} \cdot r_{i+M} \in (l) + I_GR[H]$.
\end{Proof} 

\begin{Corollary} The groups $Tor^{R[G]}_t (\pi _k (\frak{M}^G_i({\cal A}_{\eta} \times W   ; \frak{o})/  \frak{M}_{i+1}^G({\cal A}_{\eta} \times W )) , \Bbb{F}_l)$ vanish for all $t$. 
\end{Corollary} 
\begin{Proof} Multiplication by an element of $(l) + I_G$ annihilates $Tor^{R[G]}_i(M, \Bbb{F}_l)$.  
\end{Proof} 

This corollary now yields the following.

\begin{Proposition} The groups $\pi _*(\frak{M}^G_i({\cal A}_{\eta} \times W   ; \frak{o})/  \frak{M}_{i+1}^G({\cal A}_{\eta} \times W ) \column{\wedge}{\frak{K}^G(k)} \Bbb{H}_l)$ are all trivial.  
\end{Proposition}
\begin{Proof} This is an immediate consequence of Proposition \ref{algcriterion}.  
 \end{Proof}

\begin{Corollary} \label{smashversion} The map $\frak{K}^G({\cal A}_{\eta} \times W ) \rightarrow \frak{K}^G({\cal T}_{\eta} \times W)$ induces an equivalence
$$\frak{K}^G(W) \column{\wedge}{K^G(k)} \Bbb{H}_l \rightarrow   \frak{K}^G({\cal A}_{\eta} \times W ) \column{\wedge}{K^G(k)} \Bbb{H}_l \rightarrow \frak{K}^G({\cal T}_{\eta} \times W) \column{\wedge}{K^G(k)} \Bbb{H}_l
$$
\end{Corollary}
\begin{Proof} Immediate from Proposition \ref{newfiltration},  and the decomposition \ref{bardecomp}.
\end{Proof}

\section{The main theorem}

We are now in position to prove our main theorem.  

\begin{Proposition} \label{mainone} Let $G$ be a totally torsion free profinite $l$-group, and let $W$ be any Noetherian approximable continuous affine $G$-scheme.   Then the natural map of $\frak{K}^G(k)$-module spectra $\frak{K}^G(W) \rightarrow \frak{K}^G({\cal E} G \times W) $ induces an equivalence $$ \frak{K}^G(W)\column{\wedge}{K^G(k)} \Bbb{H}_l \rightarrow \frak{K}^G({\cal E} G \times W) \column{\wedge}{K^G(k)} \Bbb{H}_l
$$
\end{Proposition} 
\begin{Proof} Since ${\cal E}G$ is an infinite product of continuous affine $G$-schemes of the form ${\cal T}_{\eta}$, the equivariant $\frak{K}$-theory of ${\cal E}G \times W$ can be expressed as a colimit of spectra of the form $\frak{K}^G({\cal T}_{\eta} \times W)$.  Corollary \ref{smashversion} asserts that the map 
$$ \frak{K}^G(W) \column{\wedge}{K^G(k)} \Bbb{H}_l  \rightarrow \frak{K}^G({\cal T}_{\eta} \times W) \column{\wedge}{K^G(k)} \Bbb{H}_l
$$is an equivalence. 
The construction $-\column{\wedge}{K^G(k)} \Bbb{H}_l$ commutes with colimits, and the result therefore follows.  
\end{Proof} 

\begin{Theorem} \label{maintwo} Let $G$ and $W$  be as above, then there are equivalences of spectra 
$$  \frak{K}^G(W)^{\wedge}_{\epsilon} \stackrel{\sim}{\rightarrow} \frak{K}^G({\cal E}G \times W)^{\wedge}_{\epsilon} \cong \frak{K}({\cal E}G \column{\times}{G} W) ^{\wedge}_l
$$
well defined up to homotopy, and natural in $W$. 
\end{Theorem}
\begin{Proof} Follows directly from  Proposition \ref{algcriterion}, Corollary \ref{nonequiv} and Proposition \ref{orbitequiv}.  
\end{Proof}  
\begin{Corollary}\label{mainfinal} 
Let $G = G_F$ denote the absolute Galois group of a field  $F$ containing $k$, so that $G_F$ is a pro-$l$ group and so that $l$ is prime to the characteristic of $k$. Let $\overline{F}$ denote the algebraic closure of $G$, equipped with the defining $G_F$-action.   Then there are natural equivalences of spectra
$$  \frak{K}^G(k)^{\wedge}_{\epsilon} \stackrel{\sim}{\rightarrow} \frak{K}^G({\cal E}G )^{\wedge}_{\epsilon} \cong \frak{K}({\cal E}G/G ) ^{\wedge}_l
$$ and 
spectra 
$$  \frak{K}(F)^{\wedge}_l \cong \frak{K}^G(\overline{F})^{\wedge}_{\epsilon} \stackrel{\sim}{\rightarrow} \frak{K}^G({\cal E}G \times Spec(\overline{F}))^{\wedge}_{\epsilon} \cong \frak{K}({\cal E}G \column{\times}{G} Spec(\overline{F})) ^{\wedge}_l
$$

\end{Corollary}
\begin{Proof}  We simply apply Theorem \ref{maintwo} in the cases $W = Spec(k)$ and $W= Spec(\overline{F})$.  The left hand equivalence in the case $W = Spec(\overline{F})$ is proved in Corollary \ref{fineequiv}.
\end{Proof}

 \end{document}